\documentclass[a4paper,12pt,leqno]{amsart}

\usepackage{amsmath}
\usepackage{amssymb}
\usepackage{amsthm}
\usepackage{amsopn}
\usepackage{amscd}
\usepackage{xcolor}
\usepackage[utf8]{inputenc}
\usepackage[T1]{fontenc}
\usepackage[colorlinks=true,linkcolor=blue,citecolor=blue,urlcolor=blue]{hyperref}

\theoremstyle{plain}
\newtheorem{thm}{Theorem}[section]
\newtheorem{lemma}[thm]{Lemma}
\newtheorem{cor}[thm]{Corollary}
\newtheorem{prop}[thm]{Proposition}

\theoremstyle{definition}
\newtheorem{defi}[thm]{Definition}

\newtheorem{rmk}[thm]{Remark}

\providecommand{\caJ}{\mathcal{J}}\providecommand{\caL}{\mathcal{L}}\providecommand{\caR}{\mathcal{R}}
\providecommand{\bbZ}{\mathbb{Z}}\providecommand{\bbC}{\mathbb{C}}\providecommand{\Id}{\mathrm{Id}}
\providecommand{\frS}{\mathfrak{S}}\providecommand{\frs}{\mathfrak{s}}
\begin{document}

\numberwithin{equation}{section}

\title{Resolution of ladder representations for classical groups}
\author{David Renard}
\address{
Centre de Mathématiques \\
Laurent Schwartz
\\ Ecole Polytechnique }
\date{\today}

\begin{abstract}
Atobe has introduced ladder representations of the $p$-adic groups $\mathbf{Sp}(2n)$ and
$\mathbf{SO}(2n+1)$, which contain the strongly positive discrete series and M\oe glin's
representations attached to irreducible Arthur parameters, and has expressed their characters as
alternating sums of costandard modules. We show that all the coefficients of this formula are $\pm1$,
with the sign given by the dimension of the orbit of the Langlands parameter, and that these
representations have coresolutions of Cousin type by costandard modules, whose Euler characteristic is
Atobe's formula.
\end{abstract}

\maketitle

\section*{Introduction}

Let $F$ be a non-Archimedean local field of characteristic zero, and let $G_n$ be the split group
$\mathbf{Sp}(2n,F)$ or $\mathbf{SO}(2n+1,F)$. In \cite{Ren26}, we constructed, for a ladder
representation $\pi$ of $\mathbf{GL}_m(F)$, an exact sequence of costandard modules whose Euler
characteristic is the determinantal formula of Tadi\'c and Lapid--M\'inguez \cite{LM14}. Atobe
\cite{Ato24} has introduced a class of ladder representations $\pi(\caL)$ of the groups $G_n$, which
contains the strongly positive discrete series in the sense of M\oe glin (\cite{Moe02}; see also
\cite{MT02}) and the irreducible representations with
irreducible $A$-parameters, and has proved a formula expressing the class of $\pi(\caL)$ in the Grothendieck group as an
alternating sum of the classes of representations $I_\sigma(\caL)$, indexed by a set $\frS(\caL)$ of
permutations (equivalently, an identity between their characters). Each $I_\sigma(\caL)$ is a direct sum
of \emph{costandard modules}, i.e.\ of representations $J=\Delta_1\times\cdots\times\Delta_r\rtimes\tau$
(essentially square integrable $\Delta_i$ of general linear groups with increasing exponents, $\tau$
tempered) as in the subrepresentation form of the Langlands classification, which have a simple socle
$L(J)$ (\S\ref{ssec:std}; Atobe calls them standard modules).

 For $\mathbf{GL}_m(F)$, the Speh
representations are the irreducible representations with an irreducible Arthur parameter
$\rho\boxtimes S_a\boxtimes S_b$, and every ladder representation is obtained from a Speh representation
by derivatives; more generally, the products of Speh representations attached to distinct unitary
cuspidal $\rho$ (Arthur parameters $\bigoplus_\rho\rho\boxtimes S_{a_\rho}\boxtimes S_{b_\rho}$) give in the
same way the products of ladder representations attached to these $\rho$. For $G_n$, Atobe's ladder
representations are exactly the representations obtained by derivatives from the representations of
M\oe glin's construction in the Arthur packets of the parameters
$\bigoplus_\rho\rho\boxtimes S_{a_\rho}\boxtimes S_{b_\rho}$, over distinct self-dual $\rho$ (Proposition
\ref{prop:arthurclosure}). In particular, the class is stable under the derivatives which govern the
induction.

Grouping the isomorphic summands, Atobe's formula reads $[\pi(\caL)]=\sum_Jc_J[J]$, the sum running over
costandard modules $J$ with integer coefficients $c_J$. Let $\Phi(\caL)$ be the set of the $J$ with
$c_J\neq0$; it contains the summand $J_0=I_{\Id}(\caL)$, whose socle is $\pi(\caL)$. The \emph{degree}
$\deg J$ of $J\in\Phi(\caL)$ is the dimension of the orbit of its Langlands parameter in the degree $-1$
part $\widehat{\mathfrak g}_{-1}$ of a graded Lie algebra attached to the infinitesimal character, minus
that of $J_0$ (\S\ref{sec:geom}). Only the differences of these dimensions matter, as for the lengths
in Kazhdan--Lusztig theory; the normalization puts $J_0$ in degree $0$, and then $\deg J\geq0$, with
equality only for $J_0$. We prove (Theorem
\ref{conj:geom}) that
\begin{enumerate}
\item $c_J=(-1)^{\deg J}$ for every $J\in\Phi(\caL)$; in particular, all the coefficients of Atobe's
formula are $\pm1$;
\item $\pi(\caL)$ has a coresolution by costandard modules of Cousin type, whose Euler characteristic is
Atobe's formula:
\[0\to\pi(\caL)\to\caJ^0\to\caJ^1\to\caJ^2\to\cdots,\qquad \caJ^k=\bigoplus_{J\in\Phi(\caL),\ \deg J=k}J.\]
\end{enumerate}

\subsection*{Resolutions of irreducible representations by standard modules}
Theorem \ref{conj:geom} belongs to a family of results in which an irreducible representation is
resolved by standard (or costandard) modules, one for each stratum of a geometric stratification,
placed in the degree given by the codimension of the stratum:
\begin{itemize}
\item the Bernstein--Gelfand--Gelfand resolution of a finite-dimensional representation by Verma
modules \cite{BGG75}, whose dual is the Cousin resolution of the structure sheaf of the flag variety
along the Bruhat cells (Kempf \cite{Kem78}, Mili\v{c}i\'c \cite{Mil}); its $k$-th term is the sum of
the Verma modules indexed by the Weyl group elements of length $k$, the length being the dimension of
the Bruhat cell;
\item Johnson's resolutions of the cohomologically induced modules $A_\mathfrak{q}(\lambda)$ by
standard Harish-Chandra modules \cite{Joh84}, graded by Vogan's length, which is, up to a constant,
the dimension of the corresponding $K$-orbit in the flag variety (for an integral infinitesimal
character; \cite{LV83}), and the
resolutions of $\theta$-twisted Speh representations of \cite{AMR18}, \S10.3, graded by the
$\theta$-length;
\item for $p$-adic $\mathbf{GL}_m$, the image of the dual BGG resolution under the Arakawa--Suzuki
functor \cite{BC15}, and the coresolutions of the ladder representations of \cite{Ren26}, in particular
of the Speh representations, where the terms are indexed by a lower Bruhat interval of the symmetric
group and the degree is the length, which is the codimension of the corresponding orbit in
Zelevinsky's space of graded nilpotent endomorphisms \cite{Zel81};
\item the ladder representations of Atobe for $\mathbf{Sp}(2n,F)$ and $\mathbf{SO}(2n+1,F)$
(Theorem \ref{conj:geom}), graded by the dimensions of the orbits of the Langlands parameters.
\end{itemize}

In all these cases, the mechanism is that of a Cousin resolution. There is a stratified variety, and
the irreducible representation corresponds to a sheaf whose costalks are concentrated in a single
degree along each stratum, the degree given by its codimension, as for the constant sheaf of a
(rationally) smooth variety; the Cousin complex of this sheaf with respect to the stratification then
gives the coresolution, its terms being the costandard objects ($*$-extensions of the local systems on
the strata). Here the stratified variety is the space $\widehat{\mathfrak g}_{-1}$ of the graded Lie
algebra attached to the infinitesimal character, stratified by the orbits of the group $G_0$, which
parametrize the Langlands parameters (\S\ref{sec:geom}); the multiplicities of the irreducible
representations in the standard modules are stalks of intersection cohomology complexes on the orbit
closures (the $p$-adic Kazhdan--Lusztig conjecture, \cite{Sol25}). As for Harish-Chandra modules, the
stabilizers need not be connected, an orbit may carry several local systems, and the extensions
between objects attached to the same orbit involve its equivariant cohomology; this is controlled by
Kato's structure theorem for extension algebras \cite{Kat17}, for the costandard modules of discrete
type.

The present paper and \cite{Ren26} are part of a project which aims at understanding the extent
of this phenomenon: which irreducible representations of real or $p$-adic reductive groups admit such
resolutions, and how the grading is determined by the geometry of the Langlands (or Arthur)
parameters.

\subsection*{Structure of the proof}
The proof combines combinatorics (Atobe's formula and his involution, Tadi\'c's formula for Jacquet
modules), representation theory (derivatives, socles, extensions) and geometry (orbits of Langlands
parameters, intersection cohomology, Kato's theory of extension algebras). We describe it backwards,
starting from the construction of the coresolution.

\emph{Derivatives.} The main tool of the induction, as in Atobe's proof of his formula, is the
derivative. For a self-dual cuspidal representation $\rho$ of $\mathbf{GL}_d(F)$ and $\lambda\in\mathbb
R$, the derivative $r^{\rho,\lambda}(\pi)$ of a representation $\pi$ of $G_n$ is the generalized
$\rho|\cdot|^\lambda$-isotypic part, for the Bernstein centre of $\mathbf{GL}_d(F)$, of the normalized
Jacquet module of $\pi$ with respect to the maximal parabolic subgroup with Levi factor
$\mathbf{GL}_d(F)\times G_{n-d}$ (\S\ref{ssec:deriv}). The functor $r^{\rho,\lambda}$ is exact, and on
the Grothendieck group it is the $\rho|\cdot|^\lambda$-derivative $D_{\rho|\cdot|^\lambda}$ used by Atobe,
the analogue for classical groups of the $\rho$-derivatives of Jantzen and M\'inguez for general linear
groups (\cite{Jan07}, \cite{AM23}). The derivatives detect the representations which are not
supercuspidal, and the derivative of a ladder representation is either $0$ or again a ladder
representation, of a smaller group (\cite{Ato24}, Proposition 3.6); this makes the induction on the
rank possible.

\emph{The construction} (\S\ref{sec:johnson}). Following Johnson \cite{Joh84}, the coresolution is built by successive cokernels
and the  differentials are not explicitly  written down. Recall that its $k$-th term
$\caJ^k$ is the direct sum of the costandard modules $J\in\Phi(\caL)$ of degree $k$, each with multiplicity
one (so $\caJ^0=J_0$). Starting from $Y_0=\pi(\caL)\subset\caJ^0$, one chooses a map
$f_k:Y_k\to\caJ^k$ all of whose components are non-zero, and puts $Y_{k+1}=\mathrm{coker}(f_k)$. Theorem
\ref{conj:geom}(2) then reduces to three statements:
\begin{itemize}
\item[(a)] \emph{the sequence stops}: by a computation in the Grothendieck group
(Lemma \ref{lem:bookkeeping}), this follows from Atobe's formula once its coefficients are known to be
$c_J=(-1)^{\deg J}$, i.e.\ from Theorem \ref{conj:geom}(1);
\item[(b)] \emph{the maps $f_k$ exist}, and are unique up to scalars (Lemma \ref{lem:existence}): this
follows from (a) and from the property (Inj) below, the analogue of \cite{Joh84}, Theorem III.1. Define the
degree of an irreducible representation $L(J')$ with the infinitesimal character of $\pi(\caL)$ as that
of $J'$, and let $\mathcal C_{\geq m}$ be the category of the representations of finite length on which
the Bernstein centre acts by the character of $\pi(\caL)$, all of whose constituents have degree
$\geq m$. Then (Inj) states that every $J\in\Phi(\caL)$ is injective in $\mathcal C_{\geq\deg J}$. These
$J$ are \emph{of discrete type}: their tempered part $\tau$ is square integrable and their segments are
pairwise distinct; part of the proof of (Inj) applies to all the costandard modules of this type;
\item[(c)] \emph{the maps $f_k$ are injective} (Proposition \ref{prop:johnson}).
\end{itemize}

\emph{Injectivity of the $f_k$} (\S\ref{sec:johnson}, last subsections). By induction on the rank, it
suffices to show that the kernel of $f_k$ is killed by every derivative $r^{\rho,\lambda}$, since these
functors are exact and detect the non-zero representations with the cuspidal support of $\pi(\caL)$
(\S\ref{sec:step1}, by Frobenius reciprocity). This criterion excludes the case where $\pi(\caL)$ is
supercuspidal, which is trivial: then $\Phi(\caL)=\{J_0\}$ and $J_0=\pi(\caL)$.
The ladder datum $\caL$ attaches to each $\rho$ a finite
set $X_\rho$ of half-integers (\S\ref{sec:atobe}), and only the $\lambda\in X_\rho$ matter; such a
$\lambda$ is \emph{unpaired} if $\lambda-1\notin X_\rho$, and \emph{paired} otherwise. For an unpaired
$\lambda$, $r^{\rho,\lambda}$ sends the Johnson sequence of $\caL$ to a Johnson sequence of a ladder
datum of smaller rank (Lemma \ref{lem:inj_unpaired}). For a paired $\lambda$, $r^{\rho,\lambda}$ kills
$\pi(\caL)$, and Atobe's proof of his formula provides an involution of the set of the $\sigma$ with
$r^{\rho,\lambda}(I_\sigma(\caL))\neq0$, which changes the sign of $\sigma$ and preserves the derivative
(\S\ref{sec:step1}). We show that the two members of a pair of summands of the $I_\sigma(\caL)$ exchanged by this involution have
consecutive degrees, and $r^{\rho,\lambda}$ of the component of $f_k$
between them is an isomorphism (Lemmas \ref{lem:degpaired}, \ref{lem:degtilde} and
\ref{lem:matchediso}); the other components are controlled by the closure order of the orbits
(Proposition \ref{prop:closure}).

\emph{The signs} (\S\ref{ssec:signs}, combinatorics and representation theory).
Theorem \ref{thm:signs}, $c_J=(-1)^{\deg J}$, is proved by induction on the rank and on the degree,
and answers \cite{Ato24}, Remark 4.6. It rests on the following dichotomy
(Lemma \ref{lem:reach}): for every $J\in\Phi(\caL)$ other than $J_0$,
\begin{itemize}
\item either there is an unpaired $\lambda$ such that neither $\pi(\caL)$ nor $J$ is killed by
$r^{\rho,\lambda}$; then $r^{\rho,\lambda}(J)$ is, up to the factor $\rho|\cdot|^\lambda$, a summand of the
ladder datum $\caL_\lambda$ of smaller rank such that $r^{\rho,\lambda}(\pi(\caL))\simeq
\rho|\cdot|^\lambda\otimes\pi(\caL_\lambda)$, with the same coefficient and the same degree (Lemmas
\ref{lem:coefunpaired} and \ref{lem:degunpaired}), and one concludes by induction on the rank;
\item or there is a paired $\lambda$ such that $J$ is the member of higher degree (the \emph{upper member}) of a pair $\{J,J^\natural\}$
exchanged by Atobe's involution; then $c_{J^\natural}=-c_J$ and $\deg J^\natural=\deg J-1$ (Lemmas
\ref{lem:coefpaired} and \ref{lem:degpaired}), and one concludes by induction on the degree.
\end{itemize} The proof of
Lemma \ref{lem:reach} uses the behaviour of the socles of the summands of the $I_\sigma(\caL)$ under the derivatives (Lemmas
\ref{lem:lowersocle} and \ref{lem:lowersocle0}).

\emph{The degrees} (\S\S\ref{sec:geom} and \ref{sec:degsigns}, geometry and combinatorics). The degree of
a costandard module $J\in\Phi(\caL)$ is, up to the normalization above, the dimension of the orbit of
its Langlands parameter in the graded Lie algebra $\widehat{\mathfrak g}_{-1}$. A formula for the dimensions of the
centralizers in terms of the strings of the nilpotent element (Lemma \ref{lem:dimformula}) allows one to
compute how the degrees change under the derivatives: they are preserved for the unpaired exponents (Lemma
\ref{lem:degunpaired}), and they differ by one between the members of a pair (Lemmas \ref{lem:degpaired},
\ref{lem:degtilde}).

\emph{The derivatives of the costandard modules} (\S\S\ref{sec:step1} and \ref{sec:local},
representation theory). Atobe proves his formula by comparing the derivatives of both sides, using two
lemmas (\cite{Ato24}, Lemmas 4.8 and 4.9), which compute the derivatives of the representations
$I_\sigma(\caL)$ in the Grothendieck group. We prove module versions of these lemmas: the derivative of
each costandard module occurring in $I_\sigma(\caL)$ is computed by Tadi\'c's formula and the results of \cite{Xu17} (Propositions \ref{prop:cat49} and \ref{prop:cat48}), and the derivatives of their socles are controlled
through the explicit Zelevinsky--Aubert duality of \cite{AM23} (Proposition \ref{thm:SAM}, Lemmas
\ref{lem:matchediso}, \ref{lem:smallpaired}).

\emph{The geometric inputs.} Two results come from the geometry of $\widehat{\mathfrak g}_{-1}$.

The first one is the closure order (Proposition \ref{prop:closure}): if $L(J')$ is a constituent of $J$,
then the orbit of the Langlands parameter of $J$ lies in the closure of the orbit of that of $J'$. This
follows from the $p$-adic Kazhdan--Lusztig conjecture, which Solleveld has proved for the classical
groups \cite{Sol25}. It is used in (a), in (c) and in the proof of the signs.

The second one is the property (Inj) (\S\S\ref{ssec:inj} and \ref{ssec:clifford}). Since
$\mathcal C_{\geq m}$ is a Serre subcategory, (Inj) amounts to the vanishing of $\mathrm{Ext}^1(L',J)$
for the irreducible representations $L'$ of degree $\geq\deg J$, and we distinguish two cases. When
$L'\neq L(J)$, the vanishing follows from Kato's structure theorem for extension algebras \cite{Kat17}
(Theorem \ref{thm:injprime}). The hypotheses of Kato's theorem are purity statements, which are provided
by the work of Lusztig and Yun (Lemma \ref{lem:K}). Kato's theorem also requires the centralizers to be
connected; when they are not, we reduce to the connected case by Clifford theory, following
\cite{AMS18} (\S\ref{ssec:clifford}). When $L'=L(J)$, the question concerns the self-extensions of
$L(J)$ by $J$. These come from the unramified deformations of the segments of $J$: we compute them
with the results of \cite{OS13}, and we show that, when $J\in\Phi(\caL)$, none of them is an extension on
which the Bernstein centre acts by the character of $\pi(\caL)$ (Proposition \ref{prop:selfext}, Corollary \ref{cor:selfext}).

\subsection*{On the use of artificial intelligence}
This paper was written with the help of the AI assistant Claude (Anthropic). The starting point is a
draft by the author, which set the problem (a resolution of Atobe's ladder representations by
costandard modules, extending \cite{Ren26}) and tried to transpose the combinatorial approach of
\cite{Ren26}, indexed by a Bruhat interval; this approach turned out not to apply as such. The author
then proposed to follow Johnson's construction by successive cokernels \cite{Joh84} and the Cousin
point of view, and suggested the use of the geometry of the Langlands parameters. The assistant
formulated the grading by orbit dimensions, wrote the proofs (module versions of Atobe's lemmas,
socles, orbit dimension computations, signs of Atobe's coefficients, the injectivity statements via
Kato's theory, and the Johnson construction), searched and checked the references, corrected an
erratum in \cite{Ato24}, Definition 4.1, and checked numerically several intermediate statements. The
author directed the work (strategy, choice of the references, presentation) and revised the text. The
author takes full responsibility for its content.

\section{Preliminaries}\label{sec:prelim}

For $m\geq0$, let $\caR_m$ be the Grothendieck group of the category of smooth representations of
finite length of $\mathbf{GL}_m(F)$, and $\caR=\bigoplus_m\caR_m$. It is a graded Hopf algebra
(\cite{Zel80}): the product $m$ is given by normalized parabolic induction, denoted $\times$, and
the coproduct $m^*$ by the normalized Jacquet functors with respect to the maximal standard
parabolic subgroups. Let $\caR(G)=\bigoplus_n\caR(G_n)$, where $\caR(G_n)$ is the Grothendieck
group of the category of smooth representations of finite length of $G_n$. For a representation
$\tau$ of $\mathbf{GL}_m(F)$ and a representation $\sigma$ of $G_{n_0}$, we write
$\tau\rtimes\sigma$ for the representation of $G_{m+n_0}$ obtained by normalized parabolic
induction from the standard parabolic subgroup with Levi factor $\mathbf{GL}_m(F)\times G_{n_0}$.
This makes $\caR(G)$ an $\caR$-module, and the Jacquet functors with respect to the maximal
standard parabolic subgroups define a map $\mu^*:\caR(G)\rightarrow\caR\otimes\caR(G)$. By
\cite{Tad95b}, Theorem 5.2,
\[\mu^*(\tau\rtimes\sigma)=M^*(\tau)\rtimes\mu^*(\sigma),\qquad M^*=(m\otimes1)\circ(\widetilde{\ }\otimes m^*)\circ s\circ m^*,\]
where $s$ is the flip $x\otimes y\mapsto y\otimes x$ and $\widetilde{\ }$ is the contragredient.

All the irreducible subquotients of an induced representation $\tau\rtimes\sigma$, with $\tau$ and
$\sigma$ irreducible, have the same cuspidal support. For a cuspidal support $\frs$, we denote by
$p_\frs:\caR(G_n)\rightarrow\caR(G_n)$ the projection which keeps the irreducible representations
with cuspidal support $\frs$ and kills the others. It is the map used in \cite{Ato24}, Theorem 4.2.
On a costandard module, $p_\frs$ is either the identity or $0$.

\subsection{Langlands parameters}\label{ssec:types}
Let $q$ be the cardinality of the residue field of $F$, $W_F$ the Weil group of $F$, and $|\cdot|$ its character corresponding, by local class field
theory, to the normalized absolute value of $F^\times$. For $a\geq1$, let $S_a$ be the irreducible
algebraic representation of $\mathrm{SL}_2(\bbC)$ of dimension $a$ (the $(a-1)$-th symmetric power of the
standard representation; it is denoted $\nu_a$ in \cite{Art13}). By the local Langlands correspondence for $\mathbf{GL}_d(F)$ (\cite{HT01},
\cite{Hen00}), we identify the irreducible cuspidal representations $\rho$ of $\mathbf{GL}_d(F)$ with the
irreducible $d$-dimensional representations of $W_F$; the twist $\rho|\cdot|^x$ corresponds to
$\rho\otimes|\cdot|^x$, and $\rho$ is unitary if and only if its image is bounded. A self-dual $\rho$
preserves a non-degenerate bilinear form, unique up to a scalar, which is symmetric or alternating; we
put $\mathrm{sym}(\rho)=1$ or $-1$ accordingly.

Let $\widehat G$ be the dual group of $G_n$, $\mathrm{SO}(2n+1,\bbC)$ for $\mathbf{Sp}(2n,F)$ and
$\mathrm{Sp}(2n,\bbC)$ for $\mathbf{SO}(2n+1,F)$, let $V$ be its standard representation, with its
invariant form $B$, and put $\mathrm{sym}(B)=1$, resp.\ $-1$. An \emph{$L$-parameter} of $G_n$ is a
homomorphism $\phi:W_F\times\mathrm{SL}_2(\bbC)\to\widehat G$, continuous and with semisimple image on
$W_F$ and algebraic on $\mathrm{SL}_2(\bbC)$, up to conjugacy by $\widehat G$. Through $V$, it is a
representation of $W_F\times\mathrm{SL}_2(\bbC)$ preserving $B$, which we write as a sum
$\phi=\bigoplus_i\rho_i\boxtimes S_{a_i}$ with $\rho_i$ irreducible representations of $W_F$. The
parameter $\phi$ is \emph{tempered} if the $\rho_i$ are unitary, and \emph{discrete} if its centralizer
in $\widehat G$ is finite, i.e.\ if the $\rho_i\boxtimes S_{a_i}$ are pairwise distinct and self-dual of
the symmetry $\mathrm{sym}(B)$; the representation $\rho\boxtimes S_a$, for $\rho$ self-dual, is self-dual
of symmetry $\mathrm{sym}(\rho)(-1)^{a-1}$. The component group of $\phi$ is
$\pi_0(Z_{\widehat G}(\phi))$. By Arthur (\cite{Art13}, Theorem 1.5.1), the irreducible tempered representations of $G_n$
are the $\pi(\phi,\varepsilon)$, where $\phi$ is a tempered $L$-parameter and $\varepsilon$ an irreducible
character of its component group trivial on the image of the centre of $\widehat G$ (the other characters
correspond to the non-split inner form when $G_n=\mathbf{SO}(2n+1,F)$), $\pi(\phi,\varepsilon)$ is square integrable if and only if $\phi$ is
discrete, and every tempered $\phi$ gives a finite set of tempered representations (its $L$-packet).

We say that a self-dual $\rho$ is of \emph{integral type} if $\mathrm{sym}(\rho)=\mathrm{sym}(B)$ and of
\emph{half-integral type} otherwise: the summands $\rho\boxtimes S_{2x+1}$ of discrete parameters have
$x\in\bbZ$ in the first case and $x\in\frac12+\bbZ$ in the second. We put $E_\rho=\bbZ$ in the first case and $E_\rho=\frac12+\bbZ$
in the second: the exponents $x$ of the twists $\rho|\cdot|^x$ which occur below all lie in $E_\rho$. For instance, the trivial character of
$F^\times$ is of integral type for $\mathbf{Sp}(2n,F)$ and of half-integral type for
$\mathbf{SO}(2n+1,F)$.

\begin{rmk}[Comparison with \cite{Ato24}]\label{rmk:good} Atobe formulates his definitions with the
\emph{good} cuspidal representations (\cite{Ato24}, Definition 3.1): $\rho'$ is good if
$\rho'^\vee\simeq\rho'|\cdot|^a$ for some $a\in\bbZ$ and $\rho'|\cdot|^m\rtimes\mathbf 1$ is reducible for
some $m\in\bbZ$, where $\mathbf 1$ is the trivial representation of the trivial group $G_0$. For $\rho$
self-dual, $\rho|\cdot|^x$ is good if and only if $x\in E_\rho$, so that the two formulations agree.
Indeed, $\rho|\cdot|^s\rtimes\mathbf 1$ is induced from the Siegel parabolic subgroup of $G_d$, and by
Shahidi (\cite{Sha90}; see \cite{MT02}, Remark 2.3(i)) it is reducible for exactly one $s_\rho\geq0$,
and $s_\rho\in\{0,\frac12,1\}$. Moreover, $s_\rho=\frac12$ if and only if the $L$-function
$L(s,\rho,R)$ has a pole at $s=0$, where $R$ is the exterior square for $\mathbf{Sp}(2n,F)$ and the
symmetric square for $\mathbf{SO}(2n+1,F)$ (\cite{MT02}, (J-1) and Remark 2.3(i)). Since the local
Langlands correspondence for $\mathbf{GL}_d(F)$ preserves these $L$-functions (\cite{Hen10}), this pole
exists if and only if $\rho$ is symplectic, resp.\ orthogonal, that is, if and only if $\rho$ is of
half-integral type. Hence $s_\rho\in E_\rho$, and $\rho|\cdot|^x$ is good if and only if
$x-s_\rho\in\bbZ$, i.e.\ $x\in E_\rho$.
\end{rmk}

A tempered parameter $\phi$ is \emph{of good parity} if all its summands are of the form
$\rho\boxtimes S_a$ with $\rho$ self-dual (hence unitary) and $\frac{a-1}2\in E_\rho$, i.e.\ with
$\rho\boxtimes S_a$ self-dual of symmetry $\mathrm{sym}(B)$; following M\oe glin, these summands are
called its \emph{Jordan blocks}. All the tempered representations met in this paper have parameters of
good parity. Then the component group of $\phi$ is an elementary abelian $2$-group, and $\varepsilon$ is
described by signs $\varepsilon(\rho\boxtimes S_a)\in\{\pm1\}$ attached to the Jordan blocks, equal for
equal blocks and with product $1$ (\cite{Ato24}, \S2.4). We use the convention of \cite{Ato24}, \S2.4:
when $\rho$ is of half-integral type, one allows formally a block $\rho\boxtimes S_0$ ($S_0$ being the
zero representation, which corresponds to $x=-\frac12$ in $S_{2x+1}$), and sets
$\pi(\phi,\varepsilon)=\pi(\phi-\rho\boxtimes S_0,\varepsilon)$ if $\varepsilon(\rho\boxtimes S_0)=1$ and
$\pi(\phi,\varepsilon)=0$ otherwise. This convention makes the formulas uniform at the bottom of the
chains of Jordan blocks. The derivative $D_{\rho|\cdot|^x}$ (\S\ref{ssec:deriv}) replaces a block
$\rho\boxtimes S_{2x+1}$ by $\rho\boxtimes S_{2x-1}$, and for $x=\frac12$ the block $\rho\boxtimes S_2$
becomes $\rho\boxtimes S_0$: the derivative is then the representation without this block, or $0$,
according to its sign, which is what the convention expresses. Similarly, the ladder data of
\S\ref{sec:atobe} may contain the exponent $-\frac12$, and the signs of the Jordan blocks in Atobe's
formula alternate along the chains, starting, when this exponent occurs, with the formal block
$\rho\boxtimes S_0$.

\subsection{Standard and costandard modules}\label{ssec:std}
Let $\rho$ be an irreducible cuspidal representation of $\mathbf{GL}_d(F)$. For real numbers $x\geq y$
with $x-y\in\mathbb{Z}$, let $\Delta_\rho[x,y]$ be the unique irreducible subrepresentation of
$\rho|\cdot|^x\times\rho|\cdot|^{x-1}\times\cdots\times\rho|\cdot|^y$; it is essentially square
integrable, with central exponent $(x+y)/2$ (\cite{Ato24}, \S2.4; in the notation of \cite{Zel80} it is
attached to the segment $\{\rho|\cdot|^y,\ldots,\rho|\cdot|^x\}$). We also put $\Delta_\rho[x,x+1]=\mathbf{1}$
(the trivial representation of $\mathbf{GL}_0$). When $\rho$ is fixed, we write $\Delta[x,y]$ for $\Delta_\rho[x,y]$.

By the Langlands classification, in the subrepresentation form used by Atobe (\cite{Ato24}, \S2.4),
every irreducible representation of $G_n$ is the unique irreducible subrepresentation
$L(\Delta_1,\ldots,\Delta_r;\tau)$ of a representation
\[J=\Delta_{\rho_1}[x_1,y_1]\times\cdots\times\Delta_{\rho_r}[x_r,y_r]\rtimes\tau,\qquad
x_1+y_1\leq\cdots\leq x_r+y_r<0,\]
with $\tau$ irreducible and tempered, and the data $(\Delta_i;\tau)$ are unique up to the order of
the segments with equal exponents. We call such a $J$ a \emph{costandard module}: it has a simple
socle $L(J)=L(\Delta_1,\ldots,\Delta_r;\tau)$, which occurs in $J$ with multiplicity one. (Atobe calls
it a standard module. The contragredient of $J$ is a standard module in the quotient form of the
Langlands classification, with simple head.) The costandard modules form a basis of $\caR(G)$. The
tempered representations are written $\tau=\pi(\phi,\varepsilon)$ as in \S\ref{ssec:types}. We say
that a costandard module $J$ is \emph{of discrete type} if its segments $\Delta_1,\ldots,\Delta_r$ are
pairwise distinct and its tempered part $\tau$ is square integrable; in the parametrization of
\S\ref{sec:geom}, this means that its Langlands parameter is multiplicity free.

\subsection{Jacquet modules and derivatives}\label{ssec:deriv}
For $1\leq d\leq n$, let $P_d$ be the standard maximal parabolic subgroup of $G_n$ with Levi factor
$\mathbf{GL}_d(F)\times G_{n-d}$, and $\mathrm{Jac}_{P_d}$ the normalized Jacquet functor. Let $\rho$ be
an irreducible cuspidal representation of $\mathbf{GL}_d(F)$. Following \cite{Ato24}, \S2.3, the
\emph{$\rho$-derivative} of a representation $\pi$ of $G_n$ of finite length is the semisimple
representation $D_\rho(\pi)$ of $G_{n-d}$ such that
\[[\mathrm{Jac}_{P_d}(\pi)]=[\rho]\otimes[D_\rho(\pi)]+\sum_i[\tau_i]\otimes[\pi_i],\qquad
\tau_i\text{ irreducible},\ \tau_i\not\simeq\rho,\]
and $D_\rho$ is extended by linearity to $\caR(G)$. It is the analogue of the $\rho$-derivatives of
Jantzen and M\'inguez for the general linear groups, which refine the Bernstein--Zelevinsky
derivatives: instead of the whole derivative, one keeps the part of the Jacquet module on which the
factor $\mathbf{GL}_d(F)$ acts by the given cuspidal $\rho$ (\cite{Jan07}). Iterating, the
\emph{highest $\rho$-derivative} of an irreducible $\pi$ is $D_\rho^{k}(\pi)$, with $k$ maximal such
that $D_\rho^{k}(\pi)\neq0$ (\cite{AM23}, \S3.1). We shall use $D_\rho$ for $\rho=\rho_0|\cdot|^\lambda$, $\rho_0$ self-dual.

We also need a module-level version. The Jacquet module $\mathrm{Jac}_{P_d}(\pi)$ is a
representation of finite length of $\mathbf{GL}_d(F)\times G_{n-d}$, on which the Bernstein centre
$\mathfrak Z_d$ of $\mathbf{GL}_d(F)$ acts. Let $\chi_\rho$ be the character of $\mathfrak Z_d$ by which it
acts on $\rho$, and let $r^\rho(\pi)$ be the generalized eigenspace
\[r^\rho(\pi)=\{v\in\mathrm{Jac}_{P_d}(\pi):(z-\chi_\rho(z))^Nv=0\text{ for all }z\in\mathfrak
Z_d\text{ and }N\gg0\}.\]
It is a direct summand of $\mathrm{Jac}_{P_d}(\pi)$, stable under $\mathbf{GL}_d(F)\times G_{n-d}$: the
largest subrepresentation all of whose irreducible subquotients are of the form $\rho\otimes\tau$. As a
representation of $\mathbf{GL}_d(F)$, it is a successive extension of copies of $\rho$, which need not
be semisimple (the self-extensions of $\rho$ come from its unramified twists). Since $\rho$ is
cuspidal, its infinitesimal character determines it among the irreducible representations of
$\mathbf{GL}_d(F)$, so that $r^\rho$ is an exact functor, and $[r^\rho(\pi)]=[\rho]\otimes[D_\rho(\pi)]$. For $\rho_0$ self-dual and $\lambda\in\mathbb R$ we put
$r^{\rho_0,\lambda}=r^{\rho_0|\cdot|^\lambda}$, and we write $r^\lambda$ and $D_\lambda$ when $\rho_0$ is fixed.

\section{Ladder representations of Atobe}\label{sec:atobe}

\subsection{Ladder data}
We recall \cite{Ato24}, Definition 3.1, in the formulation of \S\ref{ssec:types} (Remark
\ref{rmk:good}). A \emph{ladder datum} for $G_n$ is a triple $\caL=(X,l,\eta)$,
$X=(X_\rho)_\rho$, $l=(l_\rho)_\rho$, $\eta=(\eta_\rho)_\rho$, indexed by the self-dual cuspidal
representations $\rho$ of general linear groups, such that:
\begin{enumerate}
\item $X_\rho=\{\rho|\cdot|^{x_1},\ldots,\rho|\cdot|^{x_{t_\rho}}\}$ with $x_1<\cdots<x_{t_\rho}$ in
$E_\rho$, $t_\rho=0$ for almost all $\rho$, and
$\sum_\rho(2x_1+\cdots+2x_{t_\rho}+t_\rho)\dim\rho$ equal to $2n$ for $\mathbf{SO}(2n+1,F)$ and to
$2n+1$ for $\mathbf{Sp}(2n,F)$;
\item $0\leq2l_\rho\leq t_\rho$, $x_j+x_{t_\rho-j+1}\geq0$ for $1\leq j\leq l_\rho$, and
$x_i\geq-\tfrac12$ for $l_\rho<i\leq t_\rho-l_\rho$;
\item $\eta_\rho\in\{\pm1\}$, with $\eta_\rho=1$ if $X_\rho=\emptyset$ or $x_{l_\rho+1}=-\tfrac12$,
$\eta_\rho=-1$ if $2l_\rho=t_\rho$, and $\prod_\rho(-1)^{[t_\rho/2]+l_\rho}\eta_\rho^{t_\rho}=1$.
\end{enumerate}
The \emph{ladder representation} $\pi(\caL)$ is the unique irreducible subrepresentation of the
costandard module
\[I(\caL)=\mathop{\times}_\rho\big(\Delta_\rho[x_1,-x_{t_\rho}]\times\Delta_\rho[x_2,-x_{t_\rho-1}]
\times\cdots\times\Delta_\rho[x_{l_\rho},-x_{t_\rho-l_\rho+1}]\big)\rtimes\pi(\phi,\varepsilon),\]
where $\phi=\bigoplus_\rho\bigoplus_{i=l_\rho+1}^{t_\rho-l_\rho}\rho\boxtimes S_{2x_i+1}$ and
$\varepsilon(\rho\boxtimes S_{2x_i+1})=(-1)^{i-l_\rho-1}\eta_\rho$. The map $\caL\mapsto\pi(\caL)$ is
injective. If $l_\rho=0$ for all $\rho$, $\pi(\caL)$ is a strongly positive discrete series, and every
strongly positive discrete series is obtained in this way; if each $X_\rho$ is a segment, $\pi(\caL)$ is
one of M\oe glin's representations with an $A$-parameter (\cite{Ato24}, Remark 3.2; see
\S\ref{ssec:arthur}).

When a self-dual $\rho$ is fixed, we drop it from the notation: $X=\{x_1<\cdots<x_t\}$, $l$, $\eta$.
For $\mathbf{Sp}(2n,F)$ and $\rho=\mathbf 1$ only, this is \cite{Ato24}, Definition 1.2.

\subsection{Ladder representations and Arthur packets}\label{ssec:arthur}
An $A$-parameter (Arthur parameter) of $G_n$
is a homomorphism $\psi:W_F\times\mathrm{SL}_2(\bbC)\times\mathrm{SL}_2(\bbC)\to\widehat G$, bounded on
$W_F$, the second $\mathrm{SL}_2(\bbC)$ being Arthur's; it determines the $L$-parameter
\[\phi_\psi(w,g)=\psi\big(w,g,\mathrm{diag}(|w|^{1/2},|w|^{-1/2})\big),\]
and Arthur attaches to it a finite
packet $\Pi_\psi$ of unitary representations, containing the $L$-packet of $\phi_\psi$ (\cite{Art13}). For
$\mathbf{GL}_m(F)$, the packets are singletons, and the irreducible $A$-parameters
$\rho\boxtimes S_a\boxtimes S_b$ give the Speh representations: by M\oe glin and Waldspurger, they are
the Langlands quotients attached to $\phi_\psi=\bigoplus_{k=0}^{b-1}\rho|\cdot|^{(b-1)/2-k}\boxtimes
S_a$, i.e.\ the ladder representations whose $b$ segments $\Delta_\rho[\frac{a-1}2+c,-\frac{a-1}2+c]$,
$c=\frac{b-1}2,\ldots,-\frac{b-1}2$, are all shifted by one; and every ladder representation of
$\mathbf{GL}_m(F)$ is obtained from a Speh representation by derivatives (\cite{Ato24},
Introduction). For $G_n$, the analogous
parameters are the $\psi=\rho\boxtimes S_a\boxtimes S_b$ with $\rho$ self-dual and
$\mathrm{sym}(\rho)(-1)^{a+b}=\mathrm{sym}(B)$ (notation of \S\ref{ssec:types}), or more generally the
sums $\psi=\bigoplus_\rho\rho\boxtimes S_{a_\rho}\boxtimes S_{b_\rho}$ over distinct $\rho$ of this kind. We
attach to such a $\psi$ the segments
\[X_\rho=\{\rho|\cdot|^{x}:\tfrac{a_\rho-b_\rho}2\leq x\leq\tfrac{a_\rho+b_\rho}2-1\},\]
which determine it ($t_\rho=b_\rho$, and $a_\rho=x_1+x_{t_\rho}+1$).

\emph{The packet $\Pi_\psi$.}  Since the $\rho$ are distinct, $\psi$ has
discrete diagonal restriction, and $\Pi_\psi$ consists of the representations $\pi_M(\psi,l,\eta)$ of
M\oe glin's construction, indexed by pairs $(l,\eta)=(l_\rho,\eta_\rho)_\rho$ (M\oe glin's theorem, as
stated in \cite{Xu17b}, Theorem 7.8, with the non-vanishing criterion of Proposition 7.10 there). These
are exactly Atobe's $\pi(\caL)$ for the ladder data $\caL$ with the above segments $X_\rho$, the pairs
$(l_\rho,\eta_\rho)$ of $\caL$ corresponding to M\oe glin's parameters (\cite{Ato24}, Remark 3.2(2); for
$a_\rho\geq b_\rho$, $l_\rho$ runs over $[0,[b_\rho/2]]$ in both parametrizations). Hence:
\begin{center}
\emph{the members of the packets $\Pi_\psi$, $\psi=\bigoplus_\rho\rho\boxtimes S_{a_\rho}\boxtimes
S_{b_\rho}$, are exactly the ladder representations $\pi(\caL)$ all of whose $X_\rho$ are segments.}
\end{center}

\emph{Two families inside $\Pi_\psi$.} The packet contains the $L$-packet of
\[\phi_\psi=\bigoplus_\rho\bigoplus_{k=0}^{b_\rho-1}\rho|\cdot|^{(b_\rho-1)/2-k}\boxtimes S_{a_\rho}\]
(\cite{Art13}), which the Langlands classification describes as follows. The summands
$\rho|\cdot|^{c}\boxtimes S_{a_\rho}$ with $c<0$ give the segments $\Delta_\rho[c+\frac{a_\rho-1}2,c-\frac{a_\rho-1}2]$,
those with $c>0$ their duals, and the summand with $c=0$, present if $b_\rho$ is odd, gives the tempered
parameter $\phi_0=\bigoplus_\rho\rho\boxtimes S_{a_\rho}$ ($b_\rho$ odd). The $L$-packet of $\phi_\psi$ consists of the
unique irreducible subrepresentations of the costandard modules
\[\mathop{\times}_\rho\mathop{\times}_{c<0}\Delta_\rho[c+\tfrac{a_\rho-1}2,c-\tfrac{a_\rho-1}2]\rtimes\pi(\phi_0,\varepsilon),\]
where
$\varepsilon$ runs over the characters of the component group of $\phi_0$ (\S\ref{ssec:types}). These are
the members of $\Pi_\psi$ with $l_\rho=[\frac{b_\rho}2]$ for all $\rho$: the segment
$\Delta_\rho[x_j,-x_{t_\rho-j+1}]$, $j\leq l_\rho$, is the one of exponent $c=-\frac{b_\rho-1}2+j-1$, and
$\varepsilon$ is given by the signs $\eta_\rho$. At the other extreme, when $a_\rho\geq b_\rho$ for all
$\rho$, the members with $l_\rho=0$ for all $\rho$ are the discrete series in $\Pi_\psi$, whose Jordan
blocks $\rho\boxtimes S_{a_\rho-b_\rho+1},\ldots,\rho\boxtimes S_{a_\rho+b_\rho-1}$ carry alternating
signs. For instance, $b=1$ gives the discrete series $\pi(\rho\boxtimes S_a,\varepsilon)$, and
$\psi=\mathbf 1\boxtimes S_1\boxtimes S_{2n+1}$ gives the trivial representation of $\mathbf{Sp}(2n,F)$, with
$X=\{-n,\ldots,n\}$ and $l=n$.

\emph{Derivatives.} The class of the members of the packets $\Pi_\psi$, $\psi$ as above, is not stable
under the derivatives, which govern the induction in \cite{Ato24} (\cite{Ato24}, Introduction). The
ladder representations are the representations obtained from this class by derivatives, as the
products of ladder representations attached to distinct unitary cuspidal $\rho$ of $\mathbf{GL}_m(F)$ are
obtained from the products of Speh representations attached to these $\rho$. (For $\mathbf{GL}_m(F)$, the
latter statement reduces to the case of one $\rho$. A product $\pi_1\times\cdots\times\pi_k$ of irreducible
representations whose cuspidal supports lie on distinct lines $\rho_i|\cdot|^{\mathbb Z}$ is irreducible
(\cite{Zel80}). By the geometric lemma, its $\rho_i|\cdot|^x$-derivative is $\pi_1\times\cdots\times
D_{\rho_i|\cdot|^x}(\pi_i)\times\cdots\times\pi_k$, since the other factors have no cuspidal support on
the line of $\rho_i$; so the derivatives act on each factor separately.)

\begin{prop}\label{prop:arthurclosure} The ladder representations of Atobe form the smallest class of
irreducible representations of the groups $G_n$ which contains the representations $\pi(\caL)$ with
$X_\rho$ a segment for every $\rho$ (that is, the representations of M\oe glin's construction in the
packets $\Pi_\psi$, $\psi=\bigoplus_\rho\rho\boxtimes S_{a_\rho}\boxtimes S_{b_\rho}$) and which is stable
under the non-zero derivatives $D_{\rho|\cdot|^x}$, $\rho$ self-dual and $x\in E_\rho$.
\end{prop}

\begin{proof}
By \cite{Ato24}, Proposition 3.6, a non-zero derivative $D_{\rho|\cdot|^x}(\pi(\caL))$ is irreducible and
equal to $\pi(\caL')$, where $\caL'$ is obtained by replacing an element $x=x_j$ of $X_\rho$, with
$x-1\notin X_\rho$, by $x-1$ (see Lemma \ref{lem:ato49} below); so the class of ladder representations
is stable under the derivatives. Conversely, let $\caL'=(X',l,\eta)$ be a ladder datum. Since
$D_{\rho|\cdot|^x}$ only modifies the component $\caL'_\rho$, we may treat one $\rho$ at a time. Write
$X'_\rho=\{x'_1<\cdots<x'_t\}$ (exponents), and let $\caL$ be the datum with the same $l$ and $\eta$ and
with the segment $X_\rho=\{x'_t-t+1,\ldots,x'_t\}$. Then $x_j\geq x'_j$ for all $j$. The conditions of
\S\ref{sec:atobe} on the exponents are lower bounds, and those on $\eta$ only involve $t$, $l$ and the
equality $x_{l+1}=-\frac12$, which forces $x'_{l+1}=-\frac12$; so $\caL$ is a ladder datum. We lower
$x_1$ to $x'_1$, one unit at a time, then $x_2$ to $x'_2$, and so on. At each step, the current datum $\caL''=(X'',l,\eta)$, with
$X''_\rho=\{x''_1<\cdots<x''_t\}$, is modified by replacing an element $x$ of $X''_\rho$ by $x-1$, with
$x-1\notin X''_\rho$, and the new datum is again a ladder datum, since its exponents lie between those of $\caL'$ and $\caL$. It remains to see that the
corresponding derivative $D_{\rho|\cdot|^x}$ is non-zero. By \cite{Ato24}, Proposition 3.6, this holds
if the graph $\mathcal E(\caL''_\rho)$ of \cite{Ato24}, Definition 3.3, has a minimal vertex $(x,y)$
with $f_\rho(x,y)=0$. Let $x=x''_{l+1+i}$; the row $-i$ of the graph runs from $x$ down to
$-x''_{t-l-i}$. Its top vertex $(x,-i)$ is minimal: the only arrow which could end there is the
diagonal arrow from $(x-1,-(i-1))$, and this vertex exists only if $x-1\leq x''_{l+i}$, i.e.\ $x-1\in
X''_\rho$. Moreover $f_\rho(x,-i)=0$: this is clear if $i<0$ or $i\geq t-2l$, and otherwise
$f_\rho(x,-i)\neq0$ would mean $x\leq i$ (resp.\ $x\leq i-\frac\eta2$ if $\rho$ is of half-integral
type), whereas the new datum, being a ladder datum, satisfies $x-1\geq i$ (resp.\ $x-1\geq
i-\frac\eta2$) (\cite{Ato24}, Definition 3.3).
\end{proof}

\subsection{Atobe's formula}
We recall \cite{Ato24}, Definition 4.1 and Theorem 4.2. We write the data for one self-dual $\rho$,
dropped from the notation; in general, $\frS(\caL)=\prod_\rho\frS(\caL_\rho)$, where
$\caL_\rho=(X_\rho,l_\rho,\eta_\rho)$, the sign is the product of the signs, and the segments, the Jordan
blocks and their signs of $I_\sigma(\caL)$, $\sigma=(\sigma_\rho)_\rho$, are the unions of those attached to
the $\sigma_\rho$, the maps $\delta$ being defined on the disjoint union of the sets $J^-_{\sigma_\rho}$.
\begin{enumerate}
\item $\frS(\caL)$ is the set of permutations $\sigma$ of $\{1,\ldots,t\}$ such that $\sigma$ is
increasing on $\{1,\ldots,l\}$ and on $\{l+1,\ldots,t-l\}$, such that $\sigma^{-1}(i)\leq l$ whenever
$x_i\leq-1$, and such that $\sigma^{-1}(i)\notin\{l+1,\ldots,t-l\}$ whenever $x_i=-\tfrac12$ and
$\eta=-1$ (see Remark \ref{rmk:typo} for the last condition).
\item $\mathrm{sgn}:\frS(\caL)\to\{\pm1\}$ is the restriction of the sign character.
\item For $\sigma\in\frS(\caL)$, let 
\[ J^+_\sigma=\{1\leq j\leq l:\sigma(j)<\sigma(t-j+1)\}, \quad  J^-_\sigma=\{1\leq j\leq l:\sigma(j)>\sigma(t-j+1)\},
\]
 and
\[I_\sigma(\caL)=\mathop{\times}_{j\in J^+_\sigma}\Delta[x_{\sigma(j)},-x_{\sigma(t-j+1)}]\rtimes
\Big(\bigoplus_{\delta}\pi(\phi_\sigma,\varepsilon_{\sigma,\delta})\Big),\]
where $\delta$ runs over the maps $J^-_\sigma\to\{\pm1\}$,
\[\phi_\sigma=\bigoplus_{i=l+1}^{t-l}S_{2x_{\sigma(i)}+1}\oplus\bigoplus_{j\in J^-_\sigma}
\big(S_{2x_{\sigma(j)}+1}\oplus S_{2x_{\sigma(t-j+1)}+1}\big),\]
and
\[\varepsilon_{\sigma,\delta}(S_{2x_{\sigma(i)}+1})=(-1)^{i-l-1}\eta\ \ (l<i\leq t-l),\]
\[\varepsilon_{\sigma,\delta}(S_{2x_{\sigma(j)}+1})=\varepsilon_{\sigma,\delta}(S_{2x_{\sigma(t-j+1)}+1})
=\delta(j)\ \ (j\in J^-_\sigma).\] (A segment $\Delta[x,x+1]$ is $\mathbf{1}$, a segment $\Delta[x,y]$
with $y>x+1$ makes the summand $0$, and the convention on $S_0$ of \S\ref{ssec:types} applies.)
\end{enumerate}
In (3), each $j\leq l$ pairs the exponents $x_{\sigma(j)}$ and $x_{\sigma(t-j+1)}$. For $j\in J^+_\sigma$
they are the ends of the segment $\Delta[x_{\sigma(j)},-x_{\sigma(t-j+1)}]$. For $j\in J^-_\sigma$ they are
in the opposite order, and they give instead two Jordan blocks $\rho\boxtimes S_{2x_{\sigma(j)}+1}$ and
$\rho\boxtimes S_{2x_{\sigma(t-j+1)}+1}$ of the tempered part, which carry the same sign; $\delta(j)\in\{\pm1\}$
is this common sign. Thus $\delta$ runs over the set $\{\pm1\}^{J^-_\sigma}$ of sign choices, and
$I_\sigma(\caL)$ is the direct sum of the $2^{|J^-_\sigma|}$ representations attached to them (some of
which may be $0$).
Each $I_\sigma(\caL)$ is a direct sum of costandard modules, some of which may be $0$ or have the wrong
cuspidal support.

\begin{thm}[\cite{Ato24}, Theorem 4.2]\label{thm:atobe} Let $\frs$ be the cuspidal support of
$\pi(\caL)$. Then, in $\caR(G)$,
\[[\pi(\caL)]=p_\frs\Big(\sum_{\sigma\in\frS(\caL)}\mathrm{sgn}(\sigma)\,[I_\sigma(\caL)]\Big).\]
\end{thm}

\begin{rmk}\label{rmk:typo} As printed, the last condition of \cite{Ato24}, Definition 4.1(1), reads:
if $x_i=-\tfrac12$ and $\eta=-1$, then $\sigma^{-1}(i)\leq l$. This is not compatible with
\cite{Ato24}, Example 4.3: for the trivial representation of $\mathbf{SO}(5,F)$, $X=\{-\tfrac32,
-\tfrac12,\tfrac12,\tfrac32\}$, $l=2$, $\eta=-1$. The printed conditions force $\sigma^{-1}(1),
\sigma^{-1}(2)\leq2$, hence $\sigma(1)=1$, $\sigma(2)=2$ ($\sigma$ is increasing on $\{1,2\}$), and give
two permutations, whereas Example 4.3 uses the six permutations with $\sigma(1)=1$, and Theorem 4.2
needs all of them. These six are exactly the permutations allowed by the condition in (1) above:
for $x_i=-\tfrac12$ and $\eta=-1$, $i$ must not be sent to the range $\{l+1,\ldots,t-l\}$ of the
Jordan blocks of $\phi$, where it would give the block $S_0$ with a sign which may be $-1$ (here this
range is empty). This is the reading we use. Similarly, in the proof of \cite{Ato24}, Lemma 4.8, case
(2), the segment of $I_\sigma(\caL)$ printed as $\Delta[\lambda,-(\lambda-1)]$ is
$\Delta[x_{j-1},-x_j]=\Delta[\lambda-1,-\lambda]$ (Lemma \ref{lem:ato48}), and case (1) is written when
$\lambda-1$ and $\lambda$ are the tops of two segments; Proposition \ref{prop:cat48} gives the module version
of the lemma in all cases.
\end{rmk}

\noindent The right-hand side of Theorem \ref{thm:atobe} may contain cancelling summands
(\cite{Ato24}, Examples 4.4 and 4.5), and Atobe expects that, after cancellation, the coefficients of
$[\pi(\caL)]$ in the basis of costandard modules are $\pm1$ (\cite{Ato24}, Remark 4.6). Section
\ref{sec:conj} explains where the signs and the cancellations come from, and Theorem \ref{thm:signs}
proves Atobe's expectation.

\section{The geometric grading and the main theorem}\label{sec:conj}

\subsection{The summands of Atobe's formula}\label{ssec:summands}
For $\sigma\in\frS(\caL)$ and a sign choice $\delta\in\{\pm1\}^{J^-_\sigma}$ such that
$\pi(\phi_\sigma,\varepsilon_{\sigma,\delta})\neq0$, the representation $\times_{j\in
J^+_\sigma}\Delta[x_{\sigma(j)},-x_{\sigma(t-j+1)}]\rtimes\pi(\phi_\sigma,\varepsilon_{\sigma,\delta})$ is a
summand of $I_\sigma(\caL)$, and it is a costandard module (\cite{Ato24}, \S4.1); we keep those with cuspidal 
support $\frs$. Grouping the isomorphic ones, Theorem \ref{thm:atobe} reads
\[[\pi(\caL)]=\sum_Jc_J[J],\qquad c_J=\sum_{(\sigma,\delta)}\mathrm{sgn}(\sigma),\]
the second sum running over the pairs $(\sigma,\delta)$ whose summand is isomorphic to $J$. Let
$\Phi(\caL)$ be the set of costandard modules $J$ with $c_J\neq0$; since the costandard modules are a
basis of $\caR(G)$, it only depends on $\pi(\caL)$. The module $J_0=I(\caL)$ ($\sigma=\Id$) is in
$\Phi(\caL)$, with $c_{J_0}=1$, and $\pi(\caL)=L(J_0)$.

Different pairs $(\sigma,\delta)$ may give isomorphic summands, with equal or opposite signs, so that
$c_J$ is in general a sum of several signs. These summands then cancel partially, or completely, in
which case $J\notin\Phi(\caL)$ and $J$ does not occur in the coresolution of Theorem \ref{conj:geom}.
Remark \ref{rem:twins} gives an example where three summands are isomorphic, with coefficient
$2-1=1$. By Theorem \ref{conj:geom}(1), $c_J=\pm1$ for every $J\in\Phi(\caL)$.

We use the following notation. We write $J_{\sigma,\delta}=J_{\sigma,\delta}(\caL)$ for the summand of $I_\sigma(\caL)$ attached to the sign choice $\delta\in\{\pm1\}^{J^-_\sigma}$; we call the non-zero ones with cuspidal support $\frs$ the \emph{summands} of $\caL$ (or of
Atobe's formula). When $\rho$ is fixed, the
positions $1,\ldots,t$ are divided into the three intervals $\mathrm L=\{1,\ldots,l\}$,
$\mathrm M=\{l+1,\ldots,t-l\}$ and $\mathrm R=\{t-l+1,\ldots,t\}$: for $\sigma\in\frS(\caL)$, the
positions $j\in\mathrm L$ and $t-j+1\in\mathrm R$ carry the two ends of the $j$-th segment (or the two
members of a pair of Jordan blocks if $j\in J^-_\sigma$), and the positions of $\mathrm M$ carry the
middle Jordan blocks. More precisely, $J_{\sigma,\delta}$ is determined by $\delta$ and by the sequence
$(x_{\sigma(1)},\ldots,x_{\sigma(t)})$ of the exponents placed at the positions $1,\ldots,t$: for $j\in
J^+_\sigma$, the positions $j$ and $t-j+1$ give the segment $\Delta[x_{\sigma(j)},-x_{\sigma(t-j+1)}]$, and every
other position $p$ gives the Jordan block $\rho\boxtimes S_{2x_{\sigma(p)}+1}$, with the sign attached to $p$
(given by $\eta$ if $p\in\mathrm M$, by $\delta$ otherwise). Accordingly, each exponent $x_i$ is exactly
one of the following in $J_{\sigma,\delta}$: a \emph{top}, if $\sigma^{-1}(i)\in J^+_\sigma$ ($x_i$ is the top
of a segment); a \emph{bottom}, if $t-\sigma^{-1}(i)+1\in J^+_\sigma$ ($-x_i$ is the bottom of a segment);
a \emph{block} otherwise ($x_i$ gives the Jordan block $\rho\boxtimes S_{2x_i+1}$). To \emph{replace} an
exponent $x_i$ by a value $y$ in $J_{\sigma,\delta}$ means to apply the same rule with $y$ at the position
$\sigma^{-1}(i)$, the segments being then reordered by their exponents.

\subsection{The geometric degree}\label{sec:geom}
Let $\widehat G$, $V$ and $B$ be as in \S\ref{ssec:types}. The infinitesimal character of $\pi(\caL)$,
which is also that of the costandard modules $J\in\Phi(\caL)$ (they have the same cuspidal support), is
the representation $\lambda_\caL$ of $W_F$ on $V$ defined by
$\lambda_\caL(w)=\phi_J(w,\mathrm{diag}(|w|^{1/2},|w|^{-1/2}))$, for any $J\in\Phi(\caL)$ with Langlands
parameter $\phi_J$. It has the form
\[\lambda_\caL=\bigoplus_\rho\bigoplus_{x}\rho|\cdot|^x\otimes V_{\rho,x},\]
where $\rho$ runs over the self-dual irreducible representations of $W_F$, $x$ over $\frac12\bbZ$, and the
$V_{\rho,x}$ are multiplicity spaces. Put $M_\rho=\bigoplus_xV_{\rho,x}$, so that $V=\bigoplus_\rho M_\rho
\otimes\rho$. The form $B$ induces on $M_\rho$ a non-degenerate form $B_\rho$, of symmetry
$\mathrm{sym}(B)\,\mathrm{sym}(\rho)$, which pairs $V_{\rho,x}$ with $V_{\rho,-x}$. Let $s_\rho$ be the
semisimple automorphism of $M_\rho$ acting by $q^x$ on $V_{\rho,x}$, and put
\[G_0=Z_{\widehat G}(\lambda_\caL(W_F))=\widehat G\cap\prod_\rho Z_{\mathrm{I}(M_\rho)}(s_\rho),\]
\[\widehat{\mathfrak g}_{-1}=\bigoplus_\rho\{N\in\mathrm{Lie}\,\mathrm{I}(M_\rho):N(V_{\rho,x})\subset
V_{\rho,x-1}\text{ for all }x\},\]
where $\mathrm{I}(M_\rho)$ is the isometry group of $(M_\rho,B_\rho)$, acting on $M_\rho\otimes\rho$ through
the first factor. The Langlands parameters $\phi$ of $G_n$ with infinitesimal character $\lambda_\caL$
are determined, up to conjugacy by $\widehat G$, by the nilpotent element
$N_\phi=\log\phi(1,\left(\begin{smallmatrix}1&1\\0&1\end{smallmatrix}\right))\in\widehat{\mathfrak g}_{-1}$, up to
conjugacy by $G_0$, and $G_0$ has finitely many orbits on $\widehat{\mathfrak g}_{-1}$. The irreducible
representations of $G_n$ with this infinitesimal character, and the costandard modules, are
parametrized by pairs $(O,\xi)$ of a $G_0$-orbit $O\subset\widehat{\mathfrak g}_{-1}$ and an irreducible
$G_0$-equivariant local system $\xi$ on $O$, relevant for $G_n$ (Arthur's local Langlands correspondence \cite{Art13} and
the Langlands classification; see \cite{AMS22}, Theorems C and D, and \cite{Sol25}, \S5).

\begin{rmk}[Infinitesimal characters]\label{rmk:infchar} This definition of $\lambda_\caL$ is Vogan's
(\cite{Vog93}; see also \cite{Hai14}). The Bernstein centre (\cite{BD84}) plays the role of
$Z(\mathfrak g)$, and $\lambda_\caL$ only depends on the cuspidal support, by the compatibility of the
local Langlands correspondence with cuspidal supports, proved for $G_n$ by M\oe glin and Xu
(\cite{Xu17}). The space $\widehat{\mathfrak g}_{-1}$ with the action of $G_0$ is Vogan's variety of
$L$-parameters with infinitesimal character $\lambda_\caL$, the $p$-adic analogue of the variety of
geometric parameters of Adams, Barbasch and Vogan \cite{ABV92}.
\end{rmk}

For a costandard module $J=\Delta_1\times\cdots\times\Delta_k\rtimes\pi(\phi,\varepsilon)$, the nilpotent
element $N_J$ of its Langlands parameter is described by its Jordan strings in the spaces $M_\rho$:
\begin{itemize}
\item a segment $\Delta_\rho[x,y]$ gives the string of weights $x,x-1,\ldots,y$ in $M_\rho$, together with
the dual string $-y,-y-1,\ldots,-x$;
\item a Jordan block $\rho\boxtimes S_{2x+1}$ of $\phi$ gives the string $x,x-1,\ldots,-x$ in $M_\rho$.
\end{itemize}
Let $O_J=G_0\cdot N_J$, and $d(J)=\dim O_J$. It is the rank of the map $X\mapsto[X,N_J]$ from
$\widehat{\mathfrak g}_0=\mathrm{Lie}(G_0)$ to $\widehat{\mathfrak g}_{-1}$.

For a string $a$ (a sequence of weights $\alpha_a,\alpha_a-1,\ldots,\omega_a$), let $a^\vee$ be the dual string
($-\omega_a,\ldots,-\alpha_a$), and for two strings put
\[h(a,b)=1\text{ if }\omega_a\leq \omega_b\leq \alpha_a\leq \alpha_b,\qquad h(a,b)=0\text{ otherwise}.\]

\begin{lemma}\label{lem:dimformula} Assume that the strings of $N=N_J$ are pairwise distinct (this is
the case if $J$ is of discrete type). Then
\[d(J)=\dim\widehat{\mathfrak g}_0-\tfrac12\Big(\sum_{a,b}h(a,b)-\sum_a h(a,a^\vee)\Big),\]
the first sum running over the pairs of strings of $N$ in the same space $M_\rho$, the second over the
strings of $N$.
\end{lemma}

\begin{proof}
$d(J)=\dim\widehat{\mathfrak g}_0-\dim Z_{\widehat{\mathfrak g}_0}(N)$, and $Z_{\widehat{\mathfrak g}_0}(N)$ is
the direct sum over $\rho$ of the centralizers of the components $N_\rho$ in $\mathrm{Lie}\,\mathrm{I}(M_\rho)$
of degree $0$. Fix $\rho$. A degree-preserving map from the string $a$ to the string $b$ commuting with
$N$ is determined by the image $w\in b$ of the top vector of $a$, of weight $\alpha_a$, subject to
$N^{\alpha_a-\omega_a+1}w=0$; so the space of such maps has dimension $h(a,b)$, and $\dim
Z_{\mathfrak{gl}(M_\rho)_0}(N_\rho)=\sum_{a,b}h(a,b)$, the sum running over the strings in $M_\rho$. Now
the centralizer in $\mathrm{Lie}\,\mathrm{I}(M_\rho)$ is the fixed space of the involution $\theta(X)=-X^*$
($X^*$ the adjoint for $B_\rho$) of $Z_{\mathfrak{gl}(M_\rho)_0}(N_\rho)$, of dimension
$\frac12(\dim+\mathrm{tr}\,\theta)$. The involution sends the maps $a\to b$ to the maps $b^\vee\to
a^\vee$; since the strings are distinct, only the pairs $(a,a^\vee)$ contribute to the trace. For
$X:a\to a^\vee$, the form $B_X(v,w)=B_\rho(Xv,w)$ on $a$ satisfies $B_X(Nv,w)=-B_X(v,Nw)$, so it has
the symmetry $(-1)^{2\alpha_a}$ (compare $B_X(e_k,N^{2k}e_k)$ and $B_X(N^{2k}e_k,e_k)$); on the other hand
$\theta(X)=\epsilon X$ gives the symmetry $-\epsilon\cdot\mathrm{sym}(B_\rho)$. The weights $\alpha_a$ lie in
$E_\rho$, so $(-1)^{2\alpha_a}=\mathrm{sym}(B)\,\mathrm{sym}(\rho)=\mathrm{sym}(B_\rho)$
(\S\ref{ssec:types}); hence $\epsilon=-1$, and $\mathrm{tr}\,\theta=-\sum_ah(a,a^\vee)$.
\end{proof}

\begin{prop}[Closure order]\label{prop:closure} Let $J,J'$ be costandard modules with the same
infinitesimal character. If $L(J')$ is a constituent of $J$, then $O_J\subset\overline{O_{J'}}$; in
particular $d(J')\geq d(J)$, with equality only if $L(J')=L(J)$. This applies when
$\mathrm{Hom}(J,J')\neq0$, since the image then contains $L(J')$.
\end{prop}

\begin{proof}
By the $p$-adic Kazhdan--Lusztig conjecture, proved by Solleveld for the classical groups
(\cite{Sol25}, Theorem 5.4, with Arthur's correspondence \cite{Art13}, the Hecke algebras of the
Bernstein blocks of \cite{Hei17} and the comparison results of \cite{AMS22}), the multiplicity
$[J:L(J')]$ is the multiplicity of the local system $\xi$ of $J$ in the cohomology sheaves, restricted to
$O_J$, of the intersection cohomology complex $\mathrm{IC}(O_{J'},\xi')$. In \cite{Sol25}, the variety of
Langlands parameters is $\{y:\mathrm{Ad}(\exp\sigma)y=q^{-1}y\}$ ((5.6) there), which is
$\widehat{\mathfrak g}_{-1}$; the group acting on it is a finite cover of $Z_{\widehat
G}(\lambda_\caL(W_F))$, with the same orbits; and the standard modules are the analytic standard
modules of \cite{Sol25}, Theorem 3.4, which correspond to the standard representations of $G_n$
(\cite{Sol25}, Theorem 5.4(a)). The multiplicities vanish unless $O_J\subset\overline{O_{J'}}$. If
$O_J=O_{J'}$, the stalk is the local system $\xi'$ itself, so $\xi=\xi'$ and $L(J')=L(J)$. Taking duals
exchanges the standard modules with simple head and the costandard modules with simple socle, and
preserves the multiplicities, so the statement holds in either convention.
\end{proof}

For $J\in\Phi(\caL)$, put $\deg J=d(J)-d(J_0)$. By Proposition \ref{prop:closure}, the matrix of
the multiplicities $[J:L(J')]$ is unitriangular for the closure order; so is its inverse, and
$c_J\neq0$ implies $O_{J_0}\subset\overline{O_J}$. Hence $\deg J\geq0$ for $J\in\Phi(\caL)$, with
equality only for $J=J_0$.

\subsection{The main theorem}
The resolution we are looking for is the analogue of a Cousin resolution: a coresolution of
$\pi(\caL)$ by costandard modules, in which the costandard modules sit in the degree given by the
codimension of their orbits.

\begin{thm}[Main theorem]\label{conj:geom} Let $\caL$ be a ladder datum.
\begin{enumerate}
\item For every $J\in\Phi(\caL)$, $c_J=(-1)^{\deg J}$.
\item There is an exact sequence
\[0\to\pi(\caL)\to J_0\to\bigoplus_{\deg J=1}J\to\cdots\to\bigoplus_{\deg J=k}J\to\cdots,\]
the sums running over $J\in\Phi(\caL)$.
\end{enumerate}
\end{thm}

\noindent The Euler characteristic of the sequence in (2) is $\sum_{J\in\Phi(\caL)}(-1)^{\deg J}[J]$,
which is Theorem \ref{thm:atobe} by (1). Atobe's formula is recovered as follows. Its terms are the
summands of the $I_\sigma(\caL)$, with the sign $\mathrm{sgn}(\sigma)=(-1)^{\ell(\sigma)}$, and $c_J$ is the
sum of these signs over the summands isomorphic to $J$.
 Part (1) contains Atobe's expectation that the coefficients are $\pm1$; the sign
$(-1)^{\deg J}$ is the one predicted by the geometry: when the stalks of the intersection cohomology
complexes have the parity property, the inverse of the matrix of multiplicities has signs
$(-1)^{\dim O_J-\dim O_{J'}}$.

\begin{rmk}\label{rem:twins} The Bruhat grading and the geometric grading do differ. For
$\mathbf{Sp}(8)$, $X=\{0,1,2\}$, $l=1$, one summand of $I_{(321)}(\caL)$ (length $3$) is the discrete series
with blocks $0,1,2$ and signs $(+,+,+)$, which also occurs in $I_{(231)}(\caL)$ and $I_{(312)}(\caL)$
(length $2$). Its coefficient is $2-1=1$, and its degree is $2$.
\end{rmk}

\section{Atobe's involution and derivatives}\label{sec:step1}

\subsection{Atobe's two lemmas}
We recall the structure of the proof of \cite{Ato24}, Theorem 4.2 (\S 4.3), whose module version is
the subject of the following sections. Let $\frs$ be the cuspidal support of $\pi(\caL)$, and put
\[\Sigma(\caL)=p_\frs\Big(\sum_{\sigma\in\frS(\caL)}\mathrm{sgn}(\sigma)[I_\sigma(\caL)]\Big),\]
so that Theorem \ref{thm:atobe} states that $[\pi(\caL)]=\Sigma(\caL)$. The proof is by induction on
$n$. If $\pi(\caL)$ is supercuspidal, then $\frS(\caL)$ has one element and the formula is trivial
(\cite{Ato24}, Corollary 3.7(1)). Otherwise, Atobe shows that $D_{\rho|\cdot|^\lambda}[\pi(\caL)]=
D_{\rho|\cdot|^\lambda}(\Sigma(\caL))$ for every self-dual $\rho$ and every $\lambda\in\mathbb{R}$, by the two
lemmas below. Applying them repeatedly, $[\pi(\caL)]$ and $\Sigma(\caL)$ have the same Jacquet module
with respect to the parabolic subgroup from which $\pi(\caL)$ is supercuspidally induced; since both
are linear combinations of standard modules with cuspidal support $\frs$, they are equal (\cite{Ato24},
end of \S4.3).

\begin{lemma}[\cite{Ato24}, Lemma 4.9 and Proposition 3.6]\label{lem:ato49} Assume that
$D_{\rho|\cdot|^\lambda}[\pi(\caL)]\neq0$. Then $\lambda=x_j$ for some $j$ and $\lambda-1\notin X_\rho$.
Let $\caL'$ be the ladder datum obtained by replacing $x_j$ by $x_j-1$. Then
$D_{\rho|\cdot|^\lambda}[\pi(\caL)]=[\pi(\caL')]$, $\frS(\caL')=\frS(\caL)$,
$D_{\rho|\cdot|^\lambda}[I_\sigma(\caL)]=[I_\sigma(\caL')]$ for all $\sigma$, and
$D_{\rho|\cdot|^\lambda}\circ p_\frs=p_{\frs'}\circ D_{\rho|\cdot|^\lambda}$, where $\frs'$ is the cuspidal
support of $\pi(\caL')$. In particular, by the induction hypothesis for $\caL'$,
$D_{\rho|\cdot|^\lambda}(\Sigma(\caL))=\Sigma(\caL')=[\pi(\caL')]=D_{\rho|\cdot|^\lambda}[\pi(\caL)]$.
\end{lemma}

\begin{lemma}[\cite{Ato24}, Lemma 4.8]\label{lem:ato48} Assume that
$D_{\rho|\cdot|^\lambda}[\pi(\caL)]=0$. Then $D_{\rho|\cdot|^\lambda}(\Sigma(\caL))=0$. More precisely,
if some $D_{\rho|\cdot|^\lambda}[I_\sigma(\caL)]$ is non-zero, then $\lambda=x_j$ with $x_{j-1}=x_j-1$, by
Tadi\'c's formula and \cite{Xu17}, Lemma 7.3, and Atobe defines an involution
$\sigma\mapsto\sigma'$ of the set $\frS_\lambda(\caL)$ of the $\sigma\in\frS(\caL)$ with $D_{\rho|\cdot|^\lambda}[I_\sigma(\caL)]\neq0$,
with $\mathrm{sgn}(\sigma')=-\mathrm{sgn}(\sigma)$ and
$D_{\rho|\cdot|^\lambda}[I_\sigma(\caL)]=D_{\rho|\cdot|^\lambda}[I_{\sigma'}(\caL)]$:
\begin{itemize}
\item[(1)] if $j-1=\sigma(i)$ and $j=\sigma(i')$ with $i<i'\leq l$, then
$\sigma'=\sigma\cdot(t-i+1,t-i'+1)$: the indices paired with $j-1$ and $j$ are exchanged (when
$\lambda-1$ and $\lambda$ are the tops of two segments, $I_\sigma(\caL)$ and $I_{\sigma'}(\caL)$ differ by
exchanging the other ends of these segments);
\item[(2)] otherwise $\sigma'=(j-1,j)\cdot\sigma$. In the delicate subcase where
$j-1=\sigma(i)$ and $j=\sigma(t-i+1)$, one has
\[I_\sigma(\caL)=\Delta[\lambda-1,-\lambda]\times\tau\rtimes\bigoplus_k\pi(\phi,\varepsilon_k),\qquad
I_{\sigma'}(\caL)=\tau\rtimes\bigoplus_k\big(\pi(\phi',\varepsilon'_{k,+})\oplus\pi(\phi',\varepsilon'_{k,-})\big),\]
with $\phi'=\phi\oplus\rho\boxtimes S_{2\lambda-1}\oplus\rho\boxtimes S_{2\lambda+1}$. Both
derivatives equal
\[\Delta[\lambda-1,-(\lambda-1)]\times\tau\rtimes\bigoplus_k\pi(\phi,\varepsilon_k).\]
\end{itemize}
The contributions of $\sigma$ and $\sigma'$ to $D_{\rho|\cdot|^\lambda}(\Sigma(\caL))$ cancel, whence the
first assertion.
\end{lemma}

\noindent At the level of modules, the involution has the following meaning, which is established in
Sections \ref{sec:local} and \ref{sec:degsigns}. Fix a paired exponent $\lambda$, i.e.\ $\lambda,\lambda-1\in
X_\rho$, so that $r^{\rho,\lambda}(\pi(\caL))=0$. A summand $J$ with $r^{\rho,\lambda}(J)\neq0$ and its
partner $J^\natural$, the summand which corresponds to it under Atobe's involution (\S\ref{ssec:paired}), are two
\emph{distinct} costandard modules with opposite coefficients, unless $\lambda-1$ and $\lambda$ are Jordan
blocks with the same sign in both, in which case $J^\natural\simeq J$ and $c_J=0$ (Lemma
\ref{lem:coefpaired}); so if one of them is in $\Phi(\caL)$, both are, their degrees differ by one (Lemma
\ref{lem:degpaired}), and both occur in the coresolution of Theorem \ref{conj:geom}. What they have
in common is their derivative: $r^{\rho,\lambda}(J)\simeq r^{\rho,\lambda}(J^\natural)$ (Proposition \ref{prop:cat48}). When $r^{\rho,\lambda}$ is applied to the coresolution, the
component of the differential between $J$ and $J^\natural$ becomes an isomorphism (Lemma
\ref{lem:matchediso} and the proof of Proposition \ref{prop:johnson}), so that the members of each pair
cancel each other and the complex obtained by applying $r^{\rho,\lambda}$ is exact, in accordance with $r^{\rho,\lambda}(\pi(\caL))=0$.

\subsection{Exactness through derivatives}
Let $r^{\rho,\lambda}=r^{\rho|\cdot|^\lambda}$ be the exact functor of \S\ref{ssec:deriv}. Let $C^\bullet$ be a complex of
representations of finite length, all of whose irreducible subquotients have a non-supercuspidal
cuspidal support. Then $C^\bullet$ is exact as soon as $r^{\rho,\lambda}(C^\bullet)$ is exact for all
pairs $(\rho,\lambda)$ (with $\rho$ in the cuspidal supports involved). This applies to the
representations with the cuspidal support $\frs$ of $\pi(\caL)$ as soon as $\pi(\caL)$ is not
supercuspidal; when $\pi(\caL)$ is supercuspidal, $\Phi(\caL)=\{J_0\}$ with $J_0=\pi(\caL)$, and Theorem
\ref{conj:geom} is trivial. Indeed:
\begin{itemize}
\item $r^{\rho,\lambda}$ is exact, so $r^{\rho,\lambda}(H^j(C^\bullet))=H^j(r^{\rho,\lambda}(C^\bullet))=0$;
\item the family of the $r^{\rho,\lambda}$ detects the non-zero representations of this kind: if an
irreducible subquotient $\pi$ embeds into $\rho_1|\cdot|^{\lambda_1}\times\cdots\rtimes\sigma$ ($\sigma$
supercuspidal), then $r^{\rho_1,\lambda_1}(\pi)\neq0$ by Frobenius reciprocity.
\end{itemize}
In Sections \ref{sec:local} to \ref{sec:johnson}, we study one functor $r^{\rho,\lambda}$ at a time, and
we fix $\rho$; we write $X,l,\eta$ for $X_\rho,l_\rho,\eta_\rho$, $\Delta[x,y]$ for $\Delta_\rho[x,y]$, and
$r^\lambda$ for $r^{\rho,\lambda}$. We identify $X$ with the set of exponents $x_1<\cdots<x_t$; an exponent
$\lambda\in X$ is \emph{unpaired} if $\lambda-1\notin X$, and \emph{paired} if $\lambda-1\in X$. The data attached to the other self-dual $\rho'$ are carried along
unchanged: a segment $\Delta_{\rho'}[x,y]$, $\rho'\not\simeq\rho$, contains no $\rho|\cdot|^\lambda$ (two
self-dual cuspidal representations which are twists of each other are equal), so the
$\rho|\cdot|^\lambda$-part of its Jacquet modules vanishes, and, by Tadi\'c's formula
(\S\ref{sec:prelim}), $r^{\rho,\lambda}$ only modifies the $\rho$-components of the summands.

\section{Derivatives of the summands}\label{sec:local}
In this section we compute the derivatives $r^\lambda(J)$ of the summands $J$ of $\caL$, and we compare
them with the derivatives of their socles $L(J)$. For the latter, we use the following terminology. For
a summand $J$ and an exponent $\lambda$, we say that $J$ \emph{satisfies (S) at $\lambda$} if
\[r^\lambda(J)\neq0\ \Longrightarrow\ r^\lambda(L(J))\neq0,\]
i.e.\ if the derivative $r^\lambda$ does not kill the socle of $J$ when it does not kill $J$. This is not
true in general: it is the property established in Proposition \ref{thm:SAM}
and Lemma \ref{lem:smallpaired}, in the cases needed below.

\subsection{Unpaired exponents}\label{ssec:unpaired}
Recall (\S\ref{sec:step1}) that $\rho$ is fixed, that we write $\caL=(X,l,\eta)$ for its component at
$\rho$, the other components being carried along unchanged, and that $X=\{x_1<\cdots<x_t\}$ is the set of
exponents. For an unpaired exponent $\lambda$ ($\lambda\in X$, $\lambda-1\notin X$), we denote by
$\caL_\lambda=(X_\lambda,l,\eta)$, $X_\lambda=(X\smallsetminus\{\lambda\})\cup\{\lambda-1\}$, the datum obtained
by replacing the exponent $\lambda$ by $\lambda-1$. By Lemma \ref{lem:ato49}, if
$D_{\rho|\cdot|^\lambda}(\pi(\caL))\neq0$, then $\caL_\lambda$ is a ladder datum and
$D_{\rho|\cdot|^\lambda}(\pi(\caL))=\pi(\caL_\lambda)$.

Fix $\sigma\in\frS(\caL)$, a summand
$J=\Delta_1\times\cdots\times\Delta_k\rtimes T$ of $I_\sigma(\caL)$, with
$\Delta_j=\Delta[x_{\sigma(j)},-x_{\sigma(t-j+1)}]$ ($j\in J^+_\sigma$) and $T=\pi(\phi_\sigma,
\varepsilon_{\sigma,\delta})$, and $\lambda\in\mathbb{R}$. By the geometric lemma, $r^\lambda(J)$ has
a filtration whose subquotients are of three kinds:
\begin{itemize}
\item[(t)] $\rho|\cdot|^\lambda\otimes(\cdots\times\Delta[\lambda-1,y]\times\cdots\rtimes T)$, for
a segment $\Delta[\lambda,y]$ of $J$ (top);
\item[(b)] $\rho|\cdot|^\lambda\otimes(\cdots\times\Delta[x,-\lambda+1]\times\cdots\rtimes T)$, for
a segment $\Delta[x,-\lambda]$ of $J$ (bottom, through the contragredient in $M^*$);
\item[(d)] $\rho|\cdot|^\lambda\otimes(\Delta_1\times\cdots\times\Delta_k\rtimes
D_{\rho|\cdot|^\lambda}(T))$ (derivative of the tempered part).
\end{itemize}

\begin{lemma}\label{lem:onepiece} At most one subquotient of the filtration of $r^\lambda(J)$ is
non-zero. Hence $r^\lambda(J)$ is $0$, or isomorphic to the corresponding induced representation.
\end{lemma}

\begin{proof}
Each exponent $x_i$ is exactly one of a top, a bottom or a block of $J$ (\S\ref{ssec:summands}). A
subquotient of type (t) (resp.\ (b)) requires $\lambda=x_i$ with $x_i$ a top (resp.\ a bottom) of $J$. A
subquotient of type (d) requires $\rho\boxtimes S_{2\lambda+1}\subset\phi_\sigma$
(\cite{Xu17}, Lemma 7.3, as in \cite{Ato24}, proof of Lemma 4.8), i.e.\ $\lambda=x_i$ with $x_i$ a block
of $J$. Since the $x_i$ are distinct, at most one of these occurs.
\end{proof}

\begin{prop}[Module version of \cite{Ato24}, Lemma 4.9]\label{prop:cat49} Let $\lambda$ be an
unpaired exponent with $\lambda\geq1$ or $\lambda\leq-1$, and assume that $\caL_\lambda$ is a ladder
datum. Then $\frS(\caL_\lambda)=\frS(\caL)$. Moreover, for every $\sigma$ and every sign choice $\delta$,
the summands $J=J_{\sigma,\delta}(\caL)$ and $J^\lambda=J_{\sigma,\delta}(\caL_\lambda)$ satisfy
\[r^\lambda(J)\simeq\rho|\cdot|^\lambda\otimes J^\lambda.\]
Here $J^\lambda=0$ if one of its segments is $0$. If $\lambda\leq-1$, then $\lambda$ is a top in every
summand, and $r^\lambda(J)=0$ exactly when $J$ contains the empty segment $\Delta[\lambda,\lambda+1]$,
i.e.\ when $-\lambda-1\in X$ is paired with $\lambda$ in the same segment.
\end{prop}

\begin{proof}
Assume first $\lambda\geq1$. Since $\lambda-1\notin X$, the order of the elements of $X_\lambda$ is that of $X$: $\lambda-1$ takes the
index $c$ of $\lambda=x_c$. So the conditions of \cite{Ato24}, Definition 4.1, which only involve
the indices and the elements $x_i\leq-1$ or $x_i=-\tfrac12$, are the same for $\caL$ and $\caL_\lambda$
(note that $\lambda-1\geq0$). Hence $\frS(\caL_\lambda)=\frS(\caL)$. For the same reason, the sets
$J^\pm_\sigma$, which compare $\sigma(j)$ and $\sigma(t-j+1)$, are the same. The signs of the
middle Jordan blocks only depend on the positions, and those of the pairs in $J^-_\sigma$ on
$\delta$. Thus $J^\lambda$ is obtained from $J$ by replacing the exponent $\lambda=x_c$ by $\lambda-1$
(\S\ref{ssec:summands}). According as $\lambda$ is a top, a bottom or a block of $J$:
\begin{itemize}
\item[(t)] $\lambda$ is the top of a segment $\Delta[\lambda,y]$, which becomes $\Delta[\lambda-1,y]$;
\item[(b)] $\lambda$ is a bottom, of a segment $\Delta[x,-\lambda]$, which becomes
$\Delta[x,-\lambda+1]$;
\item[(d)] $\rho\boxtimes S_{2\lambda+1}$ is a Jordan block of $T$, which becomes
$\rho\boxtimes S_{2\lambda-1}$, with the same sign.
\end{itemize}
By Lemma \ref{lem:onepiece}, $r^\lambda(J)$ is isomorphic to the unique non-zero subquotient of the
geometric lemma filtration, or is $0$. In the three cases this subquotient is, respectively:
\begin{itemize}
\item[(t)] $\rho|\cdot|^\lambda\otimes(\Delta_1\times\cdots\times\Delta[\lambda-1,y]\times\cdots
\rtimes T)$, since the $\rho|\cdot|^\lambda$-part of the Jacquet module of $\Delta[\lambda,y]$ with
respect to $\mathbf{GL}_d\times\mathbf{GL}_{*}$ is $\rho|\cdot|^\lambda\otimes\Delta[\lambda-1,y]$;
\item[(b)] $\rho|\cdot|^\lambda\otimes(\cdots\times\Delta[x,-\lambda+1]\times\cdots\rtimes T)$, the
element $\rho|\cdot|^\lambda$ being the contragredient of the last element $\rho|\cdot|^{-\lambda}$ of
$\Delta[x,-\lambda]$ in $M^*$;
\item[(d)] $\rho|\cdot|^\lambda\otimes(\Delta_1\times\cdots\rtimes D_{\rho|\cdot|^\lambda}(T))$.
Since $\rho\boxtimes S_{2\lambda+1}\subset\phi$ and $\rho\boxtimes S_{2\lambda-1}\not\subset\phi$
(because $\lambda-1\notin X$), \cite{Xu17}, Lemma 7.3, gives
$D_{\rho|\cdot|^\lambda}(\pi(\phi,\varepsilon))=\pi(\phi',\varepsilon')$, where $S_{2\lambda+1}$ is
replaced by $S_{2\lambda-1}$ with the same sign. This is the discrete series of $J^\lambda$.
\end{itemize}
In all cases, the segments appear in the order of $J$, whereas $J^\lambda$ is defined with its segments in
non-decreasing order of $x+y$. The two orders can only differ for a pair formed by the modified
segment $A'$ and another segment $B$ which was tied with $A$ (since the sums $x+y$ are integers and
change by $1$). In case (t), $A=[y,\lambda]$ and $B=[y',x']$ have the same centre, hence are
nested. If $x'>\lambda$, then $B\supset A\supset A'$. If $x'<\lambda$, then $x'\leq\lambda-1$ and
$y'\geq y$, so $A'\supset B$. Case (b) is symmetric. In both cases $A'$ and $B$ are nested, hence not
linked, so $A'\times B\simeq B\times A'$ (\cite{Zel80}, Theorem 9.7). This proves
$r^\lambda(J)\simeq\rho|\cdot|^\lambda\otimes J^\lambda$.

Finally, a segment of $J^\lambda$ is $0$ if and only if it is $\Delta[\lambda-1,\lambda+1]$ or
$\Delta[-\lambda-1,-\lambda+1]$. The corresponding segment of $J$ is then the empty segment
$\Delta[\lambda,\lambda+1]$ or $\Delta[-\lambda-1,-\lambda]$. In both cases $\rho|\cdot|^\lambda$ does not
occur in the segments of $J$, and $r^\lambda(J)=0$. Both cases require $-\lambda-1\in X$: in the first, the bottom
$\lambda+1$ is $-x$ for some $x\in X$; in the second, $-\lambda-1$ is a top. So if
$-\lambda-1\notin X$, then $r^\lambda(J)\neq0$ for all summands.

Assume now $\lambda\leq-1$. An element $\leq-1$ lies at a position of L, so $\lambda$ is a top in every
summand, and $\lambda$ and $\lambda-1\leq-2$ are subject to the same conditions of \cite{Ato24},
Definition 4.1; hence $\frS(\caL_\lambda)=\frS(\caL)$. Only case (t) occurs, and the argument above
applies; the summand is killed exactly when the segment with top $\lambda$ is the empty segment
$\Delta[\lambda,\lambda+1]$, whose bottom is $-x$ with $x=-\lambda-1\in X$.
\end{proof}

\begin{rmk} The cuspidal supports satisfy $\mathrm{supp}(J)=\mathrm{supp}(J^\lambda)+\rho|\cdot|^\lambda$.
So the projections $p_\frs$ and $p_{\frs'}$ of \cite{Ato24}, Theorem 4.2, keep or discard $J$ and
$J^\lambda$ simultaneously. The unpaired exponents $\lambda=\tfrac12$ ($\rho$ of half-integral type) and
$\lambda=0$ ($\rho$ of integral type) are treated in Lemma \ref{lem:small}.
\end{rmk}

\begin{lemma}\label{lem:small} Let $\mu=\tfrac12$ ($\rho$ of half-integral type, with $-\tfrac12\notin X$) or
$\mu=0$ ($\rho$ of integral type, with $-1\notin X$), and assume $D_{\rho|\cdot|^\mu}(\pi(\caL))\neq0$.
Then $r^\mu(J)\simeq\rho|\cdot|^\mu\otimes J^\mu$ for every summand $J$, where $J^\mu$ is the
corresponding summand of $\caL_\mu$, or $0$. The summands with $r^\mu(J)\neq0$ correspond bijectively to
the summands of $\caL_\mu$.
\end{lemma}

\begin{proof}
By Lemma \ref{lem:onepiece}, $r^\mu(J)$ is $0$ or isomorphic to the corresponding induced representation.

\emph{Case $\mu=\tfrac12$.} The exponent $\tfrac12$ is not a bottom: the segment would be $\Delta[x,-\tfrac12]$
with $x<\tfrac12$ in $X$, hence $x\leq-\tfrac32$, and the segment would be $0$.
\begin{itemize}
\item If $\tfrac12$ is a top, the summand of $\caL_\mu$ has $\Delta[-\tfrac12,y]$ instead of
$\Delta[\tfrac12,y]$, as in Proposition \ref{prop:cat49}.
\item If $\tfrac12$ is a Jordan block, i.e.\ $S_2\subset\phi$ with sign $s$, then by \cite{AM23},
Theorem 5.3 (with $m'=1$ and $\eta(S_0)=+1$), $D_{\rho|\cdot|^{1/2}}(\pi(\phi,\varepsilon))\neq0$ if and
only if $s=+1$. In that case it is $\pi(\phi-S_2,\varepsilon)$, i.e.\ $S_2$ is replaced by $S_0$, which
is dropped. On the side of $\caL_\mu$, $-\tfrac12$ replaces $\tfrac12$. There are two cases.
\begin{itemize}
\item If $\tfrac12$ is in the middle part, it is at the first middle position, since the other middle
elements are $>\tfrac12$. So its sign is $s=\eta$. The corrected Definition 4.1 of \cite{Ato24}
(Remark \ref{rmk:typo}) excludes these $\sigma$ from $\frS(\caL_\mu)$ exactly when $\eta=-1$, and
drops $S_0$ when $\eta=+1$.
\item If $\tfrac12$ is the smaller element of a pair of $J^-$, its sign is $s=\delta$. In
$I_\sigma(\caL_\mu)$, the summand with $\delta=-1$ vanishes, and $S_0$ is dropped for $\delta=+1$
(\cite{Ato24}, Definition 4.1).
\end{itemize}
\end{itemize}
In all cases, the vanishing of $r^\mu(J)$ matches the absence of the corresponding summand of
$\caL_\mu$. The order of the segments is handled as in Proposition \ref{prop:cat49}.

\emph{Case $\mu=0$.} The derivative $r^0$ of a discrete series vanishes, by Casselman's criterion, so
$r^0(J)=0$ if $0$ is a block of $J$. A bottom $0$ only occurs in the empty segment $\Delta[-1,0]$, which is
excluded since $-1\notin X$. So $0$ is a top, of $\Delta[0,y]$ with $y<0$, which becomes
$\Delta[-1,y]$ in $\caL_\mu$. In $\caL_\mu$ the element $-1$ must be at a position of L
(\cite{Ato24}, Definition 4.1). So $\frS(\caL_\mu)$ consists of the $\sigma$ where $0$ is a top, which
are those with $r^0(I_\sigma(\caL))\neq0$.
\end{proof}

\subsection{Paired exponents}\label{ssec:paired}
\begin{prop}[Module version of \cite{Ato24}, Lemma 4.8]\label{prop:cat48} Let $\lambda=x_j\geq1$ be a
paired exponent ($x_{j-1}=\lambda-1$), let $\sigma\mapsto\sigma'$ be Atobe's involution of
$\frS_\lambda(\caL)$ (Lemma \ref{lem:ato48}), and $\delta\in\{\pm1\}^{J^-_\sigma}$ a sign choice.
\begin{enumerate}
\item[(i)] Assume that $j-1$ and $j$ are not the two members $\sigma(k),\sigma(t-k+1)$ of a pair,
$k\leq l$. Let $\delta'=\delta$ in case (2) of Lemma \ref{lem:ato48}, and, in case (1),
$\delta'=\delta\circ(i\ i')$, where $j-1=\sigma(i)$ and $j=\sigma(i')$. Then $\delta'$ is a sign choice for
$\sigma'$ (in case (2), $J^\pm_{\sigma'}=J^\pm_\sigma$), $J_{\sigma',\delta'}$ is obtained from $J_{\sigma,\delta}$ by
exchanging the values $\lambda-1$ and $\lambda$ (\S\ref{ssec:summands}), and
\[r^\lambda(J_{\sigma,\delta})\simeq r^\lambda(J_{\sigma',\delta'}).\]
\item[(ii)] Assume that $j-1=\sigma(k)$ and $j=\sigma(t-k+1)$ for some $k\leq l$, so that $k\in J^+_\sigma$ and
$J_{\sigma,\delta}$ contains the segment $\Delta[\lambda-1,-\lambda]$. Then $\sigma'=(j-1,j)\sigma$,
$J^-_{\sigma'}=J^-_\sigma\cup\{k\}$, and, if $\delta_\pm$ denotes the sign choice for $\sigma'$ which extends
$\delta$ by $\delta_\pm(k)=\pm1$ (the common sign of the Jordan blocks $S_{2\lambda-1},S_{2\lambda+1}$),
\[r^\lambda(J_{\sigma,\delta})\simeq r^\lambda(J_{\sigma',\delta_+})\oplus r^\lambda(J_{\sigma',\delta_-}).\]
\end{enumerate}
Since $\sigma\mapsto\sigma'$ is an involution, (i) and (ii) cover all the $\sigma\in\frS_\lambda(\caL)$: the case
$j=\sigma(k)$, $j-1=\sigma(t-k+1)$ is (ii) applied to $\sigma'$. The cases $\lambda=\tfrac12$ and $\lambda=0$
are treated in Lemma \ref{lem:smallpaired}.
\end{prop}

\begin{proof}
\emph{Step 1: comparing $J_{\sigma,\delta}$ and $J_{\sigma',\delta'}$ in (i).} We show that
$J_{\sigma',\delta'}$ is obtained from $J_{\sigma,\delta}$ by exchanging the values $\lambda-1$ and $\lambda$
(\S\ref{ssec:summands}), the signs of the Jordan blocks staying at their places.
\begin{itemize}
\item In case (2), $\sigma'=(j-1,j)\sigma$. Since $j-1$ and $j$ are adjacent values, exchanging them does
not change any comparison $\sigma(k)<\sigma(t-k+1)$, except for a pair made of $j-1$ and $j$ themselves,
which is excluded. So $J^\pm_{\sigma'}=J^\pm_\sigma$, each position gives a top, a bottom or a block in the
same way for $\sigma$ and $\sigma'$, with the same sign, and only the values at the positions
$\sigma^{-1}(j-1)$ and $\sigma^{-1}(j)$ are exchanged.
\item In case (1), let $b=\sigma(t-i+1)$ and $a=\sigma(t-i'+1)$ be the indices paired with $j-1$ and $j$ in
$\sigma$; then $\sigma'=\sigma\cdot(t-i+1,t-i'+1)$ pairs $j-1$ with $a$ (at $i$) and $j$ with $b$ (at
$i'$), and coincides with $\sigma$ elsewhere. Since $a,b\notin\{j-1,j\}$, the comparison of $j$ with $b$
(resp.\ of $j-1$ with $a$) is that of $j-1$ with $b$ (resp.\ of $j$ with $a$). So the pair $(j,b)$ of
$\sigma'$ gives the segment $\Delta[\lambda,-x_b]$ if the pair $(j-1,b)$ of $\sigma$ gives
$\Delta[\lambda-1,-x_b]$, and two Jordan blocks with the sign $\delta'(i')=\delta(i)$ if the latter gives two
Jordan blocks with the sign $\delta(i)$; similarly for the pairs $(j-1,a)$ and $(j,a)$. In particular
$J^-_{\sigma'}=(i\ i')(J^-_\sigma)$, and $\delta'$ is a sign choice for $\sigma'$.
\end{itemize}

\emph{Step 2: the derivatives in (i).} By Lemma \ref{lem:onepiece},
$r^\lambda(J_{\sigma,\delta})\simeq\rho|\cdot|^\lambda\otimes\widetilde J_\sigma$, where $\widetilde J_\sigma$ is
obtained from $J_{\sigma,\delta}$ by replacing $\lambda$ by $\lambda-1$, its factors being taken in the
order of $J_{\sigma,\delta}$, and similarly for $\sigma'$. By Step 1, in both $\widetilde J_\sigma$ and
$\widetilde J_{\sigma'}$ the value $\lambda-1$ is placed at the two places of $\lambda-1$ and $\lambda$ in
$J_{\sigma,\delta}$. When one of these places is a Jordan block, we use \cite{Xu17}, Lemma 7.3: for a
discrete series $\pi(\phi,\varepsilon)$ with $\rho\boxtimes S_{2\lambda+1}\subset\phi$,
\begin{itemize}
\item if $S_{2\lambda-1}\not\subset\phi$, then $D_{\rho|\cdot|^\lambda}\pi(\phi,\varepsilon)$ is
the discrete series where $S_{2\lambda+1}$ is replaced by $S_{2\lambda-1}$, with the same sign;
\item if $S_{2\lambda-1}\subset\phi$, the derivative is non-zero if and only if the two signs agree.
It is then the irreducible tempered representation with parameter
$\phi-S_{2\lambda+1}+S_{2\lambda-1}$ and the common sign.
\end{itemize}
In all cases the result only depends on the places and their signs. So $\widetilde J_\sigma$ and
$\widetilde J_{\sigma'}$ have the same tempered part and the same segments, possibly in different
orders, and $r^\lambda(J_{\sigma,\delta})\neq0$ if and only if $r^\lambda(J_{\sigma',\delta'})\neq0$.

\emph{Step 3: the order of the segments in (i).} We show that the two orders differ only by exchanges
of non-linked segments, which give isomorphic products (\cite{Zel80}, Theorem 9.7). Order the segments by
their centre $\frac{x+y}2\in\frac12\bbZ$. At most two segments are modified, and the centre of a modified
segment changes by $\frac12$, so its position relative to an unmodified segment $C$ can only change if
their centres are equal before or after the modification.
\begin{itemize}
\item If $C$ has the same centre as a modified segment $A'$, they are nested.
\item If $C$ has the same centre $c$ as the unmodified segment and $A'$ has centre $c\pm\frac12$, write
$C=[y_C,x_C]$. Its top is not $\lambda-1$ or $\lambda$, and its bottom is not $-\lambda$ or
$-\lambda+1$, since these values are taken by the modified segments. One checks as in the proof
of Proposition \ref{prop:cat49} that $C$ and $A'$ are nested. For instance, if $\lambda$ is a bottom,
$A'=[-\lambda+1,x_A]$ and $x_C+y_C=x_A-\lambda$. If $x_C>x_A$, then $y_C<-\lambda$ and
$C\supset A'$. If $x_C<x_A$, then $y_C\geq-\lambda+2$ and $C\subset A'$.
\end{itemize}
If two segments $A',B'$ are modified, then either $\lambda-1$ and $\lambda$ are two tops (resp.\ two
bottoms), and $A',B'$ share their top $\lambda-1$ (resp.\ their bottom $-\lambda+1$), hence are nested; or
one is a top and the other a bottom, and the centres of $A'$ and $B'$ change by the same amount, so that
their relative order is the same. This proves (i).

\emph{Step 4: case (ii).} Here $\widetilde J_\sigma=\tau_1\times\Delta[\lambda-1,-\lambda+1]\times\tau_2\rtimes T$,
where $\tau_1,\tau_2$ are the other segments. A segment $C=[y_C,x_C]$ of $\tau_2$ has centre
$\geq-\frac12$ and $x_C\notin\{\lambda-1,\lambda\}$.
\begin{itemize}
\item If $x_C\geq\lambda+1$, then $y_C<-\lambda$, and $C\supset[-\lambda+1,\lambda-1]$.
\item If $x_C\leq\lambda-2$, then $y_C\geq-1-x_C\geq-\lambda+1$, and $C\subset[-\lambda+1,\lambda-1]$.
\end{itemize}
So $\Delta[\lambda-1,-\lambda+1]$ commutes with $\tau_2$, and $r^\lambda(J_{\sigma,\delta})\simeq
\rho|\cdot|^\lambda\otimes\big(\tau_1\times\tau_2\rtimes(\Delta[\lambda-1,-\lambda+1]\rtimes
T)\big)$. The representation $\Delta[\lambda-1,-\lambda+1]\rtimes T$ is unitary, hence semisimple. It
is the sum $T_+\oplus T_-$ of the two tempered representations with parameter
$\phi\oplus\rho\boxtimes S_{2\lambda-1}^{\oplus2}$ whose characters extend $\varepsilon$; there are two
of them, since $S_{2\lambda-1}$ is of good parity and does not occur in $\phi$ (M\oe glin's construction
of the tempered $L$-packets, \cite{Ato20}, Theorem 3.3; \cite{MT02}). On the other side, as in Step 1,
$J_{\sigma',\delta_\pm}$ has the same segments $\tau_1,\tau_2$ (up to order) and the discrete series
$\pi(\phi',\varepsilon'_\pm)$, $\phi'=\phi\oplus\rho\boxtimes(S_{2\lambda-1}\oplus S_{2\lambda+1})$, the two
new blocks having the sign $\pm$. By \cite{Xu17}, Lemma 7.3,
$D_{\rho|\cdot|^\lambda}(\pi(\phi',\varepsilon'_\pm))=T_\pm$. Hence $r^\lambda(J_{\sigma',\delta_+})\oplus
r^\lambda(J_{\sigma',\delta_-})\simeq\rho|\cdot|^\lambda\otimes(\tau_1\times\tau_2\rtimes(T_+\oplus T_-))$,
the segments of $\tau_1\times\tau_2$ being in the same order up to commutations, as in Step 3.

Finally, since $\mathrm{supp}(J)=\mathrm{supp}(\widetilde J)+\rho|\cdot|^\lambda$, the projections
$p_\frs$ act in the same way on both sides.
\end{proof}

\noindent\emph{Pairs.} Let $\lambda\geq1$ be a paired exponent. Atobe's involution exchanges two
permutations which differ by a transposition, hence have different lengths; we denote by $v$ the shorter one
and by $u$ the longer one, and write $v<u$. In case (ii) of Proposition \ref{prop:cat48}, $v$ is the
permutation whose summands contain the segment $\Delta[\lambda-1,-\lambda]$ (Lemma \ref{lem:betapair}). 
For a sign choice $\delta$ for $v$, the summands $J_v=J_{v,\delta}$ and $J_u=J_{u,\delta'}$, with $\delta'$
as in Proposition \ref{prop:cat48}(i) (in case (ii), $J_v$ and the two summands
$J_u^\pm=J_{u,\delta_\pm}$), form a \emph{pair}; $J_v$ is its \emph{lower member} and $J_u$ its \emph{upper
member}, and each is the \emph{partner} of the other, denoted $J^\natural$. By Proposition \ref{prop:cat48},
$r^\lambda(J_v)\simeq r^\lambda(J_u)$ (resp.\ $r^\lambda(J_v)\simeq r^\lambda(J_u^+)\oplus r^\lambda(J_u^-)$), and
we shall see that $\deg J_u=\deg J_v+1$ when $J_v$ and $J_u$ do not have the same strings
(Lemma \ref{lem:degpaired}).

\subsection{Socles}\label{ssec:socles}
We now study property (S), defined at the beginning of this section.

\begin{lemma}\label{lem:nonvanishing} Let $f:J\to J'$ be non-zero, where $J'$ has simple socle
$L'$. If $r^\lambda(L')\neq0$, then $r^\lambda(f)\neq0$. If moreover $r^\lambda(J)\simeq
r^\lambda(J')$, and $r^\lambda(J')$ has a simple socle occurring once in it, then $r^\lambda(f)$
is an isomorphism.
\end{lemma}

\begin{proof}
The image of $f$ contains $L'$, and $r^\lambda$ is exact, so the image of $r^\lambda(f)$
contains $r^\lambda(L')\neq0$. A non-zero endomorphism $\psi$ of a module of finite length whose
socle is simple and occurs once is injective: otherwise its kernel and its image would both contain
the socle. Hence it is bijective.
\end{proof}

\begin{lemma}\label{lem:embed} Let $\lambda\in X$ with $\lambda-1\notin X$, and let
$J=\Delta_1\times\cdots\times\Delta_k\rtimes T$ be a summand in which $\lambda$ is the top of a segment
$\Delta_i=\Delta[\lambda,y]$. Then
\[J\hookrightarrow\rho|\cdot|^\lambda\times(\Delta_1\times\cdots\times\Delta[\lambda-1,y]\times\cdots\rtimes T).\]
In particular $r^\lambda(L(J))\neq0$, by Frobenius
reciprocity.
\end{lemma}

\begin{proof}
$\Delta_i\hookrightarrow\rho|\cdot|^\lambda\times\Delta[\lambda-1,y]$. Let $j<i$ and
$\Delta_j=\Delta[x',y']$, so that $x'+y'\leq\lambda+y$. The segments $\{\lambda\}$ and $[y',x']$ are
linked only if $x'=\lambda-1$ or $y'=\lambda+1$. The first case is excluded since $x'\in X$. In the
second case, $\Delta_j\neq0$ gives $x'\geq\lambda$, whereas $x'\leq y-1\leq\lambda-1$; this is
impossible. Hence $\rho|\cdot|^\lambda$ commutes with $\Delta_1,\ldots,\Delta_{i-1}$ (\cite{Zel80},
Theorem 9.7).
\end{proof}

\noindent The paper \cite{AM23} computes the highest $\rho|\cdot|^x$-derivative of an irreducible
representation from its Langlands data: Theorem 7.1 for $x>0$, and Proposition 6.1 for $x<0$. This
gives the analogue of Lemma \ref{lem:embed} for tops, bottoms and blocks, including $\lambda=\tfrac12$,
$\lambda<0$, and the paired exponents.

Let $J=\Delta_1\times\cdots\times\Delta_k\rtimes T$ be a summand of $\caL$, so that
$L(J)=L(\Delta_1,\ldots,\Delta_k;T)$, and let $x=\lambda>0$. In the notation of \cite{AM23}, \S7.1:
\begin{itemize}
\item $t_0\in\{0,1\}$ (denoted $t$ there) is the number of segments equal to $\Delta[\lambda-1,-\lambda]$;
\item $A_x$ and $A_{x-1}$ are the sets of segments with top $\lambda$, resp.\ $\lambda-1$;
\item $B_x$ and $B_{x-1}$ are the sets of segments with bottom $-\lambda$, resp.\ $-(\lambda-1)$
(excluding $\Delta[\lambda-1,-\lambda]$);
\item $m$ and $m'$ are the multiplicities of $S_{2\lambda+1}$ and $S_{2\lambda-1}$ in the parameter
$\psi$ of \cite{AM23}, \S7.1. For $x=\tfrac12$, the convention there is $m'=1$ and
$\eta(S_0)=+1$ (\cite{AM23}, Theorem 5.3).
\end{itemize}
Since each exponent is exactly one of a top, a bottom or a block of $J$, all these sets have at most one element,
and $m,m'\leq1$.
The best matching function (\cite{AM23}, \S6.1), as used in \S7.1 there, matches some segments of $A_x$ with segments of
$A_{x-1}$, and some segments of $B_x$ with segments of $B_{x-1}$; $A^c_x$, $A^c_{x-1}$, $B^c_x$ denote the
sets of unmatched segments.
\begin{itemize}
\item The two segments of $A_x$ and $A_{x-1}$ are matched by the best matching function if and only
if the bottom of the segment of top $\lambda$ is larger than that of the segment of top
$\lambda-1$.
\item The two segments of $B_x$ and $B_{x-1}$ are matched if and only if the top of the segment of
bottom $-\lambda$ is smaller than that of the segment of bottom $-(\lambda-1)$.
\end{itemize}
By \cite{AM23}, Theorem 7.1, the order of the highest derivative is
\[k=|A^c_x|+\max\big\{m+\max\{|B^c_x|-m',0\}-|A^c_{x-1}|,0\big\}.\]

\begin{lemma}\label{lem:vacuous} Let $\lambda=x_j$ with $x_{j-1}=\lambda-1$, and let $u$ be the upper
member of a pair of Atobe's involution. Then the following never happens:
\begin{itemize}
\item[$(\ast)$] $\lambda$ is a bottom, and $\lambda-1$ is the smaller
element of a pair of $J^-_u$ placed at a position of R after that of $\lambda$.
\end{itemize}
The same holds for the lower member, with $\lambda-1$ and $\lambda$ exchanged.
\end{lemma}

\begin{proof}
In subcase $(\ast)$, let $p<q$ be positions of R. Here $u(p)=j$ is the bottom of the segment
$k=t-p+1\in J^+_u$, so $u(k)<j$. And $u(q)=j-1$ is the smaller member of a pair
$k'=t-q+1\in J^-_u$, so $u(k')>j-1$, hence $u(k')\geq j+1$, since $j$ is at $p$. Now $k'<k$ are
positions of L, where $u$ is increasing, so $u(k)>u(k')\geq j+1$. This contradicts $u(k)<j$.

The dual subcase exchanges $j-1$ and $j$: we have $v(k)<j-1$ and $v(k')>j$, which is again
incompatible with $v(k)>v(k')$.
\end{proof}

\begin{prop}\label{thm:SAM} Let $\lambda\in X$, and $J$ a summand of $\caL$ with
$r^\lambda(J)\neq0$. Then $r^\lambda(L(J))\neq0$, i.e.\ $J$ satisfies (S) at $\lambda$, in each of the
following cases:
\begin{enumerate}
\item[(a)] $\lambda-1\notin X$ (unpaired $\lambda$, including $\lambda=\tfrac12$, $\lambda=0$ and $\lambda<0$);
\item[(b)] $\lambda-1\in X$, $\lambda\neq0$, and $J$ is the upper member of a pair (\S\ref{ssec:paired}).
\end{enumerate}
\end{prop}

\begin{proof}
Assume first $\lambda>0$. We compute $k$ according as $\lambda$ (and $\lambda-1$) are tops, bottoms or blocks.

(a) Here $A_{x-1}=B_{x-1}=\emptyset$, $t_0=0$, and $m'=0$ unless $\lambda=\tfrac12$.
\begin{itemize}
\item If $\lambda$ is a top, then $|A^c_x|=1$.
\item If $\lambda$ is a Jordan block, then $D_{\rho|\cdot|^\lambda}(T)\neq0$, so $m=1$ in $\psi$
(for $\lambda=\tfrac12$ this means $\eta(S_2)=+1$). Then $k=\max\{1+\max\{0-m',0\},0\}=1$.
\item If $\lambda$ is a bottom, then $k=\max\{1-m',0\}$. This is $1$ for
$\lambda\geq1$. For $\lambda=\tfrac12$ this case does not occur when $-\tfrac12\notin X$: the top of
the segment would be an element of $X$ smaller than $\tfrac12$, hence $\leq-\tfrac32$, and the
segment $\Delta[x,-\tfrac12]$ would be $0$.
\end{itemize}

(b) In the upper member, $j$ precedes $j-1$ (case (2) of Lemma \ref{lem:ato48}), or the tops
$\lambda-1,\lambda$ are associated with the bottoms in decreasing order (case (1) of Lemma \ref{lem:ato48}). The possible configurations are the following.
\begin{itemize}
\item Two tops: the upper member pairs $\lambda$ with the smaller bottom, so $A_x$ and
$A_{x-1}$ are not matched, and $|A^c_x|=1$.
\item $\lambda$ a top and $\lambda-1$ a bottom or a block: $|A^c_x|=1$.
\item $\lambda$ a block and $\lambda-1$ a bottom: $m=1$, $m'=0$, $A_{x-1}=\emptyset$, so $k=1$.
\item $\lambda$ a bottom and $\lambda-1$ a bottom: in the upper member the segment of bottom
$-\lambda$ has the larger top, so $B_x,B_{x-1}$ are not matched, $|B^c_x|=1$, $m'=0$, and $k=1$.
\item $\lambda$ and $\lambda-1$ both blocks, including case (ii) of Proposition \ref{prop:cat48}: $r^\lambda(J)\neq0$ forces equal
signs, and then $\psi=\phi$, $m=m'=1$, and $k=1$.
\item $\lambda$ a bottom and $\lambda-1$ a block: this is subcase $(\ast)$, which does not occur
(Lemma \ref{lem:vacuous}).
\end{itemize}
For $\lambda<0$, $\lambda$ is a top in every summand, and \cite{AM23}, Proposition 6.1, gives $k=|A^c_x|$.
This is $1$ in (a), and also in (b), since $\lambda$ and $\lambda-1$ are then two tops, not matched in the
upper member.

For $\lambda=0$ in (a) ($\rho$ of integral type, $-1\notin X$), $0$ is a top in every summand $J$
with $r^0(J)\neq0$ (proof of Lemma \ref{lem:small}), and Lemma \ref{lem:embed} applies.
\end{proof}

\noindent For the lower members, $r^\lambda(L(J))=0$ (Lemmas \ref{lem:lowersocle} and
\ref{lem:lowersocle0}). The case $\lambda=0$ is treated directly in Lemma \ref{lem:smallpaired}. The
normalizations of \cite{AM23} and \cite{Ato24} agree: both papers use Arthur's correspondence for the
split group $G_n$ in the form of \cite{Ato20} (\cite{AM23}, \S\S2.4 and 5), and \cite{AM23}, Theorem 5.3, uses
the convention of \S\ref{ssec:types} for the blocks $\rho\boxtimes S_0$.

\begin{lemma}\label{lem:matchediso} Let $\lambda=x_j\neq0$ with $x_{j-1}=\lambda-1$, and let
$\{J_v,J_u\}$ be a pair of Atobe's involution (\S\ref{ssec:paired}), so that $r^\lambda(J_v)\simeq
r^\lambda(J_u)$ (Proposition \ref{prop:cat48}). Then for
every non-zero $f:J_v\to J_u$, the map $r^\lambda(f)$ is an isomorphism. In case (ii) of Proposition \ref{prop:cat48}, the same holds for $(f_+,f_-):J_v\to J_u^+\oplus J_u^-$ with
$f_\pm\neq0$.
\end{lemma}

\begin{proof}
We have $r^\lambda(J_u)=\rho|\cdot|^\lambda\otimes\widetilde J_u$, where $\widetilde J_u$ is
obtained by replacing $\lambda$ by $\lambda-1$.
\begin{itemize}
\item In the upper member, this does not create a segment of centre $0$. That would require the
segment $\Delta[\lambda-1,-\lambda]$, which belongs to the lower member in case (ii) of Proposition \ref{prop:cat48}.
\item $D_{\rho|\cdot|^\lambda}(T)$ is irreducible tempered (\cite{Xu17}, Lemma 7.3).
\end{itemize}
So $\widetilde J_u$ is a costandard module, with simple socle occurring once (Langlands
classification). Lemma \ref{lem:nonvanishing} and Proposition \ref{thm:SAM}(b) conclude.

In case (ii) of Proposition \ref{prop:cat48}, we have $r^\lambda(J_v)\simeq\rho|\cdot|^\lambda\otimes(\tau\rtimes T_+\oplus
\tau\rtimes T_-)$. Moreover $\mathrm{Hom}(\tau\rtimes T_\mp,\tau\rtimes T_\pm)=0$, since the
socle $L(\tau;T_\pm)$ of $\tau\rtimes T_\pm$ is not a constituent of the costandard module $\tau\rtimes T_\mp$:
the leading exponent of the latter only carries $T_\mp$. So $r^\lambda(f_\pm)$ is zero on one
summand, and an isomorphism on the other one.
\end{proof}

\begin{lemma}\label{lem:smallpaired} Let $\lambda=\tfrac12$ with $-\tfrac12\in X$ ($\rho$ of half-integral type), or
$\lambda=0$ with $-1\in X$ ($\rho$ of integral type). Then Proposition \ref{prop:cat48} holds for $\lambda$,
and so do property (S) for the upper summands and Lemma \ref{lem:matchediso}.
\end{lemma}

\begin{proof}
\emph{Case $\lambda=0$.} Since $-1\leq-1$, the element $-1$ is always a top. Moreover $r^0(J)\neq0$ forces
$0$ to be a top, as in Lemma \ref{lem:small}. So every $\sigma\in\frS_\lambda(\caL)$ is in case (1), Atobe's involution
is $\sigma'=\sigma\cdot(t-i,t-i+1)$, and Proposition \ref{prop:cat48}(i) applies. In the upper member
the two tops are not matched, so the segment of top $-1$ comes after that of top $0$, and the
argument of Lemma \ref{lem:embed} (top case) gives (S). Lemma \ref{lem:matchediso} follows, since $\widetilde J_u$
is again a costandard module.

\emph{Case $\lambda=\tfrac12$.} The element $-\tfrac12$ can be a top, or a Jordan block $S_0$. In the
latter case its sign is $+1$, since the summands where it would be $-1$ vanish (\cite{Ato24},
Definition 4.1). So $S_0$ behaves like a block with sign $+1$, which is the convention of
\cite{AM23}, Theorem 5.3. The proof of Proposition \ref{prop:cat48} then goes
through, with the blocks $S_0$ dropped at the end.
\begin{itemize}
\item In case (ii) of Proposition \ref{prop:cat48}, the segment $\Delta[\lambda-1,-\lambda]=\Delta[-\tfrac12,-\tfrac12]$
corresponds to the pair of blocks $S_2\oplus S_0$. Only the sign $+1$ survives, and both
derivatives equal $\rho|\cdot|^{1/2}\otimes\tau\rtimes T$, since $\Delta[-\tfrac12,\tfrac12]$ is
empty. So there is no splitting into $T_\pm$.
\item Property (S) for the upper summands is Proposition \ref{thm:SAM}(b), whose computation with
$m'=1$ is valid for $x=\tfrac12$.
\end{itemize}
\end{proof}

\section{Degrees and signs}\label{sec:degsigns}
In this section, we compare the degrees of the summands of $\caL$ and of the data $\caL_\mu$ obtained by
the derivatives, using Lemma \ref{lem:dimformula}, and we prove Theorem \ref{conj:geom}(1).

\subsection{Degrees of the unpaired summands}
\begin{lemma}\label{lem:degunpaired} Let $\mu$ be an unpaired exponent ($\mu=x_i$, $\mu-1\notin X$),
$J$ a summand with $r^\mu(J)\neq0$, and $J^\mu$ the summand of $\caL_\mu$ with
$r^\mu(J)\simeq\rho|\cdot|^\mu\otimes J^\mu$. Then $d(J)-d(J^\mu)=\dim\widehat{\mathfrak g}_0-
\dim\widehat{\mathfrak g}'_0-m$, where $m$ is the number of indices $k\neq i$ with $x_k=-\mu$. In
particular $\deg J^\mu=\deg J$.
\end{lemma}

\begin{proof}
Every index $k$ contributes exactly one string with top $x_k$ and one string with bottom $-x_k$: the
segment string and its dual, or the self-dual string of a block (the index of $S_0$, if any,
contributes nothing and is irrelevant below). By Proposition \ref{prop:cat49} and Lemma \ref{lem:small}, $N_{J^\mu}$ is
obtained from $N_J$ by two operations: remove the top weight $\mu$ of the unique string $c$ with top
$\mu$, then remove the bottom weight $-\mu$ of the unique string with bottom $-\mu$ (which is
$c^\vee$, or $c$ itself for a block). No string has top $\mu-1$ or bottom $-(\mu-1)$, since
$\mu-1\notin X$. If $\mu=\tfrac12$, then $-\mu=\mu-1\notin X$, so $m=0$ and the index of $S_0$ does not
interfere.

If $\mu=0$ ($\rho$ of integral type), only case (t) of Lemma \ref{lem:onepiece} occurs: case (b) would require
the empty segment $\Delta[-1,0]$, and case (d) a block $S_1$ with $D_\rho(\tau)\neq0$, which is
impossible for $\tau$ square integrable ($\tau\hookrightarrow\rho\rtimes\tau'$ would force $S_1$ to occur
twice in $\phi_\tau$). Then $c=\langle0,\ldots,-x'\rangle$ with $x'>0$, $c^\vee=\langle x',\ldots,0\rangle$ is the only other string with bottom $0$, and the
computation below gives a change of $-1$ for $h_{\mathrm{tot}}$ and of $-1$ for $h_\vee$, hence $0$ for
$\dim Z_{\widehat{\mathfrak g}_0}(N)$. We now assume $\mu\neq0$; then $k=i$ never satisfies
$x_k=-\mu$.

Let $h_{\mathrm{tot}}=\sum_{a,b}h(a,b)$. If $c'$ is $c$ with its top lowered from $\mu$ to $\mu-1$, then for
$b\neq c$, $h(c',b)-h(c,b)=[\omega_c\leq \omega_b]([\alpha_b=\mu-1]-[\omega_b=\mu])$ and $h(b,c')-h(b,c)=-[\omega_b\leq
\omega_c\leq \alpha_b][\alpha_b=\mu]$, and $h(c',c')-h(c,c)=-[\omega_c=\mu]$. Since no other string has top $\mu$ or
$\mu-1$, and $\omega_c\leq\mu$, the first operation changes $h_{\mathrm{tot}}$ by $-\#\{b:\omega_b=\mu\}=-m$. Dually
($h(a,b)=h(b^\vee,a^\vee)$), the second operation changes $h_{\mathrm{tot}}$ by minus the number of strings with top
$-\mu$ after the first operation, that is $-m-[\mu=\tfrac12]\,e$, where $e=1$ if $c'$ is non-empty with
top $-\tfrac12$.

Let $h_\vee=\sum_ah(a,a^\vee)$, with $h(a,a^\vee)=1$ if and only if $\alpha_a\geq0$ and $\omega_a\leq -\alpha_a$. Only $c$
and $c^\vee$ change. In the three cases of Lemma \ref{lem:onepiece}:
(t) $c=\langle\mu,\ldots,-x'\rangle$ with $x'>\mu$ the exponent paired with $\mu$: if $\mu>0$, then
$h(c,c^\vee)=1$, $h(c^\vee,c)=0$, and after the operations $h=[\mu\geq1]$ and $0$; if $\mu\leq-1$, both
terms vanish before and after;
(b) $c=\langle\mu,\ldots,-a\rangle$, the dual of $\Delta[a,-\mu]$ with $a<\mu$: $h(c,c^\vee)=0$ before and after (since
$a\neq\mu-1$), and $h(c^\vee,c)=[a\geq0]$ before and after;
(d) $c=c^\vee=\langle\mu,\ldots,-\mu\rangle$: $h(c,c)=1$, and after the operations $h=[\mu\geq1]$.
In each case $h_\vee$ changes by $-[\mu=\tfrac12]$, and $e=[\mu=\tfrac12]$ (in case (b) with $\mu=\tfrac12$,
$a=-\tfrac12$ and $c'$ is empty, and $h_\vee$ does not change). By Lemma \ref{lem:dimformula},
$\dim Z_{\widehat{\mathfrak g}_0}(N)=\frac12(h_{\mathrm{tot}}-h_\vee)$ changes by $-m$ in all cases. Since
$\widehat{\mathfrak g}_0$ and $\widehat{\mathfrak g}'_0$ only depend on the infinitesimal characters, we
get the formula, and $\deg J^\mu=\deg J$ by applying it to $J$ and to $J_0$ (which is not killed by $r^\mu$, since
$r^\mu(\pi(\caL))\neq0$).
\end{proof}

\subsection{Degrees of the paired summands}
We now compare the degrees of the two members of a pair of Atobe's involution. The computation is
an application of Lemma \ref{lem:dimformula}. Fix $\sigma\in\frS(\caL)$ and a summand $J$ of
$I_\sigma(\caL)$, and let $\iota_J$ be the involution of $\{1,\ldots,t\}$ which exchanges $\sigma(k)$ and
$\sigma(t-k+1)$ for $k\in J^+_\sigma$ and fixes the other indices. For each index $i$, let $s_i$ be the
string with top $\alpha_i=x_i$ and bottom $\omega_i=-x_{\iota_J(i)}$, that is, the weights $x_i,x_i-1,\ldots,
-x_{\iota_J(i)}$. It is empty if $\omega_i=\alpha_i+1$ (a segment $\Delta[x,x+1]$, or the block $S_0$).
\begin{itemize}
\item For $k\in J^+_\sigma$, $s_{\sigma(k)}$ and $s_{\sigma(t-k+1)}$ are the string of the segment
$\Delta[x_{\sigma(k)},-x_{\sigma(t-k+1)}]$ and its dual. For a fixed index $i$, $s_i$ is the string of the
block $S_{2x_i+1}$.
\item In all cases $s_{\iota_J(i)}=s_i^\vee$, and the non-empty $s_i$ are the strings of $N_J$. They are
pairwise distinct, since their tops are.
\end{itemize}
The formula defining $h(a,b)$ in terms of the endpoints gives $h(a,b)=0$ as soon as $a$ or $b$ is
empty (for $a$ empty, $\omega_a=\alpha_a+1$, and the conditions $\omega_a\leq \omega_b\leq \alpha_a$ or $\omega_b\leq \omega_a\leq
\alpha_b\leq \alpha_a$ cannot hold). So we may sum over all the indices in Lemma \ref{lem:dimformula}. Since
$\sum_ah(a,a^\vee)=\sum_ih(s_i,s_{\iota_J(i)})$, we get
\begin{equation}\label{eq:Zpi}
Z(J):=\dim\widehat{\mathfrak g}_0-d(J)=\frac12\sum_{(i,k):\ k\neq\iota_J(i)}h(s_i,s_k).
\end{equation}

Now let $\lambda=x_j$ be a paired exponent, $x_{j-1}=\lambda-1$ (so $\lambda\geq0$; $\lambda=0$ only
for $\rho$ of integral type, $\lambda=\tfrac12$ only for $\rho$ of half-integral type). For a summand $J$ put
\[\beta(J)=\omega_j=-x_{\iota_J(j)},\qquad\gamma(J)=\omega_{j-1}=-x_{\iota_J(j-1)},\]
the bottoms of the (possibly empty) strings of top $\lambda$ and $\lambda-1$.

\begin{lemma}\label{lem:betagamma} Let $J$ be a summand, and $r'$ the index $\neq j$ with $x_{r'}=-\lambda$,
if it exists.
\begin{enumerate}
\item[(a)] If $\iota_J(j-1)\notin\{j-1,j,r'\}$, then $\gamma(J)\leq\lambda-1$ and
$\gamma(J)\notin\{-\lambda,-\lambda+1\}$.
\item[(b)] If $\iota_J(j)\notin\{j-1,j,r'\}$ and the string of top $\lambda$ of $J$ is non-empty, then
$\beta(J)\leq\lambda-1$ and $\beta(J)\notin\{-\lambda,-\lambda+1\}$. The string of top $\lambda$ is non-empty
if $r^\lambda(J)\neq0$.
\end{enumerate}
\end{lemma}

\begin{proof}
All the endpoints lie in $\lambda+\bbZ$, and a string of top $\alpha$ has a bottom $\leq\alpha+1$, with
equality exactly when it is empty. So $\gamma(J)\leq\lambda$, and $\gamma(J)=\lambda$ means
$x_{\iota_J(j-1)}=-\lambda$, i.e.\ $\iota_J(j-1)\in\{j,r'\}$ ($j$ only if $\lambda=0$). Similarly
$\beta(J)\leq\lambda$ if the string of top $\lambda$ is non-empty, and $\beta(J)=\lambda$ means
$\iota_J(j)\in\{j,r'\}$. Finally, $-x_{\iota_J(j-1)}\in\{-\lambda,-\lambda+1\}$ means $\iota_J(j-1)\in\{j,j-1\}$,
and similarly for $\beta(J)$. If the string of top $\lambda$ is empty, the
corresponding factor of $J$ is $\Delta[\lambda,\lambda+1]=\mathbf 1$; then $J$ has no segment with top
$\lambda$, no segment with bottom $-\lambda$ (only the index $j$ has the value $\lambda$) and no block
$S_{2\lambda+1}$, so the three kinds of subquotients of \S\ref{ssec:unpaired} vanish, and $r^\lambda(J)=0$.
\end{proof}

\begin{lemma}\label{lem:betapair} Let $\lambda=x_j$ be a paired exponent, and $\{J_v,J_u\}$ a pair
(\S\ref{ssec:paired}) whose members do not have the same strings. Then $\beta(J_u)<\beta(J_v)$. Moreover
$\beta(J_u)<\gamma(J_u)$ and $\beta(J_v)>\gamma(J_v)$: the upper member is the one in which the string of top
$\lambda$ has the smaller bottom. In case (ii)
of Proposition \ref{prop:cat48}, $v$ is the permutation whose summands contain the segment
$\Delta[\lambda-1,-\lambda]$.
\end{lemma}

\begin{proof}
In case (ii) of Proposition \ref{prop:cat48} (with its notation $\sigma$, $k$), $\sigma'=(j-1,j)\sigma$ has $j$
at the position $k$ of L
and $j-1$ at $t-k+1$, so it is the longer one, and it gives the blocks $S_{2\lambda+1},S_{2\lambda-1}$: the
string of top $\lambda$ of $J_u^\pm$ is $\langle\lambda,\ldots,-\lambda\rangle$, whereas that of $J_v$ is
$\langle\lambda,\ldots,-\lambda+1\rangle$, so $\beta(J_u^\pm)=-\lambda<-\lambda+1=\beta(J_v)$. We now consider
case (i), and write $\beta$, $\beta'$ for the values of $\beta$ on the two members. Recall that $\iota_J$ pairs
the positions $k$ and $t-k+1$ when $\sigma(k)<\sigma(t-k+1)$ ($k\in J^+_\sigma$), and that $\sigma$ is
increasing on L and on the middle positions.

Case (1) of Atobe's involution: $j-1=\sigma(i)$, $j=\sigma(i')$, $i<i'\leq l$, and
$\sigma'=\sigma\cdot(a,b)$ with $a=t-i'+1<b=t-i+1$. Then $\ell(\sigma')>\ell(\sigma)$ if and only if
$\sigma(a)<\sigma(b)$, i.e.\ $\sigma'(a)>\sigma'(b)$: in the longer member, the index paired with $j$ is
larger. If $j$ is in $J^+$ in both members, $\beta=-x_{\sigma(a)}$ and $\beta'=-x_{\sigma(b)}$, and
we are done. If $j$ is in $J^+$ in one member only, the index $k$ paired with it is $>j$ there and $<j$ in the other,
so the first is the longer one, and there $\beta=-x_k<-\lambda$, whereas $\beta=-\lambda$
(block) in the other. If $j$ is in $J^-$ in both, so is $j-1$, and the two members have the same
strings.

Case (2): $\sigma'=(j-1,j)\sigma$, and the longer member $u$ is the one in which $j$ comes before
$j-1$. The positions of $j$ and $j-1$ are not both in L (that is case (1)), and not both in the middle
(by monotonicity, $j$ could not precede $j-1$). Let $p<q$ be the positions of $j$ and $j-1$ in $u$; by
Proposition \ref{prop:cat48}(i), the sets $J^\pm$ are the same for both members.
\begin{itemize}
\item $p\in\mathrm{L}$, $q\notin\mathrm{L}$. If $p\in J^+$, the index paired with $j$ in $u$ is $>j$, so
$\beta(J_u)<-\lambda$, whereas in $v$ the index $j$ sits at $q$: as a block ($\beta(J_v)=-\lambda$), as the
member of a pair of $J^-$ (block), or as the bottom of a segment whose top is $<j-1$
($\beta(J_v)>-\lambda+1$). If $p\in J^-$, $j$ is a block in $u$ ($\beta(J_u)=-\lambda$); if $j-1$ is a
block in $u$ as well, the two members have the same strings; otherwise $q\in\mathrm{R}$, paired with an index
$<j-1$, and in $v$ the index $j$ is the bottom of a segment with this top, so $\beta(J_v)>-\lambda+1$.
\item $p,q\in\mathrm{R}$: the indices paired with $j$ and $j-1$ are at the positions $t-p+1>t-q+1$ of L,
so the one paired with $j$ is larger. If both are in $J^+$, $\beta(J_u)=-x_{\iota(j)}<-x_{\iota(j-1)}=\beta(J_v)$. If both
are in $J^-$, the members have the same strings. The index paired with $j$ cannot be $<j$ while that paired with $j-1$
is $>j-1$. If $j$ is in $J^-$ and $j-1$ in $J^+$, then $\beta(J_u)=-\lambda$, and in $v$ the index $j$
is paired with the index paired with $j-1$ in $u$, which is $<j-1$, so $\beta(J_v)>-\lambda+1$.
\item $p$ in the middle, $q\in\mathrm{R}$: $\beta(J_u)=-\lambda$. In $v$, $j$ sits at $q$, paired with the
index paired with $j-1$ in $u$. If this index is $<j-1$, then $\beta(J_v)>-\lambda+1$; otherwise both indices are blocks
in both members.
\end{itemize}
In all cases with distinct strings, $\beta(J_u)<\beta(J_v)$.

For the last assertion, let $J$ be a member and $J^\natural$ the other one. In case (ii), $\beta(J_u^\pm)=-\lambda<
-\lambda+1=\gamma(J_u^\pm)$ and $\beta(J_v)=-\lambda+1>-\lambda=\gamma(J_v)$. In case (i), $J^\natural$ is obtained
from $J$ by exchanging the values $\lambda-1$ and $\lambda$ (Proposition \ref{prop:cat48}(i)), and $j-1$ and $j$
are not paired together; so $\beta(J^\natural)=\gamma(J)$, except when $\lambda-1$ is a block of $J$, where
$\gamma(J)=-\lambda+1$ and $\beta(J^\natural)=-\lambda$, and then $\lambda$ is not a block of $J$, so that
$\beta(J)\neq-\lambda$. In all cases $[\beta(J)<\gamma(J)]=[\beta(J)<\beta(J^\natural)]$, and the first
assertion concludes.
\end{proof}

\begin{lemma}\label{lem:degpaired} Let $\lambda$ be a paired exponent, and $\{J_v,J_u\}$ a pair
(\S\ref{ssec:paired}) whose members do not have the same strings. Then $d(J_u)=d(J_v)+1$. In particular
$\deg J_u=\deg J_v+1$.
\end{lemma}

\noindent The excluded case is that of isomorphic summands which cancel in Atobe's formula (\S\ref{ssec:summands}
and Remark \ref{rem:twins}): $\lambda-1$ and $\lambda$ are then Jordan blocks in both members, and
these summands do not contribute to $\Phi(\caL)$.

\begin{proof}
The algebras $\widehat{\mathfrak g}_0$ are the same for $J_v$ and $J_u$, so we have to show
$Z(J_v)-Z(J_u)=1$. By Proposition \ref{prop:cat48}, there are two
situations.
\begin{itemize}
\item[(E)] In case (i) of Proposition \ref{prop:cat48}, $J_u$ is obtained from $J_v$ by exchanging the values
$\lambda-1$ and $\lambda$ at their positions. In terms of strings, $\iota_{J_u}=\tau\iota_{J_v}\tau$ with
$\tau=(j-1,j)$: the tops $\lambda$ and $\lambda-1$ of $s_j$ and $s_{j-1}$ are exchanged, and so are the
bottoms $-\lambda$ and $-\lambda+1$ of $s_{\iota(j)}$ and $s_{\iota(j-1)}$, all other endpoints being
unchanged (empty strings included). We may then index the strings of both members by
$\{1,\ldots,t\}$ through the positions, i.e.\ use the same $\iota=\iota_{J_v}$ and the values $x\circ\tau$
for $J_u$.
\item[(D)] In case (ii) of Proposition \ref{prop:cat48}, $J_v$ contains the strings $\langle\lambda,\ldots,-\lambda+1\rangle$
and $\langle\lambda-1,\ldots,-\lambda\rangle$ ($\iota(j)=j-1$), and $J_u$ contains the strings of the
blocks $S_{2\lambda+1},S_{2\lambda-1}$ instead, the other strings being the same. Here
$\lambda\geq\tfrac12$.
\end{itemize}

\emph{Step 1: the comparisons which change.} Each $h(s_i,s_k)$ is determined by the three comparisons
$\omega_i\leq \omega_k$, $\omega_k\leq \alpha_i$, $\alpha_i\leq \alpha_k$. All endpoints lie in one coset of $\mathbb{Z}$, and in
(E) the moved endpoints change by $\pm1$: a top moves inside $\{\lambda-1,\lambda\}$, a bottom inside
$\{-\lambda,-\lambda+1\}$. Comparing a moved endpoint with a fixed value $c$, one finds that a
comparison can only change in the following cases: a moved top $w$ and a fixed bottom $c=\lambda$ in
$c\leq w$; a moved bottom $w$ and a fixed top $c=-\lambda$ in $w\leq c$. (A comparison $w\leq c$ of a
moved top with a fixed top changes only for $c=\lambda-1$, and $c\leq w$ only for $c=\lambda$; these
tops are moved. Similarly for two bottoms.) A fixed bottom $\lambda$ or a fixed top $-\lambda$ exists
only if $-\lambda=x_r$ for an index $r\notin\{j-1,j\}$, which forces $\lambda\geq1$; the bottom
$\lambda$ is then $\omega_{\iota(r)}$ and the top $-\lambda$ is $\alpha_r$. For two moved endpoints, the
comparisons which change are $\alpha_j$ versus $\alpha_{j-1}$, $\omega_{\iota(j)}$ versus $\omega_{\iota(j-1)}$, and
$\omega_{\iota(j)}\leq \alpha_j$, $\omega_{\iota(j-1)}\leq \alpha_{j-1}$ (only if $\lambda\leq\tfrac12$); the cross comparisons
$\omega_{\iota(j)}\leq \alpha_{j-1}$ and $\omega_{\iota(j-1)}\leq \alpha_j$ read $-\lambda\leq\lambda-1$ and
$-\lambda+1\leq\lambda$ before, and $-\lambda+1\leq\lambda$ and $-\lambda\leq\lambda-1$ after, so they
do not change. The comparisons $\omega_{\iota(i)}\leq \alpha_i$ only occur in the terms $h(s_i,s_{\iota(i)})$, which
are excluded from (\ref{eq:Zpi}). Hence, in (E), the only terms of (\ref{eq:Zpi}) which can change are
those of the pairs
\[(j,j-1),\ (j-1,j),\ (\iota(j),\iota(j-1)),\ (\iota(j-1),\iota(j)),\]
\[(j,\iota(r)),\ (j-1,\iota(r)),\ (r,\iota(j)),\ (r,\iota(j-1)),\]
those of the second line only when $r$ exists. By the symmetry $h(a,b)=h(b^\vee,a^\vee)$, each term of the
third and fourth pairs equals a term of the first two pairs, and each term of the last two equals a term
of the fifth and sixth, in both members.

In (E), put $\beta=\beta(J_v)$, $\gamma=\gamma(J_v)$, $\beta'=\beta(J_u)$, $\gamma'=\gamma(J_u)$. Then:
\begin{itemize}
\item if $\iota(j),\iota(j-1)\notin\{j-1,j\}$, then $\beta'=\gamma$ and $\gamma'=\beta$;
\item if $\iota(j)=j$, then $\beta=-\lambda$, $\gamma'=-\lambda+1$ and $\beta'=\gamma$;
\item if $\iota(j-1)=j-1$, then $\gamma=-\lambda+1$, $\beta'=-\lambda$ and $\gamma'=\beta$.
\end{itemize}
(The
cases $\iota(j)=j-1$ and $\iota(j)=j,\iota(j-1)=j-1$ give the same strings, and do not occur in (E).) We
show that, in (E),
\begin{equation}\label{eq:Zswap} Z(J_v)-Z(J_u)=[\beta'<\beta]-[\beta'>\beta],\qquad\beta\neq\beta'.
\end{equation}

\emph{Step 2: (E), when $r$ does not exist or $\iota(r)\notin\{j-1,j\}$.} By Lemma
\ref{lem:betagamma}, applied to $J_v$, we have $\beta\leq\lambda-1$ and $\beta\notin\{-\lambda,-\lambda+1\}$ if
$\iota(j)\notin\{j-1,j\}$, and $\gamma\leq\lambda-1$ and $\gamma\notin\{-\lambda,-\lambda+1\}$ if
$\iota(j-1)\notin\{j-1,j\}$. Indeed, here $\iota(j),\iota(j-1)\neq r$, and the string of top $\lambda$ of
$J_v$ is non-empty, since otherwise $J_u$ would contain the string $\langle\lambda-1,\ldots,\lambda+1\rangle$.
When $\iota(j)=j$ (resp.\ $\iota(j-1)=j-1$), $J_u$ (resp.\ $J_v$) contains the string
$\langle\lambda-1,\ldots,-\lambda+1\rangle$, which forces $\lambda\geq\tfrac12$.

\emph{The pairs $(j,\iota(r))$, $(j-1,\iota(r))$.} Note that $\iota(r)\neq r$, since $x_r=-\lambda\leq-1$
cannot be a Jordan block. Since $\omega_{\iota(r)}=\lambda$ and $\alpha_{\iota(r)}=x_{\iota(r)}\notin
\{\lambda-1,\lambda\}$, we get $h(s_j,s_{\iota(r)})=[\beta\leq\lambda][x_{\iota(r)}>\lambda]$ and
$h(s_{j-1},s_{\iota(r)})=0$ for $J_v$ (since $\alpha_{j-1}=\lambda-1<\lambda$), and the same values with
$\beta'$ in place of $\beta$ for the strings of top $\lambda$ and $\lambda-1$ of $J_u$. As
$\beta,\beta'\leq\lambda$, these terms do not change.

\emph{The pairs $(j,j-1)$, $(j-1,j)$.} For $J_v$, $h(s_j,s_{j-1})=0$ since $\alpha_j>\alpha_{j-1}$, and
$h(s_{j-1},s_j)=[\gamma\leq\beta\leq\lambda-1]$; for $J_u$ the values $\lambda$ and $\lambda-1$ are
exchanged. With the symmetric pairs, we get
\[Z(J_v)-Z(J_u)=[\gamma\leq\beta\leq\lambda-1]-[\gamma'\leq\beta'\leq\lambda-1].\]
If $\iota(j),\iota(j-1)\notin\{j-1,j\}$, this is $[\gamma<\beta]-[\beta<\gamma]$, which is
(\ref{eq:Zswap}) since $\beta'=\gamma$. If $\iota(j)=j$, it is $[\gamma\leq-\lambda]-[-\lambda+1\leq
\gamma]$ (note $-\lambda\leq\lambda-1$), and $\gamma\notin\{-\lambda,-\lambda+1\}$, which is
(\ref{eq:Zswap}) since $\beta=-\lambda$, $\beta'=\gamma$. If $\iota(j-1)=j-1$, it is $[-\lambda+1\leq\beta]
-[\beta\leq-\lambda]$ (for $\lambda=\tfrac12$ the first term vanishes, since $\beta\leq\lambda-1<-\lambda+1$),
which is (\ref{eq:Zswap}) since $\beta'=-\lambda$.

\emph{Step 3: (E), when $\iota(r)\in\{j-1,j\}$.} Exchanging $J_v$ and $J_u$ exchanges the indices of the
strings of top $\lambda$ and $\lambda-1$ and changes the sign of both sides of (\ref{eq:Zswap}), so we
may assume $\iota(r)=j$. Then $s_j=\langle\lambda\rangle$, $s_r=\langle-\lambda\rangle$, $\beta=\lambda$,
and in $J_u$ the corresponding strings $\langle\lambda-1,\lambda\rangle$ and
$\langle-\lambda,-\lambda+1\rangle$ are empty. The pairs of Step 1 are $(j,j-1)$, $(j-1,j)$,
$(r,\iota(j-1))$, $(\iota(j-1),r)$, $(j,j)$ and $(r,r)$. In $J_v$, $h(s_j,s_j)=h(s_r,s_r)=1$, and the other
terms vanish: $h(s_{j-1},s_j)=0$ since $\omega_j=\lambda>\alpha_{j-1}$, $h(s_j,s_{j-1})=0$ since $\alpha_j>\alpha_{j-1}$,
and $h(s_r,s_{\iota(j-1)})=h(s_{\iota(j-1)},s_r)=0$ since $\omega_{\iota(j-1)}=-\lambda+1>\alpha_r$ and
$\omega_{\iota(j-1)}>\omega_r$. In $J_u$ all these terms involve an empty string or vanish for the same reasons.
Hence $Z(J_v)-Z(J_u)=1$. On the other hand $\beta'\leq\lambda-1<\beta$: the string of top $\lambda$ of
$J_u$ is $\langle\lambda,\ldots,\gamma\rangle$, with $\gamma\leq\lambda-1$ by Lemma \ref{lem:betagamma}(a), or
$\langle\lambda,\ldots,-\lambda\rangle$. This is (\ref{eq:Zswap}).

\emph{Step 4: (D).} Only the bottoms of $s_j,s_{j-1}$ and the involution change. The pairs whose terms
can change are $(j,j)$, $(j-1,j-1)$ (in the sum for $J_v$ only), $(j,j-1)$, $(j-1,j)$ (in the sum for
$J_u$ only), and $(r,j)$, $(r,j-1)$ if $r$ exists (by Step 1 applied to the moved bottoms; for
$\lambda=\tfrac12$ the top $-\lambda$ is $\alpha_{j-1}$, already listed). For $J_v$,
$h(s_j,s_j)=h(s_{j-1},s_{j-1})=1$, and $h(s_r,s_j)=0$, $h(s_r,s_{j-1})=[\omega_r\leq-\lambda]$. For $J_u$,
$h(s_j,s_{j-1})=[-\lambda\leq-\lambda+1\leq\lambda\leq\lambda-1]=0$, $h(s_{j-1},s_j)=0$ (since
$-\lambda+1>-\lambda$), $h(s_r,s_j)=[\omega_r\leq-\lambda]$ and $h(s_r,s_{j-1})=0$. Hence
$Z(J_v)-Z(J_u)=1$.

\emph{Conclusion.} In (E), Lemma \ref{lem:betapair} gives $\beta'=\beta(J_u)<\beta(J_v)=\beta$, and
(\ref{eq:Zswap}) gives $Z(J_v)-Z(J_u)=1$.
\end{proof}

\noindent We now compare the degree of a member of a pair with that of its derivative; this is used in the
proofs of Lemma \ref{lem:coefpaired} and Proposition \ref{prop:johnson}. For a summand $J$ with
$r^\lambda(J)\neq0$, let $\widetilde J$ be the module of Lemma
\ref{lem:onepiece}, $r^\lambda(J)\simeq\rho|\cdot|^\lambda\otimes\widetilde J$ (in case (ii) of Proposition \ref{prop:cat48},
$\widetilde J_v=\widetilde J_u^+\oplus\widetilde J_u^-$, and we take $\widetilde J_u^\pm$). Its
strings are obtained from those of $J$ by lowering the top $\lambda$ of $s_j$ to $\lambda-1$ and raising
the bottom $-\lambda$ of $s_{\iota(j)}$ to $-\lambda+1$; that is, they are the strings $\widetilde s_i$
attached to the same involution $\iota=\iota_J$ and to the values $\widetilde x$, with $\widetilde
x_j=\lambda-1$ and $\widetilde x_i=x_i$ for $i\neq j$. Two of them may coincide (the strings of
$S_{2\lambda-1}$ in $\widetilde J_u^\pm$, case (ii) of Proposition \ref{prop:cat48}). Formula (\ref{eq:Zpi}) still holds for
$\widetilde J$, with the same $\iota$: in the proof of Lemma \ref{lem:dimformula}, the trace of $\theta$
only involves the maps from a string to its dual string, which is the string $\widetilde s_{\iota(i)}$, as
an orthogonal summand of $V$ (a block $S_{2\lambda-1}$ of $\phi$ is such a summand).

\begin{lemma}\label{lem:degtilde} Let $\lambda=x_j$ be a paired exponent, and $m$ the number of indices
$k\neq j$ with $x_k=-\lambda$. For every summand $J$ with $r^\lambda(J)\neq0$, whose pair
has members with distinct strings,
\[d(J)-d(\widetilde J)=\dim\widehat{\mathfrak g}_0-\dim\widetilde{\mathfrak g}_0-m+\mathrm{up}(J),\]
where $\widetilde{\mathfrak g}_0$ is the Lie algebra of the group $G_0$ attached to $\widetilde J$, and
$\mathrm{up}(J)=1$ if $J$ is the upper member of its pair, $\mathrm{up}(J)=0$ if it is the lower one.
In particular, if $u,u'$ are upper members of the same degree, then $d(\widetilde J_u)=d(\widetilde
J_{u'})$. Moreover, $\widetilde J_u$ determines the upper member $u$.
\end{lemma}

\begin{proof}
Write $\beta=\beta(J)$ and $\gamma=\gamma(J)$. By Lemma \ref{lem:betapair}, $\mathrm{up}(J)=[\beta<\gamma]$.

We have to show $Z(\widetilde J)-Z(J)=\mathrm{up}(J)-m$. Let $r'$ be the index $\neq j$ with
$x_{r'}=-\lambda$, if it exists: $r'=r$ if $\lambda\geq1$, $r'=j-1$ if $\lambda=\tfrac12$, and $r'$
does not exist if $\lambda=0$. As in Step 1 of the proof of Lemma \ref{lem:degpaired}, since $\alpha_j$
decreases by $1$ and $\omega_{\iota(j)}$ increases by $1$, the only comparisons which change are
$\alpha_j\leq \alpha_{j-1}$ and $\omega_{\iota(j-1)}\leq \omega_{\iota(j)}$ (from false to true), $\omega_{\iota(r')}\leq \alpha_j$
($\omega_{\iota(r')}=\lambda$) and $\omega_{\iota(j)}\leq \alpha_{r'}$ ($\alpha_{r'}=-\lambda$) (from true to false), and
$\omega_{\iota(j)}\leq \alpha_j$, which only occurs in the excluded pair $(j,\iota(j))$. So the terms of
(\ref{eq:Zpi}) which can change are those of the pairs $(j,j-1)$ and $(\iota(j-1),\iota(j))$, which are
exchanged by the symmetry $h(a,b)=h(b^\vee,a^\vee)$, and $(j,\iota(r'))$ and $(r',\iota(j))$, also exchanged
by the symmetry. (If $\iota(j)=j$, $\omega_j$ moves too, which is the same rule.)

\emph{The pairs $(j,\iota(r'))$, $(r',\iota(j))$.} Unless $\iota(r')\in\{j-1,j\}$, the string $s_{\iota(r')}$ has
bottom $\lambda$, hence top $>\lambda$, and $h(s_j,s_{\iota(r')})=[\beta\leq\lambda]=1$, since $\beta\leq
\lambda$ when $r^\lambda(J)\neq0$. For $\widetilde J$, the top of $\widetilde s_j$ is $\lambda-1<\lambda$, so
the term is $0$. These two pairs contribute $-1=-m$ (and $0=-m$ if $r'$ does not exist). They are
excluded from (\ref{eq:Zpi}) if $(j,\iota(r'))$ is a pair $(i,\iota(i))$, i.e.\ $r'=j$, which does not
happen.

\emph{The pairs $(j,j-1)$, $(\iota(j-1),\iota(j))$.} For $J$, $h(s_j,s_{j-1})=0$ since $\alpha_j>\alpha_{j-1}$. For
$\widetilde J$, $h(\widetilde s_j,\widetilde s_{j-1})=[\widetilde\beta\leq\gamma\leq\lambda-1]$, where
$\widetilde\beta$ is the bottom of $\widetilde s_j$ ($\widetilde\beta=\beta$, except $\widetilde\beta=
-\lambda+1$ if $\iota(j)=j$). If these pairs are distinct and not excluded, they contribute
$[\widetilde\beta\leq\gamma\leq\lambda-1]$, which is $[\beta<\gamma]=\mathrm{up}(J)$ in the following
cases: $\iota(j),\iota(j-1)\notin\{j-1,j\}$ and $\iota(j-1)\neq r'$ (then $\gamma\leq\lambda-1$ by Lemma
\ref{lem:betagamma}(a)); $\iota(j)=j$ (then $\widetilde\beta=-\lambda+1$, $\beta=-\lambda$,
and $\gamma\notin\{-\lambda,-\lambda+1\}$ when $\iota(j-1)\neq j-1$); $\iota(j-1)=j-1$ (then
$\gamma=-\lambda+1$).

\emph{Special cases.}
\begin{itemize}
\item $\iota(r')=j$ ($\lambda\geq1$): $s_j=\langle\lambda\rangle$ and $s_r=\langle-\lambda\rangle$, which
become empty in $\widetilde J$. Here $\beta=\lambda>\gamma$, so $\mathrm{up}(J)=0$. The changing terms are
$h(s_j,s_j)$ and $h(s_r,s_r)$ (the pairs $(j,\iota(r'))$ and $(r',\iota(j))$), from $1$ to $0$, and the
pairs $(j,j-1)$, $(\iota(j-1),r)$, which vanish on both sides. Hence $Z(\widetilde J)-Z(J)=-1=-m$.
\item $\iota(r')=j-1$ ($\lambda\geq1$): $s_{j-1}=\langle\lambda-1,\lambda\rangle$ and $s_r$ are empty, so
$\gamma=\lambda>\beta$ and $\mathrm{up}(J)=1$. All the changing terms involve $s_{j-1}$ or $s_r$, hence
vanish, and $Z(\widetilde J)-Z(J)=0=1-m$.
\item $\iota(j)=j-1$ (lower member in case (ii) of Proposition \ref{prop:cat48}): the pairs $(j,j-1)$, $(j-1,j)$ are excluded, and
$\mathrm{up}(J)=0$. If $\lambda\geq1$, only the pairs $(j,\iota(r))$, $(r,j-1)$ change, and
$Z(\widetilde J)-Z(J)=-1=-m$. If $\lambda=\tfrac12$, the strings are $s_j=\langle\tfrac12\rangle$,
$s_{j-1}=\langle-\tfrac12\rangle$, which become empty; only $h(s_j,s_j)=h(s_{j-1},s_{j-1})=1$ change, and
$Z(\widetilde J)-Z(J)=-1=-m$.
\item $\iota(j)=j$ and $\iota(j-1)=j-1$ (upper member in case (ii) of Proposition \ref{prop:cat48}): $\mathrm{up}(J)=1$, and for $\lambda\geq1$
the pairs $(j,j-1)$ and $(j-1,j)$ go from $0$ to $1$ (here $\widetilde s_j=\widetilde s_{j-1}=\langle
\lambda-1,\ldots,-\lambda+1\rangle$), and the pairs $(j,\iota(r))$, $(r,j)$ contribute $-1$:
$Z(\widetilde J)-Z(J)=1-m$. For $\lambda=\tfrac12$, $s_{j-1}$ and $\widetilde s_j$ are empty (blocks
$S_0$), $s_j=\langle\tfrac12,-\tfrac12\rangle$ meets no other string, and $Z(\widetilde J)-Z(J)=0=1-m$.
\end{itemize}
Finally, $\widetilde{\mathfrak g}_0$ only depends on $(\caL,\lambda)$, and $\deg J_u=\deg J_{u'}$ means
$d(J_u)=d(J_{u'})$, whence the assertion on the upper members of the same degree. Moreover, of the two strings
of $\widetilde J_u$ with top $\lambda-1$, the one coming from the string of top $\lambda$ of $u$ has the
smaller bottom, since $\beta(J_u)<\gamma(J_u)$; so $\widetilde J_u$ determines $u$.
\end{proof}

\subsection{The coefficients of Atobe's formula}\label{ssec:signs}
We now prove Theorem \ref{conj:geom}(1), by induction on the rank and on the degree. Recall that
$c_J$ is the coefficient of $[J]$ in Atobe's formula, after projection to the cuspidal support of
$\pi(\caL)$, and that $\Phi(\caL)$ is the set of $J$ with $c_J\neq0$.

\begin{lemma}\label{lem:coefunpaired} Let $\mu$ be an unpaired exponent with
$D_{\rho|\cdot|^\mu}(\pi(\caL))\neq0$, and $J$ a costandard module with $r^\mu(J)\neq0$,
$r^\mu(J)\simeq\rho|\cdot|^\mu\otimes J^\mu$. Then $c_J(\caL)=c_{J^\mu}(\caL_\mu)$, and $\deg J=\deg J^\mu$ when
these coefficients are non-zero.
\end{lemma}

\begin{proof}
By Proposition \ref{prop:cat49} and Lemma \ref{lem:small}, $\frS(\caL_\mu)=\frS(\caL)$
and, for every $(\sigma,\delta)$,
\[r^\mu(J_{\sigma,\delta}(\caL))\simeq\rho|\cdot|^\mu\otimes J_{\sigma,\delta}(\caL_\mu),\]
the right-hand side being $0$ for the summands killed by $r^\mu$. The module $J$ is recovered
from $J^\mu$, by putting back the value $\mu$ at the position of $\mu-1$ (since $\mu-1\notin X$). So the pairs
$(\sigma,\delta)$ with $J_{\sigma,\delta}(\caL)=J$ are those with $J_{\sigma,\delta}(\caL_\mu)=J^\mu$, with
the same signs $\mathrm{sgn}(\sigma)$, and $c_J(\caL)=c_{J^\mu}(\caL_\mu)$. The degrees are Lemma
\ref{lem:degunpaired}.
\end{proof}

\begin{lemma}\label{lem:coefpaired} Let $\lambda$ be a paired exponent, and $J$ a costandard module with
$r^\lambda(J)\neq0$. Let $J^\natural$ be its partner (\S\ref{ssec:paired}). If $J$ and $J^\natural$ do not have the same
strings, then $c_J+c_{J^\natural}=0$; if they have the same strings, then $c_J=0$. In particular, if
$J\in\Phi(\caL)$, then $J$ and $J^\natural$ do not have the same strings.
\end{lemma}

\begin{proof}
Since $\lambda-1\in X$, $D_{\rho|\cdot|^\lambda}(\pi(\caL))=0$ (\cite{Ato24}, Lemma 4.9), so
$\sum_Jc_J[r^\lambda(J)]=0$. The summands killed by $r^\lambda$ do not contribute, and the other ones are paired by
Atobe's involution, which on the modules is $J\mapsto J^\natural$, with $r^\lambda(J_v)\simeq
r^\lambda(J_u)=\rho|\cdot|^\lambda\otimes\widetilde J_u$ for $\{J_v,J_u\}=\{J,J^\natural\}$, $J_u$ the
upper member ($\beta(J_u)<\gamma(J_u)$)
(Proposition \ref{prop:cat48}, Lemma \ref{lem:smallpaired}; in case (ii) of Proposition \ref{prop:cat48},
$r^\lambda(J_v)\simeq r^\lambda(J_u^+)\oplus r^\lambda(J_u^-)$). The modules $\widetilde J_u$ are
costandard modules (proof of Lemma \ref{lem:matchediso}), and $\widetilde J_u$ determines $u$ (Lemma
\ref{lem:degtilde}). Since the costandard modules are linearly independent in the
Grothendieck group, the coefficient of $[\rho|\cdot|^\lambda\otimes\widetilde J_u]$, which is
$c_{J_v}+c_{J_u}$, vanishes. If $J$ and $J^\natural$ have the same strings, then $\lambda-1$ and $\lambda$ are
blocks in both, with the same sign since $r^\lambda(J)\neq0$, so that $J^\natural\simeq J$; then $\widetilde J$
has the block $\rho\boxtimes S_{2\lambda-1}$ twice, and $J$ is the only costandard module whose derivative is
$\rho|\cdot|^\lambda\otimes\widetilde J$ (the value $\lambda$ is put back at one of these two blocks), so that
the coefficient of $[\rho|\cdot|^\lambda\otimes\widetilde J]$ is $c_J$, which vanishes. (Pairs of permutations giving isomorphic summands contribute to $c_J$ with
opposite signs, and do not affect the argument, which only involves the modules.)
\end{proof}

\begin{lemma}\label{lem:lowersocle} Let $\lambda>0$ be a paired exponent, $J_v$ a costandard module
with $r^\lambda(J_v)\neq0$ which is the lower member of a pair $\{J_v,J_v^\natural\}$ whose members do not
have the same strings. Then $r^\lambda(L(J_v))=0$.
\end{lemma}

\begin{proof}
Let $J_u=J_v^\natural$ (or $J_u^\pm$), the upper member; by Lemma \ref{lem:betapair}, it is the upper member in the sense of Proposition \ref{thm:SAM}. Let $J=J_v$ or
$J_u$, and $L=L(J)$. The exponent $\lambda$ is not an element of the set
$\widetilde X$ attached to $\widetilde J$ ($\lambda$ has been replaced by $\lambda-1$), so
$r^\lambda(\widetilde J)=0$ by Lemma \ref{lem:onepiece}, and the
$\rho|\cdot|^\lambda\otimes\rho|\cdot|^\lambda$-part of the Jacquet module of $J$, hence of $L$,
vanishes. Since $\rho|\cdot|^\lambda$ is not self-dual, if $r^\lambda(L)\neq0$, then $r^\lambda(L)$ is
the highest derivative $\rho|\cdot|^\lambda\otimes D_{\rho|\cdot|^\lambda}(L)$, which is
irreducible, and $L$ is the unique irreducible subrepresentation of $\rho|\cdot|^\lambda\rtimes
D_{\rho|\cdot|^\lambda}(L)$ (\cite{AM23}, \S3.2, and the references there).

Assume $r^\lambda(L(J_v))\neq0$. It is an irreducible submodule of $r^\lambda(J_v)\simeq
\rho|\cdot|^\lambda\otimes\widetilde J_u$ (or of $\rho|\cdot|^\lambda\otimes(\widetilde J_u^+\oplus
\widetilde J_u^-)$), hence equal to $\rho|\cdot|^\lambda\otimes L(\widetilde J_u)$ (or to
$\rho|\cdot|^\lambda\otimes L(\widetilde J_u^\pm)$ for one sign), the costandard modules
$\widetilde J_u^{(\pm)}$ having a simple socle. By Proposition \ref{thm:SAM}(b), $r^\lambda(L(J_u))\neq0$,
and the same argument gives $r^\lambda(L(J_u))\simeq\rho|\cdot|^\lambda\otimes L(\widetilde J_u)$. So
$L(J_v)$ and $L(J_u)$ are both the unique irreducible subrepresentation of
$\rho|\cdot|^\lambda\rtimes L(\widetilde J_u)$, hence $L(J_v)\simeq L(J_u)$ and $J_v\simeq J_u$, a
contradiction.
\end{proof}

\begin{lemma}\label{lem:lowersocle0} Let $\lambda=0$ be a paired exponent ($\rho$ of integral type, with $-1\in X$),
and $J_v$ a costandard module with $r^0(J_v)\neq0$ which is the lower member of $\{J_v,J_v^\natural\}$.
Then $r^0(L(J_v))=0$.
\end{lemma}

\begin{proof}
By Lemma \ref{lem:smallpaired}, $0$ and $-1$ are the tops of two segments $\Delta_a=\Delta[0,-e]$ and
$\Delta_b=\Delta[-1,-e']$ of $J_v$, with $e,e'\in X$ placed at positions of R, $e>0$ and $e'>-1$, so
$e,e'\geq1$. In the lower member the top $0$ has the larger bottom (Lemma \ref{lem:betapair}), so
$e<e'$. Hence $\Delta_a$ and $\Delta_b$ are linked: their union is $\Delta[0,-e']$ and their intersection
$\Delta[-1,-e]$ is non-empty.

Reorder the segments of $J_v$ by non-decreasing tops. If two of them are linked, the one with the
smaller top also has the smaller bottom, hence the smaller exponent, so the relative order of linked
segments is unchanged, and the two orderings are related by transpositions of adjacent non-linked
segments, which give isomorphic products (\cite{Zel80}, Theorem 9.7). No other segment has top $0$ or
$-1$, so $J_v\simeq\cdots\times\Delta_b\times\Delta_a\times\cdots\rtimes T$. Let $\tau$ be the unique irreducible
subrepresentation of $\Delta_b\times\Delta_a$ (\cite{Zel80}). Then $J'=\cdots\times\tau\times\cdots
\rtimes T$ is a non-zero subrepresentation of $J_v$, so it contains the simple socle $L(J_v)$, and
$r^0(L(J_v))\subset r^0(J')$.

We compute $r^0(J')$ with the geometric lemma, as in Lemma \ref{lem:onepiece}. A subquotient of
$r^0(J')$ comes from a Jacquet module of one of the factors whose segment begins with $\rho|\cdot|^0$, from
one whose segment ends with $\rho|\cdot|^{-0}$ (through the contragredient), or from
$D_{\rho}(T)$. The only element of $X$ equal to $0$ is a top, so no other segment begins or ends at
$0$, and $\rho\boxtimes S_1\not\subset\phi_T$. For a representation of $\mathbf{GL}_{md}(F)$ whose
cuspidal support lies in $\rho|\cdot|^{\mathbb Z}$, let $r_{\min}$ be its Jacquet module with respect to the
standard parabolic subgroup with Levi factor $\mathbf{GL}_d(F)^m$; we record its constituents
$\rho|\cdot|^{a_1}\otimes\cdots\otimes\rho|\cdot|^{a_m}$ as words $(a_1,\ldots,a_m)$. The words of $r_{\min}(\tau)$ end with $-e$ or $-e'$,
not with $0$. It remains to see that no word of $r_{\min}(\tau)$ begins with $0$. By \cite{Zel80},
Theorem 9.7, $\Delta_b\times\Delta_a$ has length two, with constituents $\tau$ and
$\Delta[0,-e']\times\Delta[-1,-e]$. In $r_{\min}(\Delta_b\times\Delta_a)$, the words beginning with $0$
are the words $(0)\cdot w$, $w$ a shuffle of $(-1,\ldots,-e')$ and $(-1,\ldots,-e)$; the same holds for
$r_{\min}(\Delta[0,-e']\times\Delta[-1,-e])$, with the same multiplicities. Hence no word of
$r_{\min}(\tau)$ begins with $0$, $r^0(J')=0$, and $r^0(L(J_v))=0$.
\end{proof}

\begin{lemma}\label{lem:reach} Let $J\in\Phi(\caL)$, $J\neq J_0$. Then one of the following holds:
\begin{enumerate}
\item there are a self-dual $\rho$ and an unpaired exponent $\mu$ of $X_\rho$ with
$D_{\rho|\cdot|^\mu}(\pi(\caL))\neq0$ and $r^{\rho,\mu}(J)\neq0$;
\item there are a self-dual $\rho$ and a paired exponent $\lambda$ of $X_\rho$ with $r^{\rho,\lambda}(J)\neq0$, such
that $J$ is the upper member of $\{J,J^\natural\}$, i.e.\ $\beta(J)<\gamma(J)$ (Lemma \ref{lem:betapair}).
\end{enumerate}
\end{lemma}

\begin{proof}
Otherwise, we show that $r^{\rho,\mu}(L(J))=0$ for every self-dual $\rho$ and every $\mu\in\mathbb R$, which
contradicts the criterion of Section \ref{sec:step1}. We fix $\rho$ and write $X=X_\rho$. If $\mu\notin X$, then $r^\mu(J)=0$ (Lemma \ref{lem:onepiece}). If $\mu$ is unpaired
with $D_{\rho|\cdot|^\mu}(\pi(\caL))\neq0$, then $r^\mu(J)=0$ by assumption. If $\mu$ is unpaired with
$D_{\rho|\cdot|^\mu}(\pi(\caL))=0$, all the $r^\mu(J)$ vanish (Lemma \ref{lem:ato48}). If $\mu$
is paired, $J$ is either killed by $r^\mu$, or the lower member of $\{J,J^\natural\}$, whose members do not
have the same strings (Lemma \ref{lem:coefpaired}), and $r^\mu(L(J))=0$ by
Lemma \ref{lem:lowersocle} ($\mu>0$) or Lemma \ref{lem:lowersocle0} ($\mu=0$).
\end{proof}

\begin{thm}\label{thm:signs} Theorem \ref{conj:geom}(1) holds: $c_J=(-1)^{\deg J}$ for every
$J\in\Phi(\caL)$.
\end{thm}

\begin{proof}
By induction on the rank, and, for a given $\caL$, on $\deg J$. For $\deg J=0$, $J=J_0$ (Proposition
\ref{prop:closure}) and $c_{J_0}=1$. Let $J\in\Phi(\caL)$ with $\deg J>0$. In case (1) of Lemma
\ref{lem:reach}, $c_J=c_{J^\mu}(\caL_\mu)=(-1)^{\deg J^\mu}=(-1)^{\deg J}$, by Lemma \ref{lem:coefunpaired} and
the induction on the rank. In case (2), let $J_v=J^\natural$; then $c_{J_v}=-c_J\neq0$ by Lemma
\ref{lem:coefpaired}, so $J_v\in\Phi(\caL)$, with $\deg J_v=\deg J-1$ (Lemma \ref{lem:degpaired}); by the
induction on the degree, $c_J=-(-1)^{\deg J-1}=(-1)^{\deg J}$.
\end{proof}

\noindent In particular, all the coefficients of Atobe's formula are $\pm1$, as expected in
\cite{Ato24}, Remark 4.6.

\section{A Johnson-type construction}\label{sec:johnson}

\subsection{Successive cokernels}\label{ssec:cokernels}
Fix a ladder datum $\caL$, and put $\caJ^k=\bigoplus_{J\in\Phi(\caL),\,\deg J=k}J$. For an irreducible
representation $L'=L(J')$ with the infinitesimal character of $\pi(\caL)$, put $\deg L'=d(J')-d(J_0)$.

\begin{defi} A \emph{Johnson sequence} for $\caL$ consists of representations $Y_k$ and maps
$f_k:Y_k\to\caJ^k$, $k\geq0$, such that $Y_0=\pi(\caL)$, $f_0$ is the inclusion $\pi(\caL)\subset J_0$,
$Y_{k+1}=\mathrm{coker}(f_k)$, and every component $Y_k\to J$ of $f_k$ ($J\in\Phi(\caL)$, $\deg J=k$) is non-zero.
\end{defi}

\noindent If the maps $f_k$ are injective, the sequence
$0\to\pi(\caL)\to\caJ^0\to\caJ^1\to\cdots$ obtained by composing $\caJ^{k-1}\to Y_k\xrightarrow{f_k}\caJ^k$ is
exact at every $\caJ^k$, and exactness at the end follows from the next lemma. So Theorem
\ref{conj:geom}(2) amounts to the existence of a Johnson sequence with injective maps. No differential
has to be constructed between the costandard modules, and there is no sign problem.

\begin{lemma}\label{lem:bookkeeping} Assume Theorem \ref{conj:geom}(1), and that
$f_0,\ldots,f_{k-1}$ are injective. Then
\[[Y_k]=\sum_{i\geq k}(-1)^{i-k}[\caJ^i]\]
in the Grothendieck group. Every constituent of $Y_k$ has degree $\geq k$, and the constituents of
degree $k$ are the $L(J)$, $J\in\Phi(\caL)$ of degree $k$, each with multiplicity one. In particular $Y_{K+1}=0$ for $K$ the
largest degree.
\end{lemma}

\begin{proof}
The exact sequence $0\to\pi(\caL)\to\caJ^0\to\cdots\to\caJ^{k-1}\to Y_k\to0$ gives
$[Y_k]=\sum_{i<k}(-1)^{k-1-i}[\caJ^i]+(-1)^k[\pi(\caL)]$. Atobe's formula and Theorem
\ref{conj:geom}(1) give $[\pi(\caL)]=\sum_i(-1)^i[\caJ^i]$, whence the formula. By Proposition
\ref{prop:closure}, the constituents of $J\in\Phi(\caL)$ of degree $i$ have degree $\geq i$, and $L(J)$ is the only one of
degree $i$. The constituents of $Y_k$ occur with positive coefficients on the right-hand side, so they
have degree $\geq k$, and those of degree $k$ come from the socles of the $J\in\Phi(\caL)$ of degree $k$, once each.
\end{proof}

\subsection{Existence of the maps}

\noindent We need the analogue of \cite{Joh84}, Theorem III.1, and of the
injectivity of the dual Verma modules in the truncated subcategories of the category $\mathcal{O}$.
Let $\mathcal{C}_s$ be the category of representations of finite length on which the Bernstein centre
acts by the character of $\pi(\caL)$ (not only generalized), and $\mathcal{C}_{\geq m}\subset
\mathcal{C}_s$ the Serre subcategory of the representations whose constituents have degree $\geq m$ (degree as
in \S\ref{ssec:cokernels}).
Recall (\S\ref{ssec:std}) that a costandard module is of discrete type if its segments are pairwise
distinct and its tempered part is square integrable. All the $J\in\Phi(\caL)$ are of discrete type: the blocks of $\phi$ in $I_\sigma(\caL)$ are
attached to distinct points of $X$, and the segments are distinct.
We shall prove in \S\S\ref{ssec:inj} and \ref{ssec:clifford} (Theorem \ref{thm:injprime}, Lemma
\ref{lem:extcs} and Corollary \ref{cor:selfext}) the following property:
\begin{itemize}
\item[(Inj)] every $J\in\Phi(\caL)$ is injective in $\mathcal{C}_{\geq\deg J}$.
\end{itemize}
Part of the argument applies to all the costandard modules of discrete type (Theorem
\ref{thm:injprime}), but the vanishing of the self-extensions of $L(J)$ by $J$ uses $J\in\Phi(\caL)$
(Corollary \ref{cor:selfext}). Without the discrete type, the injectivity fails (Remark \ref{rem:injcex}).

\begin{lemma}\label{lem:existence} Assume (Inj) for $\caL$, Theorem \ref{conj:geom}(1), and that
$f_0,\ldots,f_{k-1}$ are injective. Then $\dim\mathrm{Hom}(Y_k,J)=1$ for every $J\in\Phi(\caL)$ of degree $k$. In
particular the map $f_k$ exists, and is unique up to the action of $\prod_{J\in\Phi(\caL),\,\deg J=k}\mathbb{C}^\times$.
\end{lemma}

\begin{proof}
The Bernstein centre acts on each costandard module by a character, so $Y_k$, a quotient of $\caJ^{k-1}$, lies in
$\mathcal{C}_s$, and in $\mathcal{C}_{\geq k}$ by Lemma \ref{lem:bookkeeping}. By (Inj),
$\mathrm{Hom}(\cdot,J)$ is exact on $\mathcal{C}_{\geq k}$. Since $\mathrm{Hom}(L',J)=0$ for $L'\neq L(J)$,
$\dim\mathrm{Hom}(Y_k,J)=[Y_k:L(J)]=1$.
\end{proof}

\subsection{Proof of (Inj)}\label{ssec:inj}
Since $\mathcal{C}_{\geq m}$ is a Serre subcategory, (Inj) says that, for $J\in\Phi(\caL)$,
\[\mathrm{Ext}^1_{\mathcal{C}_s}(L',J)=0\quad\text{for every irreducible }L'\text{ with }\deg L'\geq\deg J.\] We separate
the cases $L'\neq L(J)$ and $L'=L(J)$.

\begin{lemma}\label{lem:extcs} Let $J$ be a costandard module and $L'$ an irreducible representation
with the infinitesimal character of $J$. The natural map $\mathrm{Ext}^1_{\mathcal{C}_s}(L',J)\to
\mathrm{Ext}^1_G(L',J)$ is injective, and it is bijective if $L'\neq L(J)$.
\end{lemma}

\begin{proof}
An extension in $\mathcal{C}_s$ which splits in the category of all smooth representations splits in
$\mathcal{C}_s$, hence injectivity. Let $0\to J\to E\to L'\to0$ be an extension of smooth representations.
The Bernstein centre acts on $E$ through a finite quotient, with generalized character $\chi$. For $z$ in
the centre, $z-\chi(z)$ kills $J$ and $L'$, so it acts on $E$ through a $G$-map $L'\to J$. If
$L'\neq L(J)$, this map is zero, since $L(J)$ is the socle of $J$. Hence $E$ lies in $\mathcal{C}_s$.
\end{proof}

\noindent So, for $L'\neq L(J)$, (Inj) is a statement about the category of all smooth
representations:
\begin{itemize}
\item[(Inj$_G$)] if $J$ is of discrete type, $L'\neq L(J)$ and $\mathrm{Ext}^1_G(L',J)\neq0$, then
$\deg L'<\deg J$.
\end{itemize}
We deduce it from S.~Kato's structure theorem for extension algebras \cite{Kat17}. Let
$\mathcal{C}$ be the category of representations of finite length in the Bernstein block of $\pi(\caL)$
with the infinitesimal character $\lambda_\caL$ (generalized). Let $\mathcal K$ be Lusztig's induced
complex attached to the cuspidal datum of this block: the parabolic induction to
$\widehat{\mathfrak g}_{-1}$ of a cuspidal local system on a Levi subgroup (\cite{Lus95}, \cite{Lus02};
$\mathcal K_{N,\sigma,-r}$ in \cite{Sol25}). By the decomposition theorem, $\mathcal K$ is a direct sum of
shifts of the complexes $\mathrm{IC}_\lambda$, where $\lambda=(O,\xi)$ runs over the pairs of the series
of the block (those which parametrize the representations in the block), and by \cite{Sol25}, Theorems
5.2(b) and 5.4, $\mathcal{C}$ is equivalent to the category of finite-dimensional modules, on which
$H^{>0}_{G_0}(\mathrm{pt})$ acts nilpotently, over $\mathrm{Ext}^\bullet_{G_0}(\mathcal K,\mathcal K)$,
which is Morita equivalent to the extension algebra
\[A=\mathrm{Ext}^\bullet_{G_0}(\mathcal{L},\mathcal{L}),\qquad \mathcal{L}=\bigoplus_\lambda
L_\lambda\boxtimes\mathrm{IC}_\lambda,\]
of the $G_0$-equivariant derived category of $\widehat{\mathfrak g}_{-1}$. The irreducible $L_\lambda$
corresponds to $L(J)$ for $J$ with parameter $\lambda$, and the costandard module $J$ to the graded dual
$K_\lambda^*$ of Kato's standard module $K_\lambda=\mathrm{Hom}_{\Gamma_\lambda}(\xi_\lambda,H^\bullet
i_N^!\mathcal L)$ ($N\in O_\lambda$, $\Gamma_\lambda$ the component group of $Z_{G_0}(N)$), which has simple
head $L_\lambda$ (\cite{Kat17}, \S1; \cite{Sol25}, Theorem 3.2(a)). We write $\lambda\preceq\mu$ if
$O_\lambda\subset\overline{O_\mu}$, and $\lambda\sim\mu$ if $O_\lambda=O_\mu$.
Kato also introduces $\widetilde K_\lambda=P_\lambda/\sum_{\mu\prec\lambda}Ae_\mu P_\lambda$, the
largest quotient of the projective cover $P_\lambda$ of $L_\lambda$ with no constituent $L_\mu$,
$\mu\prec\lambda$.

\begin{lemma}\label{lem:K} Assume that $G_0$ is connected. Then the conditions $(\spadesuit)$ and
$(\clubsuit)$ of \cite{Kat17}, \S1, hold for $(G_0,\widehat{\mathfrak g}_{-1})$ and the set of all the pairs $(O,\xi)$.
\end{lemma}

\begin{proof}
$(\spadesuit)$ is the finiteness of the number of $G_0$-orbits in $\widehat{\mathfrak g}_{-1}$, with
the labels $(O,\xi)$ (\cite{Kat17}, Condition 1.1; the stabilizers need not be connected). The pair
$(G_0,\widehat{\mathfrak g}_{-1})$ is the degree-zero group and the degree $-1$ part of the
$\bbZ$-grading of the Lie algebra of the connected reductive group $H=\big(\prod_\rho\mathrm{I}(M_\rho)
\big)^\circ$ defined by the eigenvalues of $\mathrm{Ad}(s)$, $s=(s_\rho)_\rho$.
\begin{itemize}
\item $(\clubsuit)_2$, the pointwise purity of the $\mathrm{IC}_\lambda$: by Lusztig
(\cite{Lus95b}) for $\bbZ$-gradings, and by Lusztig and Yun in general (\cite{LY17}, 12.2(a) and
Theorem 14.10), for every irreducible equivariant local system on an orbit, the cohomology sheaves of
its intersection cohomology complex are pointwise pure (over an algebraic closure of a finite field,
after spreading out) and vanish in odd degrees. These results are stated for a semisimple simply
connected group with a $\bbZ/m$-grading, the degree-zero group being the connected subgroup with Lie
algebra the degree-zero part (\cite{LY16}, 0.1; \cite{LY17}, Introduction). We apply them to the simply
connected cover $\widetilde H$ of the derived group of $H$, with the $\bbZ/m$-grading induced by
$\mathrm{Ad}(s)$ for $m$ larger than twice the degrees which occur, so that the degree $-1$ part is
$\widehat{\mathfrak g}_{-1}$. Its degree-zero group $\widetilde H_0$ maps onto $G_0$ modulo the connected
centre of $H$, which acts trivially on $\widehat{\mathfrak g}_{-1}$ (both groups are connected, with the
same Lie algebra modulo the centre); so $\widetilde H_0$ and $G_0$ have the same orbits on
$\widehat{\mathfrak g}_{-1}$, every $G_0$-equivariant local system is $\widetilde H_0$-equivariant, and
the statement holds for the $\mathrm{IC}_\lambda$.
\item $(\clubsuit)_1$, the purity of the extension algebra of the $\mathrm{IC}_\lambda$: the
$G_0$-equivariant extension groups between two intersection cohomology complexes are computed, through
Verdier duality, by the compactly supported cohomology of their tensor product on the approximations
$G_0\backslash(\Gamma\times\widehat{\mathfrak g}_{-1})$ of the equivariant derived category
(\cite{BL94}). Filtering by the orbits, the pointwise purity in $(\clubsuit)_2$ and the purity of the
cohomology of the classifying spaces of the stabilizers show that these groups are pure, the spectral
sequences degenerating for weight reasons: this is \cite{LY17}, Proposition 13.3, after \cite{Lus95b},
19.4 (see also \cite{BGS96}, Theorem 3.4.1 and the remark which follows it, for the non-equivariant
analogue).
\end{itemize}
\end{proof}

\noindent By \cite{LY16}, Theorem 0.6, the $G_0$-equivariant derived category of the nilpotent cone of
$\widehat{\mathfrak g}_{-1}$ (which is $\widehat{\mathfrak g}_{-1}$ itself) is the direct sum of blocks, with
no non-zero morphisms between objects of different blocks, and the summands of $\mathcal K$ are the simple
objects of one of them. So the extension algebra of all the $\mathrm{IC}_\lambda$ is the product of the
extension algebras of the blocks, one of which is the algebra $A$ above; Kato's results for the set of
all the pairs restrict to $A$ (the modules $P_\lambda$, $\widetilde K_\lambda$, $K_\lambda$, $L_\lambda$ attached to
a pair $\lambda$ of the series are modules over this factor, $\widetilde K_\lambda$ and $K_\lambda$ being
quotients of $P_\lambda$).

\noindent By Lemma \ref{lem:K}, \cite{Kat17}, Theorem 1.8, applies and gives: $[\widetilde K_\lambda:L_\nu]\neq0$
only if $\lambda\preceq\nu$; $[K_\lambda:L_\nu]=\delta_{\lambda\nu}$ for $\lambda\sim\nu$; the kernel
$R_\lambda$ of $\widetilde K_\lambda\to K_\lambda$ has a separable filtration whose graded pieces are
shifts of the $K_\nu$, $\nu\sim\lambda$; and, for $\nu\sim\lambda$, the multiplicity
$[\widetilde K_\lambda:L_\nu]$ is the graded dimension of
$\mathrm{Hom}_{\Gamma_\lambda}(\xi_\lambda,\xi_\nu\otimes H^\bullet_{Z_\lambda^\circ}(\mathrm{pt}))$ (\cite{Kat17},
Theorem 1.8(4)), where
$Z_\lambda$ is the stabilizer of a point of $O_\lambda$ and $\Gamma_\lambda$ its component
group.

\begin{lemma}\label{lem:compgroup} Let $J$ be of discrete type, with parameter $\phi$, and
$N=N_J$. The reductive quotient of $Z_{G_0}(N)=Z_{\widehat G}(\lambda_\caL(W_F),N)$ is $Z_{\widehat G}(\phi)$, which is
commutative: $Z_{\widehat G}(\phi)^\circ$ is a torus $T\simeq(\mathbb{C}^\times)^k$, one factor for each
segment, and the component group acts trivially on $T$.
\end{lemma}

\begin{proof}
The first assertion is standard (Jacobson--Morozov). Write
\[\phi=\bigoplus_{i=1}^k(\rho|\cdot|^{a_i}\boxtimes S_{n_i}\oplus\rho|\cdot|^{-a_i}\boxtimes S_{n_i})
\oplus\phi_\tau,\]
with $a_i<0$ the exponent of $\Delta_i$ and $\phi_\tau$ the multiplicity free
parameter of $\tau$. The blocks $\rho|\cdot|^{\pm a_i}\boxtimes S_{n_i}$ are pairwise non-isomorphic,
since the segments are distinct and their exponents are non-zero, and they are not isomorphic to the
(unitary) blocks of $\phi_\tau$. So $Z_{\widehat G}(\phi)$ is the intersection with $\widehat G$ of a
product of groups $\mathrm{GL}_1$ (one for each pair of dual blocks) and $\mathrm{O}_1$ (one for each
block of $\phi_\tau$), acting by scalars on the isotypic components. It is commutative.
\end{proof}

\begin{thm}\label{thm:injprime} Let $J$ be a costandard module of discrete type and $L'$
an irreducible representation with $L'\neq L(J)$ and $\mathrm{Ext}^1_G(L',J)\neq0$. Then
$O_{L'}\subsetneq\overline{O_J}$. In particular (Inj$_G$) and, with Lemma \ref{lem:extcs} and Corollary
\ref{cor:selfext}, (Inj) hold for the $J\in\Phi(\caL)$.
\end{thm}

\begin{proof}
Assume first that $G_0$ is connected, so that Lemma \ref{lem:K} applies.
Let $\lambda$ be the parameter of $J$ and $\mu$ that of $L'$. Since $A$ is noetherian and positively
graded, and the modules involved are gradable, $\mathrm{Ext}^1$ in $\mathcal{C}$ is the graded
$\mathrm{Ext}^1_A$ of the corresponding graded modules, the grading being forgotten. Verdier duality gives an anti-involution of $A$ fixing the
idempotents $e_\nu$ (the local systems are self-dual, the component groups being elementary abelian
$2$-groups), hence a duality $M\mapsto M^*$ on finite-dimensional graded modules, with $L_\nu^*\simeq
L_\nu$ and $\mathrm{Ext}^i(M,N)\simeq\mathrm{Ext}^i(N^*,M^*)$. So
$\mathrm{Ext}^1(L_\mu,K_\lambda^*)\simeq\mathrm{Ext}^1(K_\lambda,L_\mu)$. Let $U_\lambda=\sum_{\nu\prec
\lambda}Ae_\nu P_\lambda$, so that $\widetilde K_\lambda=P_\lambda/U_\lambda$, and $R_\lambda$ as above.
The exact sequences give
\[\mathrm{Hom}(R_\lambda,L_\mu)\to\mathrm{Ext}^1(K_\lambda,L_\mu)\to\mathrm{Ext}^1(\widetilde
K_\lambda,L_\mu),\qquad \mathrm{Hom}(U_\lambda,L_\mu)\to\mathrm{Ext}^1(\widetilde K_\lambda,L_\mu)\to0.\]

\emph{Right-hand term.} $U_\lambda$ is a quotient of a sum of shifts of the $P_\nu$, $\nu\prec\lambda$,
so $\mathrm{Hom}(U_\lambda,L_\mu)\neq0$ implies $\mu\prec\lambda$.

\emph{Left-hand term.} Let $F_0=R_\lambda\supset F_1\supset\cdots$ be the filtration. A non-zero map
$R_\lambda\to L_\mu$ is non-zero on some $F_i$ and zero on $F_{i+1}$ (separability, $L_\mu$ being
finite-dimensional), so it factors through a shift of some $K_\nu$, $\nu\sim\lambda$, and $L_\mu$ is its
head $L_\nu$. Hence $\mu\sim\lambda$ and $L_\mu$ is a constituent of $\widetilde K_\lambda$. By Lemma
\ref{lem:compgroup}, $\Gamma_\lambda$ acts trivially on $H^\bullet_{Z_\lambda^\circ}(\mathrm{pt})$, so
$[\widetilde K_\lambda:L_\mu]$ is the graded dimension of $\mathrm{Hom}_{\Gamma_\lambda}(\xi_\lambda,
\xi_\mu)\otimes H^\bullet_{Z_\lambda^\circ}(\mathrm{pt})$, which vanishes for $\xi_\mu\neq\xi_\lambda$.
Hence $\mu=\lambda$, which is excluded.

So $\mathrm{Ext}^1(L',J)\neq0$ forces $\mu\prec\lambda$, that is, $O_{L'}\subsetneq\overline{O_J}$ when
$G_0$ is connected. The general case is deduced from it in \S\ref{ssec:clifford}.
\end{proof}

\noindent The proof also shows where the discrete type enters. For $\mu\sim\lambda$ the left-hand term
is controlled by the action of $\Gamma_\lambda$ on $H^2_{Z_\lambda^\circ}(\mathrm{pt})=\mathfrak{t}^*$. When it is non-trivial,
as for a symmetric segment (the block has multiplicity two, $\Gamma_\lambda$ contains the non-trivial component
of $\mathrm{O}_2$, which acts by inversion on $\mathrm{SO}_2$), one gets the extensions of Remark
\ref{rem:injcex}. For $\mu=\lambda$ one gets $\mathrm{Ext}^1(L(J),J)\simeq\mathfrak{t}^*=\mathfrak{a}^*$,
which is consistent with Proposition \ref{prop:selfext}.

\begin{prop}\label{prop:selfext} Let $J=\Delta_1\times\cdots\times\Delta_k\rtimes\tau$ be a costandard
module of discrete type, and $\mathfrak{a}=\mathbb{C}^k$ the space of unramified deformations
$\Delta_i\mapsto\Delta_i|\cdot|^{\kappa_i}$. Then $\mathrm{Ext}^1_G(L(J),J)\simeq\mathfrak{a}$, and
$\mathrm{Ext}^1_{\mathcal{C}_s}(L(J),J)$ is isomorphic to the space of $\kappa\in\mathfrak{a}$ such
that
\[\sum_{i=1}^k\kappa_i\,m_{\Delta_i}=0,\qquad m_\Delta(a)=[\rho|\cdot|^a\in\Delta]-
[\rho|\cdot|^{-a}\in\Delta]\quad(a>0).\]
\end{prop}

\begin{proof}
Group the segments with the same exponent: segments with the same centre are nested, hence not
linked, so the corresponding products are irreducible, and $J=I_P(\sigma)$ with
$P=MN$, $M=\mathbf{GL}_{N_1}\times\cdots\times\mathbf{GL}_{N_r}\times G_m$,
$\sigma=\pi_1|\cdot|^{e_1}\otimes\cdots\otimes\pi_r|\cdot|^{e_r}\otimes\tau$, where $\pi_c$ is the
(tempered, irreducible) product of the $\Delta_i|\cdot|^{-e_c}$ with $\Delta_i$ of exponent
$e_c$, and $e_1<\cdots<e_r<0$. This is the Langlands datum of $J$, in the subrepresentation
setting. By Frobenius reciprocity (the induction is exact and right adjoint to the exact functor
$r_P$), $\mathrm{Ext}^1_G(L(J),J)=\mathrm{Ext}^1_M(r_P(L(J)),\sigma)$. The inclusion $L(J)\subset J$ gives
a surjection $r_P(L(J))\to\sigma$; let $K$ be its kernel. By \cite{BJ08}, Corollary 5.4, $\sigma$ is the
only irreducible subquotient of $r_P(J)$ with the central character of $\sigma$, and it occurs with
multiplicity one. Since $r_P(L(J))\subset r_P(J)$, the centre of $M$ acts on $K$ and on $\sigma$ with
distinct generalized characters, so $\mathrm{Ext}^i_M(K,\sigma)=0$ for all $i$, and the long exact
sequence gives $\mathrm{Ext}^1_M(r_P(L(J)),\sigma)=\mathrm{Ext}^1_M(\sigma,\sigma)$.

By the K\"unneth formula, $\mathrm{Ext}^1_M(\sigma,\sigma)$ is the sum of the
$\mathrm{Ext}^1(\pi_c,\pi_c)$ and of $\mathrm{Ext}^1_{G_m}(\tau,\tau)$. By \cite{OS13}, Theorem 3
(\S7, the $p$-adic version of Theorem 5.2), for an irreducible tempered representation
$\pi$ induced from a discrete series datum $\xi$ (notation of \cite{OS13}), one has
$\mathrm{Ext}^1(\pi,\pi)=(E_\xi^*)^{R_\xi}$, with $E_\xi$ the space of unramified deformations of
$\xi$ and $R_\xi$ its $R$-group. For $\tau$ square
integrable, $E_\xi=0$. For $\pi_c$, the $R$-group is trivial ($\mathbf{GL}$), and $E_\xi$ is the space
of deformations $\Delta_i\mapsto\Delta_i|\cdot|^{\kappa_i}$ of the segments of
exponent $e_c$. Hence $\mathrm{Ext}^1_M(\sigma,\sigma)=\mathfrak{a}$, the class of $\kappa$ being
given by $I_P(\sigma_\kappa)$, where $\sigma_\kappa$ is obtained by replacing $\Delta_i$
by $\Delta_i|\cdot|^{\epsilon\kappa_i}$ (over $\mathbb{C}[\epsilon]/(\epsilon^2)$), and the
extension of $L(J)$ by $J$ is its pull-back along $L(J)\subset J$.

Let $s_\kappa$ be the infinitesimal character of $\sigma_\kappa$, in coordinates: the exponents
$x$ of the elements of $\Delta_i$ move to $x+\epsilon\kappa_i$, those of $\tau$ stay fixed. An element $z$ of
the Bernstein centre is a $W$-invariant function, and acts on $I_P(\sigma_\kappa)$ by
$z(s)+\epsilon\, dz_s(\dot s_\kappa)$. The differential at $s$ of the $W$-invariant functions
vanishes exactly on the span of the $(1-w)v$, $w\in W_s$, i.e.\ on the vectors on which every
$W_s$-invariant linear form vanishes. These forms are spanned by
$v\mapsto\sum_{s_j=a}v_j-\sum_{s_j=-a}v_j$, $a>0$, and on $\dot s_\kappa$ they give
$\sum_i\kappa_im_{\Delta_i}(a)$. So $z-\chi(z)$ acts on the pull-back as $dz_s(\dot s_\kappa)$ times the composite
$L(J)\to L(J)\subset J$ (the map $\epsilon$), and the extension lies in $\mathcal{C}_s$ if and only if
$dz_s(\dot s_\kappa)=0$ for all $z$, i.e.\ $\sum_i\kappa_im_{\Delta_i}=0$. By Lemma \ref{lem:extcs} (injectivity), this gives the statement.
\end{proof}

\begin{lemma}\label{lem:mindep} Let $J=\Delta_1\times\cdots\times\Delta_k\rtimes T$ be a summand of
$I_\sigma(\caL)$, $\sigma\in\frS(\caL)$, with $\Delta_j=\Delta[x_{\sigma(j)},-x_{\sigma(t-j+1)}]$,
$j\in J^+_\sigma$. Then the functions $m_{\Delta_j}$ are linearly independent.
\end{lemma}

\begin{proof}
Put $a=x_{\sigma(j)}$ and $e=x_{\sigma(t-j+1)}$. Since $j\in J^+_\sigma$, $\sigma(j)<\sigma(t-j+1)$, so
$a<e$; and $a\geq-e$ since the segment is non-empty (the empty segments are omitted from $J$). Hence $e>|a|$ or $a=-e<0$; in both cases $e>0$.
For $c>0$, $\rho|\cdot|^{-c}\in\Delta_j$ if and only if $-e\leq -c\leq a$, and
$\rho|\cdot|^{c}\in\Delta_j$ if and only if $-e\leq c\leq a$. So $m_{\Delta_j}(c)=0$ for $c>e$, and
$m_{\Delta_j}(e)=-1$, since $e>a$. The points $e_j=x_{\sigma(t-j+1)}$, $j\in J^+_\sigma$, are distinct,
because $\sigma$ is a permutation and the $x_i$ are distinct. Order the segments by decreasing $e_j$:
the matrix $(m_{\Delta_j}(e_{j'}))$ is triangular with diagonal $-1$.
\end{proof}

\begin{cor}\label{cor:selfext} For every $J\in\Phi(\caL)$, $\mathrm{Ext}^1_{\mathcal{C}_s}(L(J),J)=0$.
\end{cor}

\begin{proof}
Proposition \ref{prop:selfext} and Lemma \ref{lem:mindep}.
\end{proof}

\noindent The restriction to $\mathcal{C}_s$ is essential here: in the category of all smooth
representations, $L(J)$ always has the self-extensions coming from $\mathfrak{a}$.

\begin{rmk}\label{rem:injcex} The injectivity of (Inj) fails for costandard modules which are not of discrete type. For
$\mathbf{SO}(7)$, at the infinitesimal character with exponents $\tfrac32,\tfrac12,\tfrac12$, let
$T_\pm$ be the two tempered constituents of $\Delta[\tfrac12,-\tfrac12]\rtimes1$ of $\mathbf{SO}(5)$,
and $J_\pm=\rho|\cdot|^{-3/2}\rtimes T_\pm$. Then $J_+$ and $J_-$ have the same orbit (with distinct
local systems), and $\mathrm{Ext}^1_{\mathcal{C}_s}(L(J_\mp),J_\pm)\neq0$. The extensions come from the
deformation $\Delta[\tfrac12,-\tfrac12]|\cdot|^{\epsilon}$, for which $m=0$: it stays in $\mathcal{C}_s$
and mixes $T_+$ and $T_-$. This is the case $n=1$ of \cite{OS13}, Theorem 5.2: the
$R$-group $\mathbb{Z}/2$ acts by $-1$ on the one-dimensional space $E_\xi$ of unramified deformations, so
that $\mathrm{Ext}^1$ is non-zero between the two constituents, not on each of them. See also the
discussion after Theorem \ref{thm:injprime}.
\end{rmk}

\subsection{Disconnected centralizers}\label{ssec:clifford}
We complete the proof of Theorem \ref{thm:injprime} when $G_0$ is not connected, using the description
of $A$ and of its modules in terms of the identity component $G_0^\circ$ due to Aubert, Moussaoui and
Solleveld \cite{AMS18}.

\emph{The algebra.} Let $(M,C,q\mathcal E)$ be the cuspidal quasi-support of the series of the block
(\cite{Sol25}, \S1), $\mathcal E$ a cuspidal local system for $G_0^\circ$ contained in $q\mathcal E$, and
$A^\circ=\mathbb H(G_0^\circ,M^\circ,\mathcal E)$ the extension algebra of its series for $G_0^\circ$; it
is Morita equivalent to the algebra $A$ of \S\ref{ssec:inj} constructed for $G_0^\circ$ instead of $G_0$. By
\cite{Sol25}, Theorem 1.1, and \cite{AMS18}, Proposition 3.2 and (4) (see also \cite{Sol25}, (B.3)),
\[A\simeq A^\circ\rtimes\bbC[R,\natural],\]
where $R$ is a finite group acting on $A^\circ$ by graded automorphisms and $\bbC[R,\natural]$ is the
twisted group algebra of $R$ for a $2$-cocycle $\natural$. For a subgroup $R'\subset R$, put
$A_{R'}=A^\circ\rtimes\bbC[R',\natural]$. Since $A_{R'}$ is a free $A^\circ$-module of finite rank and the
invariants of a finite group are exact in characteristic $0$:
\begin{itemize}
\item[(a)] for graded $A_{R'}$-modules $X,Y$, the group $R'$ acts on $\mathrm{Ext}^i_{A^\circ}(X,Y)$ (the
cocycle cancels), and $\mathrm{Ext}^i_{A_{R'}}(X,Y)=\mathrm{Ext}^i_{A^\circ}(X,Y)^{R'}$;
\item[(b)] the induction $\mathrm{Ind}_{R'}^R=A\otimes_{A_{R'}}-$ is exact, and
$\mathrm{Ext}^i_A(\mathrm{Ind}_{R'}^RX,Y)=\mathrm{Ext}^i_{A_{R'}}(X,Y)$ (Frobenius reciprocity).
\end{itemize}

\emph{The modules.} Let $\lambda=(O,\xi)$ be a pair of the series and $N\in O$. Let $Z=Z_{G_0}(N)$,
$\Gamma_\lambda=\pi_0(Z)$ and $\Gamma'_\lambda=\pi_0(Z\cap G_0^\circ)\subset\Gamma_\lambda$; the groups $Z$ and $Z\cap
G_0^\circ$ have the same identity component $Z^\circ$. Choose a pair $\lambda^\circ=(O^\circ,\xi^\circ)$ of the
series of $\mathcal E$ for $G_0^\circ$, with $O^\circ=G_0^\circ N\subset O$ and $\xi^\circ$ an irreducible
constituent of $\xi|_{\Gamma'_\lambda}$, and let $R_\lambda\subset R$ be its stabilizer. By Clifford theory for
twisted crossed products (\cite{AMS18}, Lemmas 4.17 and 4.18 and the proof of Theorem 4.20), there is an
irreducible projective representation $(\tau_\lambda,V_\lambda)$ of $R_\lambda$ with cocycle $\natural$ such
that
\[K_\lambda\simeq\mathrm{Ind}_{R_\lambda}^R(V_\lambda\otimes K^\circ_{\lambda^\circ}),\qquad
L_\lambda\simeq\mathrm{Ind}_{R_\lambda}^R(V_\lambda\otimes L^\circ_{\lambda^\circ}),\]
where $R_\lambda$ acts on $K^\circ_{\lambda^\circ}$ and on its quotient $L^\circ_{\lambda^\circ}$ through the
geometric intertwiners of \cite{AMS18}, Lemma 4.17, and the pair $\lambda$ is determined by the
$R$-orbit of $\lambda^\circ$ and by $\tau_\lambda$. (In \cite{AMS18} this is stated after specialization at a
central character; the isomorphisms come from the decomposition $\mathcal P_y\simeq\mathcal P^\circ_y\times R$
of \cite{AMS18}, (17), which is equivariant, and hold for the graded modules as well.) In particular,
the restriction of $K_\lambda$ to $A^\circ$ is a sum of the modules $K^\circ_{\gamma\lambda^\circ}$, $\gamma\in R$,
attached to pairs on $G_0^\circ$-orbits contained in $O$, and similarly for $L_\lambda$.

\begin{proof}[End of the proof of Theorem \ref{thm:injprime}]
As in the connected case, it suffices to show that
\[\mathrm{Ext}^1_A(K_\lambda,L_\mu)\neq0\ \Longrightarrow\ O_\mu\subsetneq\overline{O_\lambda}.\] By (b),
\[\mathrm{Ext}^1_A(K_\lambda,L_\mu)=\mathrm{Ext}^1_{A_{R_\lambda}}(V_\lambda\otimes K^\circ_{\lambda^\circ},L_\mu),\]
so, by (a), $\mathrm{Ext}^1_{A^\circ}(K^\circ_{\lambda^\circ},L^\circ_{\mu'})\neq0$ for some pair $\mu'=(O'',\eta')$
with $L^\circ_{\mu'}$ a summand of the restriction of $L_\mu$, hence $O''\subset O_\mu$. The pair
$\lambda^\circ$ is of discrete type for $G_0^\circ$, since $\Gamma'_\lambda\subset\Gamma_\lambda$ acts trivially on
$Z^\circ$ (Lemma \ref{lem:compgroup}). The connected case gives $O''\subsetneq\overline{O^\circ}$ or
$\mu'=\lambda^\circ$.

If $O''\subsetneq\overline{O^\circ}$, then $\dim O_\mu=\dim O''<\dim O^\circ=\dim O_\lambda$, and $O_\mu$ meets the
$G_0$-stable closed set $\overline{O_\lambda}$, so $O_\mu\subsetneq\overline{O_\lambda}$.

Assume now that there is no such $\mu'$, so that $\mu'=\lambda^\circ$. Then $L^\circ_{\lambda^\circ}$ is a
summand of the restriction of $L_\mu$, so $O_\mu=O_\lambda$ and, by Clifford theory, we may take
$\mu^\circ=\lambda^\circ$, so that $L_\mu=\mathrm{Ind}_{R_\lambda}^R(V_\mu\otimes L^\circ_{\lambda^\circ})$ with
$\tau_\mu\not\simeq\tau_\lambda$, since $\mu\neq\lambda$. The restriction of $L_\mu$ to $A_{R_\lambda}$ is
$V_\mu\otimes L^\circ_{\lambda^\circ}$ plus a sum of modules on which $A^\circ$ acts through the modules
$L^\circ_{\gamma\lambda^\circ}$, $\gamma\notin R_\lambda$; these pairs $\gamma\lambda^\circ\neq\lambda^\circ$ are on
$G_0^\circ$-orbits of the same dimension as $O^\circ$, so $\mathrm{Ext}^1_{A^\circ}(K^\circ_{\lambda^\circ},
L^\circ_{\gamma\lambda^\circ})=0$ by the connected case. By (a),
\[\mathrm{Ext}^1_A(K_\lambda,L_\mu)=\big(\mathrm{Hom}(V_\lambda,V_\mu)\otimes E\big)^{R_\lambda},\qquad
E=\mathrm{Ext}^1_{A^\circ}(K^\circ_{\lambda^\circ},L^\circ_{\lambda^\circ}).\]
The proof of the connected case shows that $E$ is a quotient of $\mathrm{Hom}(R^\circ_{\lambda^\circ},
L^\circ_{\lambda^\circ})$, which is filtered with graded pieces
$\mathrm{Hom}_{\Gamma'_\lambda}(\xi^\circ,\xi^\circ\otimes H^{2k}_{Z^\circ}(\mathrm{pt}))=H^{2k}_{Z^\circ}(\mathrm{pt})$, $k\geq1$
(\cite{Kat17}, Theorem 1.8(4); $\Gamma'_\lambda$ acts trivially on $H^\bullet_{Z^\circ}(\mathrm{pt})=S(\mathfrak t^*)$,
$\mathfrak t=\mathrm{Lie}\,Z^\circ$). The group $R_\lambda$ acts on these pieces through its action on
$H^\bullet_{Z^\circ}(\mathrm{pt})$, which factors through $\Gamma_\lambda$ (\cite{AMS18}, Lemma 4.17), and is
therefore trivial by Lemma \ref{lem:compgroup}. So $R_\lambda$ acts trivially on $E$, and
\[\mathrm{Ext}^1_A(K_\lambda,L_\mu)=\mathrm{Hom}_{R_\lambda}(V_\lambda,V_\mu)\otimes E=0,\]
since $\tau_\lambda$ and $\tau_\mu$ are non-isomorphic irreducible projective representations with the same
cocycle. This contradicts the assumption, so the first case occurs.
\end{proof}

\noindent Two distinct pairs $\lambda\neq\mu$ on the same orbit may have the same restriction to $G_0^\circ$
(they differ by $\tau$); it is the triviality of the action of the whole group $\Gamma_\lambda$ on $\mathfrak t$
which separates them.

\subsection{Injectivity}
\begin{lemma}\label{lem:inj_unpaired} Let $\mu$ be an unpaired exponent with
$D_{\rho|\cdot|^\mu}(\pi(\caL))\neq0$, and $(Y_k,f_k)$ a Johnson sequence for $\caL$ with
$f_0,\ldots,f_{k-1}$ injective. Then $r^\mu(Y_i)\simeq\rho|\cdot|^\mu\otimes
Y'_i$ and $r^\mu(f_i)=1\otimes f'_i$ for a Johnson sequence $(Y'_i,f'_i)$ of $\caL_\mu$, for
$i\leq k$.
\end{lemma}

\begin{proof}
By Proposition \ref{prop:cat49} and Lemma \ref{lem:small}, $r^\mu$ kills the
summands $J$ with $r^\mu(J)=0$, and sends the other ones to the summands of $\caL_\mu$, summand by summand. So the
coefficients of Atobe's formula correspond, and $r^\mu(\caJ^i)=\rho|\cdot|^\mu\otimes\caJ'^i$ by Lemma
\ref{lem:degunpaired}. Since $r^\mu$ is exact, $r^\mu(\pi(\caL))=\rho|\cdot|^\mu\otimes
\pi(\caL_\mu)$ and $r^\mu(Y_{i+1})=\mathrm{coker}\,r^\mu(f_i)$. Finally, every component $Y_i\to J$
of $f_i$ is non-zero, its image contains $L(J)$, and $r^\mu(L(J))\neq0$ (Proposition \ref{thm:SAM}(a)). So the components of $r^\mu(f_i)$ are non-zero.
\end{proof}

\begin{prop}\label{prop:johnson} Assume Theorem \ref{conj:geom}(1) and (Inj) for $\caL$ and the
data of smaller rank. Then every Johnson sequence for $\caL$ has injective maps. Hence Theorem \ref{conj:geom}(2) holds,
by induction on the rank.
\end{prop}

\begin{proof}
We may assume that $\pi(\caL)$ is not supercuspidal (\S\ref{sec:step1}). By induction on $k$, assume
$f_0,\ldots,f_{k-1}$ injective; $f_k$ exists by Lemma
\ref{lem:existence}. Let $K=\ker f_k$. We show that $r^{\rho,\mu}(K)=0$ for every self-dual $\rho$ and every
$\mu\in\mathbb R$, which gives $K=0$ by the criterion of Section \ref{sec:step1}. We fix $\rho$ and write
$X=X_\rho$, $r^\mu=r^{\rho,\mu}$.
\begin{itemize}
\item If $\mu\notin X$: all the $r^\mu(J)$ vanish (Lemma \ref{lem:onepiece}), so $r^\mu(Y_k)=0$.
\item If $\mu$ is unpaired with $D_{\rho|\cdot|^\mu}(\pi(\caL))\neq0$: by Lemma \ref{lem:inj_unpaired},
$r^\mu(f_k)$ is a map of a Johnson sequence of $\caL_\mu$, injective by induction on the rank.
\item If $\mu$ is unpaired with $D_{\rho|\cdot|^\mu}(\pi(\caL))=0$: all the $r^\mu(J)$ vanish
(Lemma \ref{lem:ato48}), so $r^\mu(Y_k)=0$.
\item If $\mu$ is paired: $r^\mu(\pi(\caL))=0$, and the summands not killed by $r^\mu$ are paired by Atobe's involution;
for a pair $\{J_v,J_u\}$, we have $r^\mu(J_v)\simeq r^\mu(J_u)$ (Proposition \ref{prop:cat48}) and
$\deg J_u=\deg J_v+1$ (Lemma \ref{lem:degpaired}). Write $Q_v=r^\mu(J_v)$, and $Q_u=r^\mu(J_u)=\rho|\cdot|^\mu\otimes
\widetilde J_u$. By induction on $k$, the projection $\caJ^{k-1}\to Y_k$ induces an isomorphism
$\bigoplus_vQ_v\simeq r^\mu(Y_k)$, the sum over the lower members $v$ of degree $k-1$. Let $u$ be an
upper member of degree $k$, $v$ its partner, and $\varphi_u:Y_k\to J_u$ the component of $f_k$.
\begin{itemize}
\item[--] $r^\mu(\varphi_u)\neq0$: $\varphi_u\neq0$, so its image contains the socle $L(J_u)$, and
$r^\mu(L(J_u))\neq0$ (Proposition \ref{thm:SAM}(b)).
\item[--] For a lower member $v'\neq v$ of degree $k-1$, with partner $u'$, every map
$Q_{v'}\simeq Q_{u'}\to Q_u$ vanishes. Otherwise its image would contain the socle of $Q_u$, so that
$L(\widetilde J_u)$ would be a constituent of $\widetilde J_{u'}$, and
$O_{\widetilde J_{u'}}\subset\overline{O_{\widetilde J_u}}$ (Proposition \ref{prop:closure}, applied to
these costandard modules of a smaller group). Since $d(\widetilde J_{u'})=d(\widetilde J_u)$ (Lemma
\ref{lem:degtilde}), this forces $\widetilde J_{u'}\simeq\widetilde J_u$. But $\widetilde J_u$
determines $u$ (Lemma \ref{lem:degtilde}). So $u'=u$, a contradiction.
\end{itemize}
Hence $r^\mu(\varphi_u)$ is non-zero on $Q_v$ only: the composite $J_v\to Y_k\xrightarrow{\varphi_u}
J_u$ is non-zero, and $r^\mu$ of it is an isomorphism $Q_v\simeq Q_u$ by Lemma \ref{lem:matchediso}.
(In case (ii) of Proposition \ref{prop:cat48}, the same argument applies to $J_u^+$ and $J_u^-$, with $Q_v\simeq Q_{u^+}\oplus
Q_{u^-}$.) Therefore the composite of $r^\mu(f_k)$ with the projection onto $\bigoplus_UQ_u$ (upper
members of degree $k$) is an isomorphism. So $r^\mu(f_k)$ is injective, its image is a complement of
$\bigoplus Q_{v''}$ (lower members of degree $k$), and the induction proceeds.
\end{itemize}
\end{proof}

\noindent The non-vanishing of the components $J_v\to J_u$ between the members of a pair, which was
needed above, is the analogue of the fact that all the maps of a BGG resolution are non-zero
(\cite{Hum08}, Theorem 6.8). Here it follows from the degree computations of Lemmas
\ref{lem:degpaired} and \ref{lem:degtilde}.

\subsection{Proof of the main theorem}
By Theorem \ref{thm:signs} and Proposition \ref{prop:johnson} (whose hypotheses, Theorem
\ref{conj:geom}(1) and (Inj), hold by Theorems \ref{thm:signs} and \ref{thm:injprime}, Lemmas
\ref{lem:extcs} and \ref{lem:K}, and Corollary \ref{cor:selfext}), Theorem \ref{conj:geom} holds for every
ladder datum $\caL$.

\begin{rmk}[Cousin resolutions]\label{rmk:cousin} We make the Cousin mechanism of the introduction
precise in the present setting. In a category of geometric origin, a module $M$ has a coresolution by
costandard objects with one term for each stratum, in the degree given by its codimension, as soon as,
for each standard object, the groups $\mathrm{Ext}^i$ from it to $M$ are one-dimensional in the degree
given by the codimension of its stratum, and vanish otherwise. Here the standard objects are Kato's modules $\widetilde K_\mu$ (\S\ref{ssec:inj}), and
Theorem \ref{conj:geom}(2) is equivalent to
\[\mathrm{Ext}^i(\widetilde K_\mu,\pi(\caL))\simeq\mathbb{C}\ \text{ if }\mu\in\Phi(\caL)\text{ and
}i=\deg\mu,\qquad\mathrm{Ext}^i(\widetilde K_\mu,\pi(\caL))=0\ \text{ otherwise},\]
provided the semi-orthogonality $\mathrm{Ext}^i(\widetilde K_\mu,K_\lambda^*)=\delta_{\mu\lambda}
\delta_{i0}$ holds (here, for graded modules $M,N$ over a graded algebra,
$\mathrm{Ext}^i(M,N)=\bigoplus_j\mathrm{Ext}^i(M,N\langle j\rangle)$ is computed in the category of graded
modules, $N\langle j\rangle$ denoting the shift of the grading, as in \cite{Kat17}, \S1) for the parameters $\lambda$ of discrete type. This is \cite{Kat17}, Theorem C(2),
under the assumption $(\spadesuit)'$ that the stabilizers are connected. Here they are not, but for a
parameter of discrete type the component group acts trivially on $H^\bullet_{Z_\lambda^\circ}(\mathrm{pt})$
(Lemma \ref{lem:compgroup}), so that, by \cite{Kat17}, Theorem 1.8(4) (which does not need
$(\spadesuit)'$), $\widetilde K_\lambda$ has a filtration by shifts of $K_\lambda$ alone, as in the
connected case. One direction is immediate: if the coresolution exists, apply
$\mathrm{Hom}(\widetilde K_\mu,-)$; by the semi-orthogonality the terms are acyclic and
$\mathrm{Hom}(\widetilde K_\mu,J)$ is $\mathbb{C}$ if $J$ has parameter $\mu$ and $0$ otherwise, so the
complex computing $\mathrm{Ext}^\bullet(\widetilde K_\mu,\pi(\caL))$ has at most one non-zero term, in
degree $\deg\mu$. Conversely, these vanishing conditions are what the Johnson construction of this section needs;
the proof of Theorem \ref{conj:geom} does not use this remark.
\end{rmk}

\end{document}